\documentclass[10pt]{amsart}

\newcommand{\proclaim}[2]{\medbreak {\bf #1}{\sl #2} \medbreak}

\newcommand{\ntop}[2]{\genfrac{}{}{0pt}{1}{#1}{#2}}

\let\newpf\proof \let\proof\relax 
\newenvironment{pf}{\newpf[\proofname]}{\qed\endtrivlist}

\def\be{\begin{equation}}
\def\ee{\end{equation}}

\def\td{{\underline {\tilde d}}}

\def\ba{{\begin{align}}}
\def\ea{{\end{align}}}

\def\d{{\underline d}}
\def\j{{\iota}}
\def\n{{\eta}}
\def\vs{{\varsigma}}

\def\g{{\gamma}}
\def\x{{\bf x}}

\def\0{{\mathbf 0}}

\def\cal{\mathcal}

\newcommand{\Per}{\operatorname{Per}}
\newcommand{\Fix}{\operatorname{Fix}}

\newtheorem{thm}{Theorem}[section]
\newtheorem{cor}[thm]{Corollary}

\newtheorem{lemma}[thm]{Lemma}

\newtheorem{prop}[thm]{Proposition}

\theoremstyle{remark}
\newtheorem{rem}{Remark}[section]

\newtheorem{example}{Example}[section]

\numberwithin{equation}{section}

\def \bn {\hfill \\ \smallskip\noindent}

\theoremstyle{definition}

\newtheorem{definition}{Definition}[section]
\def\proof{\bn {\bf Proof.} }

\def\note#1
{\marginpar
{\tiny $\leftarrow$
\par
\hfuzz=20pt \hbadness=9000 \hyphenpenalty=-100 \exhyphenpenalty=-100
\pretolerance=-1 \tolerance=9999 \doublehyphendemerits=-100000
\finalhyphendemerits=-100000 \baselineskip=6pt
#1}\hfuzz=1pt}

\newcommand{\ra}{\rightarrow}

\newcommand{\dist}{\operatorname{dist}}

\newcommand{\inter}{\operatorname{int}}
\renewcommand{\mod}{\operatorname{mod}}

\newcommand{\orb}{\operatorname{orb}}

\newcommand{\id}{\operatorname{id}}

\newcommand{\la}{{\lambda}}

\newcommand{\Om}{{\Omega}}

\newcommand{\AAA}{{\cal A}}
\newcommand{\BB}{{\cal B}}

\newcommand{\DD}{{\cal D}}
\newcommand{\EE}{{\cal E}}

\newcommand{\FF}{{\cal F}}

\newcommand{\JJ}{{\cal J}}
\newcommand{\HH}{{\cal H}}

\newcommand{\LL}{{\cal L}}
\newcommand{\MM}{{\cal M}}

\newcommand{\UU}{{\cal U}}
\newcommand{\VV}{{\cal V}}
\newcommand{\WW}{{\cal W}}

\newcommand{\C}{{\mathbb C}}
\newcommand{\D}{{\mathbb D}}
\newcommand{\E}{{\mathbb E}}
\newcommand{\F}{{\mathbb F}}

\newcommand{\N}{{\mathbb N}}

\newcommand{\R}{{\mathbb R}}

\newcommand{\U}{{\mathbb U}}

\newcommand{\X}{{\mathbb X}}
\newcommand{\Z}{{\mathbb Z}}

\def\B0{{\bold{0}}}

\catcode`\@=12

\def\Empty{}
\newcommand\oplabel[1]{
  \def\OpArg{#1} \ifx \OpArg\Empty {} \else
  	\label{#1}
  \fi}
		
\newcommand{\comm}[1]{}
\newcommand{\comment}[1]{}

\begin{document}

\title[Unimodal maps, periodic orbits and physical measures]
{Statistical properties of unimodal maps:
physical measures, periodic orbits and pathological laminations}

\author{Artur Avila and Carlos Gustavo Moreira}

\address{
Coll\`ege de France -- 3 Rue d'Ulm \\
75005 Paris -- France.
}
\email{avila@impa.br}

\address{
IMPA -- Estr. D. Castorina 110 \\
22460-320 Rio de Janeiro -- Brazil.
}
\email{gugu@impa.br}

\thanks{Partially supported by Faperj and CNPq, Brazil.}   

\begin{abstract}

In this work, we relate the geometry of chaotic attractors of typical
analytic unimodal maps to the behavior of the critical orbit.
Our main result is an explicit formula relating
the combinatorics of the critical orbit with the
exponents of periodic orbits.
This connection between topological
and smooth invariants is obtained through an analysis of
the physical measure.  Since the exponents of periodic orbits form a
complete set of smooth invariants in this setting, we have ``typical
geometric rigidity'' of the dynamics of such chaotic attractors.
This unexpected result
implies that the lamination structure of spaces of analytic maps
(obtained by the partition into topological conjugacy classes,
see \cite {ALM}) has an absolutely singular nature.

\end{abstract}

\setcounter{tocdepth}{1}

\date\today

\maketitle

\tableofcontents

\section{Introduction}

A {\it unimodal map}
is a smooth (at least $C^2$) map $f:I \to I$, where $I
\subset \R$ is an interval, which has a unique critical point
$c \in \inter I$ which is a maximum.
A unimodal map $f$ is said to be {\it regular} if it is hyperbolic and if
its critical point is non-degenerate and is not periodic or
preperiodic.  This definition is such that the set of
regular maps coincide with the set of unimodal maps which are structurally
stable, see \cite {K2} Theorem~B.  The class of regular maps is open
in the $C^2$ topology and dense in any smooth, and even analytic, topology.

The main examples of unimodal maps are quadratic maps $p_a(x)=a-x^2$,
$-1/4 \leq a \leq 2$.  Behind their innocent definition, the dynamics of
quadratic maps reveals an intricate structure and has been subject of
intense research in the past few decades.

\bigskip

Recently, several works have concentrated on investigating the dynamics
of typical unimodal maps.  The most natural notion of typical in this
context is measure-theoretical: a dynamical property is said to be
typical {\it in the quadratic family} if it is satisfied by $p_a$
for Lebesgue almost every parameter $a$.
This notion easily extends to the (infinite-dimensional) setting of general
unimodal maps: a property is typical if it corresponds to a full measure
set of parameters in an {\it ample class of families} of unimodal
maps\footnote{This notion of typical is inspired by Kolmogorov.}.

The dynamics of regular maps is quite well understood.  Due to the works of
Jakobson and Benedicks-Carleson, non-regular unimodal maps correspond to a
positive measure set of parameters in a large ($C^2$ open) set of
parametrized families.  In the works \cite {regular},
\cite {AM1}, the dynamics of typical non-regular
quadratic maps was described in great
detail from the {\it statistical point of view}.
Those results were subsequently
extended to typical analytic (and even smooth) unimodal maps in \cite {ALM},
\cite {AM2} (in the quasiquadratic\footnote{A $C^3$ unimodal map is
said to be quasiquadratic if any $C^3$ perturbation is conjugate
to a quadratic map.} case), and finally in all generality
in \cite {Av}: a typical non-regular unimodal
map $f$ possess a unique {\it non-trivial chaotic attractor} $A_f$,
which is a transitive
finite union of intervals (where periodic orbits are dense).
Moreover, this attractor is the
support of an absolutely continuous invariant measure $\mu_f$,
with excellent stochastic properties (due, notably, to the Collet-Eckmann
condition).  The measure-theoretical dynamics of
$f$ can be described by $\mu_f$ and finitely many {\it trivial attractors}
(hyperbolic periodic orbits).  The attractor $A_f$ of $f$ can be defined
also on topological grounds: it is simultaneously a metric and
topological attractor in the sense of Milnor (see \cite {attractors}).

\bigskip

Our aim in this paper is to establish much finer geometric properties of
the non-trivial attractor of a typical non-regular analytic unimodal map
$f$.  Roughly speaking, we will show how topological invariants of $f$
(coded using the theory of Milnor-Thurston) can
be used to determine (and actually compute)
a complete set of smooth invariants of $A_f$.

In the proof of this connection between topological and smooth invariants,
the physical measure $\mu_f$ will play an important role.  One of our most
important steps is to show how the information contained in
the physical measure is enough to compute some geometric invariants
of hyperbolic Cantor sets.

\bigskip

Our main theorem can be seen as a proof of ``geometric rigidity'' in the
typical setting, which is rather unexpected and even looks
paradoxical at first.
Fortunately, it is possible to visualize this consequence using the results
of \cite {ALM}.  The resulting rather amusing picture is related
to some recently discovered examples of
measure-theoretical pathological laminations (Katok's ``Fubini Foiled''
phenomena presented by Milnor \cite {Mi},
and the examples in \cite {SW} and \cite {RW}).

\comm{
Namely, the results of \cite {ALM} show that the partition of spaces of
analytic unimodal maps in topological classes has the structure of a
codimension-one lamination, with analytic leaves and quasisymmetric
transverse structure.  In particular, the parameter space of the quadratic
family has a universal quasisymmetric structure.
This metric universality result was used in \cite {ALM} and \cite {AM2}
to transfer certain results from the quadratic family to other
analytic families of (quasiquadratic) unimodal maps.

More recently, the results of \cite {ALM} were used in a different way in
\cite {Av} to produce tools that allow to work directly in analytic
families without first transferring to the quadratic setting.  The main
consequence is that the estimates thus obtained are much sharper, and are
exactly the ones used in the present work.  Indeed, our results can not
be obtained by metric universality considerations.  On the contrary, they
show the limitations of this method:
{\it the lamination by topological classes
fails to be absolutely continuous in a drastic way}.  This singular
measure-theoretical behavior is comparable to Katok's ``Fubini Foiled''
counterexample as presented by Milnor \cite {Mi}, see also the examples of
\cite {SW} and \cite {RW}.
}

\subsection{Statement of the results}

In this work, the ample set of families we will consider for the definition
of typical is very explicit: the set of non-trivial analytic families of
unimodal maps,
that is, families which contain a dense set of regular parameters.  The set
of non-trivial families is very large (its complement has infinite
codimension).  Moreover, among families of quasiquadratic maps (a $C^3$ open
condition) it is much
easier to check for non-triviality: it is enough to show existence of one
regular parameter (which is a $C^2$ open condition).  In particular,
analytic families $C^3$ close to the quadratic family are non-trivial.

\subsubsection{Typical unimodal maps and their invariants: relation between
topological and combinatorial invariants}

To each point $x \in I$, let us associate an infinite
sequence (the {\it itinerary})
of $0$s and $1$s as follows.  The $k$-th element is $0$ if
$f^k(x)$ is to the left of the critical point, and $1$ otherwise.
Itineraries are clearly invariant under topological conjugacy.
The itinerary of the critical point of $f$ is called the {\it kneading
sequence} of $f$, and it is a particularly important invariant:
the work of Milnor-Thurston shows that the kneading sequence
determines the set of itineraries of all points $x \in I$. 

The kneading sequence is actually an ``essentially'' complete topological
invariant in the sense that it determines the topological conjugacy
class up to some well understood obstructions corresponding to
trivial dynamics.  A simpler (and perhaps more basic, as it
applies in all dimensions)
example of topological invariant is the
set of periodic orbits of the system, together with their periods.
If $p$ is a periodic point, its itinerary is clearly periodic.

To a periodic orbit $p$ of period $n$ we can associate its exponent
$Df^n(p)$.  This quantity is easily checked to be invariant by a
diffeomorphic change of coordinates, thus providing the simplest example
of a smooth invariant.  By the work of Li\v vsic \cite {Li},
see also Shub-Sullivan \cite
{ShSu}, in some circumstances (say, expanding maps of the circle)
exponents of periodic orbits form a complete set of
smooth invariants, in the sense that a topological conjugacy which preserves
exponents is necessarily smooth.  In the unimodal case, the same result
holds due to the work of Martens-de Melo \cite {MM},
at least for the cases that appear in our considerations
(non-trivial attractor of a typical non-regular
unimodal map).

The main result of this paper relates the above combinatorial and smooth
invariants for typical non-regular analytic unimodal maps.

\proclaim{Theorem A.}
{
Let $f_\lambda$ be a non-trivial analytic family of unimodal maps. 
Then, for almost every non-regular parameter $\lambda$, and for every
periodic orbit $p$ in the non-trivial attractor $A_{f_\lambda}$,
the exponent of $p$ is determined by an
explicit combinatorial formula involving the kneading sequence of $f$ and
the itinerary of $p$.
}

The formula goes as follows: let $\beta$ be the kneading sequence of $f$
and let $\alpha$ be the periodic part of the itinerary of a periodic
point $p$ in $A_f$.  Let us consider the asymptotic
frequency $r(\alpha^k,\beta)$ of $\alpha^k$ ($k$ repetitions of $\alpha$)
inside $\beta$.  Ignoring for a moment the problem of existence of this
asymptotic frequency (which is part of Theorem B below), we obtain a
non-increasing sequence of numbers between $0$ and $1$.  It turns out that
this sequence decreases to $0$ geometrically at some precise rate (this is
related to Theorem C below).  The inverse of this rate is the absolute
value of the exponent of $p$ (the sign being given by $(-1)^s$ where $s$ is
the number of $1$s in $\alpha$).

\comm{
\subsubsection{Preliminaries}

Before stating our results precisely, we must first clarify more what
we understand for a typical analytic unimodal map and introduce some known
facts about the dynamics of those maps.

A unimodal map $f$ is said to be regular if it has all periodic orbits
hyperbolic, its critical point is non-degenerate and is not periodic or
preperiodic.  This definition is so that regular maps coincide with
structurally stable maps (by \cite {K}).  The class of regular maps is open
in the $C^2$ topology and dense in any smooth, and even analytic, topology.

The ample class of analytic families of
unimodal maps which we consider for the definition of
typical is the set of non-trivial families: families with a dense set of
regular parameters.  This ample class of families is readily seen to be
dense (in the analytic topology), and its complement has infinite
codimension.

In \cite {Av} it was proved that a typical non-regular
analytic unimodal maps can be renormalized to a quasiquadratic map,
that is, are $C^3$-persistently conjugate to a quadratic map.

To simplify our exposition, we will
describe our results first in the case of quasiquadratic maps
(which includes maps with negative Schwarzian derivative), which is a $C^3$
open class of unimodal maps.
The ample class of analytic families of unimodal maps
which we consider for our definition of typical is the class of non-trivial
analytic families, which can be characterized as having one regular
parameter.  This ample class is readily seem to be open (in the $C^3$
topology) and dense (in the analytic topology), and is of infinite
codimension as we will see later.  It is readily seem that the quadratic
family is non-trivial.

We give briefly some properties of typical quasiquadratic maps
obtained by \cite {ALM} and \cite {AM2}, which are needed to describe our
results.

There are two rather different kinds of behavior for typical quasiquadratic
maps.  One of them is regular behavior (since it is $C^2$ robust).  Regular
maps enjoy the many nice statistical properties of
hyperbolic dynamical systems, but are essentially understood.  The behavior
of typical non-regular maps is much more subtle.

If $f$ is a typical non-regular unimodal map, it is still possible to
obtain a good description of its dynamics from the statistical point of
view.  Those statistical properties are roughly consequence of two
estimates on the behavior of the critical orbit.  The first is a positive
Lyapunov exponent for the critical value, the Collet-Eckmann condition.  The
second is the rate of recurrence of the critical point, which is polynomial.

In the statistical description of a typical unimodal map, we make use of
physical measures, that is, an ergodic invariant probability measure
$\mu$ such that for a positive measure set of $x \in I$, and any continuous
$\phi:I \to \R$ we have
$$
\lim \frac {\sum_{k=0}^{m-1} \phi(f^k(0))} {m}=\int \phi d\mu.
$$
The set of such $x$ is called the basin of $\mu$.  For typical
quasiquadratic maps, we have existence and uniqueness of the physical
measure, which moreover has a basin of full measure.  While in the regular
case $\mu$ is concentrated in the periodic attractor, the non-regular case
is more interesting and $\mu$ is absolutely continuous.
In this case, the support of $\mu$ is
still easy to describe: it is a finite union of intervals
$T_0$, $T_1$,..., $T_{n-1}$, where $T_0$ is the smallest periodic interval
of $f$ containing the critical point and $n$ is its period.  Another
description is that $n$ is the biggest period of a renormalization of $f$
and $T_0=[f^{2n}(0),f^n(0)]$.  Moreover, the interior of the $T_i$ are
disjoint, and $f^n(T_0)=T_0$.

The support of $\mu$ is also the metric and topological attractor of $f$ in
the sense of Milnor (for a discussion, see \cite {attractors}),
and we will denote it by $A$ and call it just the
attractor of $f$.

\subsubsection{Formula for exponents of periodic orbits}

It is natural to try to understand geometric aspects of the attractor $A$ of
a typical unimodal map.  Those are related to the smooth conjugacy class of
$f|A$.  One of the most important smooth invariants in one-dimensional
dynamics is the exponent of periodic orbits: given a periodic orbit $p$ of
period $n$, $Df^n(p)$ is invariant by smooth coordinate change.

Our main result is that, for a typical non-regular analytic maps $f$, the
exponents of all periodic orbits $p$ contained in the attractor $A$ of $f$
are determined by some explicit formulas which only depend on the
combinatorial aspects of $p$ and the critical point.

More precisely,
Milnor-Thurston associate to $f$ some symbolic dynamics.  This
combinatorial model correspond to each point $x$ in $I$ some
kneading sequence $\theta(x)$.  We are able to define a combinatorial
notion of geometric frequency of (repetitions of)
one periodic kneading sequence into another which we show to correspond to
exponents of periodic orbits:

\proclaim{Theorem B.}
{
Let $f_\lambda$ be a non-trivial analytic family of quasiquadratic maps. 
Then, for almost every non-regular parameter $\lambda$, and for every
periodic orbit $p$ in $A_f$, the exponent of $p$ is determined by an
explicit combinatorial formula
$$
|Df^n(p)|=e^{-\rho}
$$
where $\rho$ denotes the geometric frequency of the (periodic part) of the
kneading sequence of $p$ in the kneading sequence of $0$.
}

Notice that while
exponents of periodic orbits are only smooth invariants, the combinatorics
of unimodal maps are topological invariants, so this result points to
an unexpected rigidity result which will elaborated later.

The proof of Theorem A involves two other Theorems of independent interest. 
}

\subsubsection{The critical orbit is typical}

Let us say that the asymptotic distribution of a point $x$ is given by
$\mu$ (or equivalently, $x$ is in the basin of $\mu$, or $x$ is typical for
$\mu$) if $\mu$ is a probability measure and for any continuous
function $\phi:I \to \R$
\be
\lim \frac {1} {n} \sum_{k=0}^{n-1} \phi(f^k(x))=\int \phi d\mu.
\ee

One important step of the proof of Theorem A is to analyze the asymptotic
distribution of the critical orbit.  The existence of an asymptotic limit
for the distribution of the critical orbit is
directly related to the existence of asymptotic
frequencies $r(\alpha,\beta)$
of an arbitrary finite sequence $\alpha$ inside the kneading sequence
$\beta$ of $f$.


\proclaim{Theorem B.}
{
Let $f_\lambda$ be a non-trivial analytic family of quasiquadratic maps. 
Then, for almost every non-regular parameter $\lambda$, the critical point
belongs to the basin of $\mu_{f_\lambda}$
(the absolutely continuous invariant measure
of $f_\lambda$).
}

In other words, for typical non-regular unimodal maps, the critical orbit is
typical for the ``correct'' measure of the system.
We are thus able to obtain the following consequence:

\begin{cor} \label {lyap}

In the setting of Theorem B, one also has equality between the Lyapunov
exponent of the critical value and the Lyapunov exponent of
$\mu_{f_\lambda}$.

\end{cor}

Recall that the Lyapunov exponent of a point $x$ is defined as
\be
\lambda(x)=\lim \frac {\ln |Df^n(x)|} {n}
\ee
provided the limit exists.
The Lyapunov exponent of $\mu_f$ is given by the formula
\be
\lambda(\mu_f)=\int \ln |Df| d\mu_f.
\ee
Some work is needed to go from Theorem B to Corollary \ref {lyap},
since $\ln |Df|$ is not continuous.

\comm{
The relevance of Theorem B in the proof of Theorem A is that it allows to
study the distribution of the critical orbit (which is combinatorially
determined by the kneading sequence of $0$) through the much more manageable
ergodic properties of $\mu$.
}

Previous progress in the direction of Theorem B was achieved (with very
different techniques) by Benedicks-Carleson \cite {BC}, who proved
typicality of the critical orbit
for a {\it positive measure set} of parameters
for the quadratic family.

\subsubsection{Regularity of the physical measure and hyperbolic sets}

In Theorem A we are interested in the exponents of (repelling) periodic
orbits.  More generally, one is led to ask about the geometry of hyperbolic
subsets $K \subset I$ (say, a Cantor set).

In order to apply Theorem B to reconstruct the geometry of $K$ from the
kneading sequence of $f$, one is led to ask: is it possible
to obtain sharp estimates for the asymptotic geometry of $K$ from
knowledge of the physical measure?

In order to do so, one should be able to relate asymptotically the physical
measure of gaps (and unions of gaps) of $K$ and their (Lebesgue) size.
Thus, behind this problem is the issue of regularity of the physical measure
$\mu_f$.

It turns out that this problem is non-trivial: indeed, if one tries to
estimate general intervals, and not just gaps of hyperbolic sets, one would
get quite negative results.  For instance, let us take $f$ to be a quadratic
map and let $T$ be an interval of radius $\epsilon$ around the critical
point.  Then
$\mu_f(T)=\mu_f(f(T))$, but $|T|$ is of order $\epsilon$ while $|f(T)|$ is
of order $\epsilon^2$.  Thus, for general intervals, estimates of
the physical measure might lead to errors of order $2$ (when taking
logarithms) on estimates of
Lebesgue measure (and thus on the formula for exponents of periodic orbits).
Connected to this fact is the following limitation on the
regularity of $\mu_f$: its density $d\mu_f$ is never in $L^2$.

So one is led to regularize the density $d\mu_f$ using the Cantor set $K$
(or view $d\mu_f$ through $K$).  Let us denote
$d\mu_f^K$ the function which is constant in each gap $T$
of $K$ and takes the average value of $d\mu_f$ on $T$.

In other words, $d\mu_f^K$ is the expectation of $d\mu_f$ with respect
to the sigma-algebra $\BB(K)$ of the gaps of $K$.
The sigma algebra $\BB(K)$ gives us enough information to
compute the exponent of periodic orbits $p$ in $A_f$ if, say, $K$ is a
Cantor set containing $p$ (any periodic point $p \in A_f$
can be included in such a Cantor set).

\comm{
The other component in the proof of Theorem A is a non-trivial estimate on
the regularity of $\mu$.  This is necessary to relate the measure of a set
with its (Lebesgue) size.  The fact that $\mu$ is absolutely continuous is
not quite enough: it is known that $d\mu$ does not belong to $L^2(I)$.  This
limitation of regularity can produce an error of up to order $2$
in the formula for the exponents of periodic orbits.

It turns out that to prove Theorem A, we only need to estimate its
regularity with respect to certain dynamically relevant sets.  A natural
choice is to consider hyperbolic invariant sets $K$: it is quite easy to
include any periodic orbit in (quite thick) such sets.  Roughly
speaking, we will prove a general estimate relating the physical
measure $\mu$ and the geometry of the hyperbolic sets $K$.  Reciprocally,
one can see this estimate (together with Theorem B)
as a generalization of Theorem A, since it allows to
compute fine geometric asymptotic properties of general
hyperbolic sets (of which periodic orbits are an example)
using $\mu$ (which, due to Theorem B can be computed
combinatorially).

Let us explain a little bit more the concept of regularity of a measure with
respect to a hyperbolic set (say, a Cantor set $K$).

In other words, we consider a sigma-algebra $\BB(K)$
which is much smaller then the Borelian
sigma-algebra and estimate the regularity of $E(d\mu|\BB(K))$
(the expectation of $d\mu$ with respect to $\BB(K)$).

We must choose $\BB(K)$ in such a way that $E(d\mu|\BB(K))$ will encode
relevant information on the geometric aspects of $K$.  The natural choice
for $\BB(K)$ is the sigma-algebra generated by the gaps of $K$.  In this
case, $E(d\mu|\BB(K))$ has a particularly simple expression: its value
restricted to a gap $\Lambda$ of $K$ is just the average value of
$d\mu$ in $\Lambda$ equal to $\mu(\Lambda)/|\Lambda|$.
Since it is obtained by an averaging
process, it is natural to think of $E(d\mu|\BB(K))$ as a regularization of
$\mu$ (obtained by viewing it through $K$).
}

\proclaim{Theorem C.}
{
Let $f_\lambda$ be a non-trivial analytic family of unimodal maps. 
For almost every non-regular parameter $\lambda$ and any hyperbolic
set $K \subset I$, we have
$d\mu_{f_\lambda}^K \in L^p$, $1 \leq p <\infty$.
}

One can see this estimate (together with Theorem B)
as a generalization of Theorem A, since it allows to
compute using $\mu_f$ (which, due to Theorem B can be computed
combinatorially), fine asymptotics of general
hyperbolic sets (of which periodic orbits are an example).

We should point out that the lack of regularity of $\mu_f$ comes from
the critical point, and essentially distributes itself along the orbit of
the critical value.  In order to show that $\mu_f$ behaves well with
respect to hyperbolic sets, we must show roughly that
``the critical orbit distributes transversely with respect to $K$''.

\subsubsection{Geometric rigidity, pathological laminations}

The main motivation for Theorem A is, as described before, the possibility
to compute, from topological information, a complete set of smooth
invariants.  This may seem at first paradoxical, since
exponents of periodic orbits {\it can} be varied without
changing the topological class, and are thus responsible for moduli of
flexibility, as opposed to rigidity (examples of geometrically rigid
systems, as diophantine irrational rotations usually do not
have periodic orbits).

In order to visualize what is really happening, we must consider the
partition of the space of unimodal maps into topological conjugacy classes.
The results of \cite {ALM} show that, in appropriate Banach spaces of
analytic unimodal maps, the set of
non-regular topological classes form a lamination
with analytic leaves and quasisymmetric holonomy, at least
almost everywhere\footnote{Almost
everywhere here is indeed stronger than our notion of typical.  More
precisely, the set of non-regular topological classes has a     
lamination structure in an open set containing all Kupka-Smale maps
(unimodal maps with a non-degenerate critical point and without
non-hyperbolic periodic orbits).  The                           
complement of this open set is clearly contained in a countably union of
codimension-one analytic varieties.}.

For each topological class of unimodal maps, the formula for
exponents of periodic orbits determines {\it at most one} ``preferred''
smooth structure on the non-trivial attractor\footnote {For a
general topological class, several things might go wrong, so that no smooth
structure is determined.
At the level of the formula,
for instance, its defining limits might not exist.  The non-trivial
attractor may not exist.  Even if both exist, the
values for exponents thus
obtained might not correspond to any smooth structure on the
non-trivial attractor.}.  In each non-regular
topological class (of codimension one by \cite {ALM}),
the set of maps with the ``correct'' smooth structure is
a tiny set (infinite codimensional, the parameters being precisely the
exponents of periodic orbits, or even empty).
However, the set of typical non-regular unimodal
maps (satisfying the conclusion of Theorem A) intersects each
topological class precisely at such a tiny set.

So ``typical rigidity'' has interesting
consequences for the regularity of the lamination by
topological classes: {\it the
stratification of the set of typical non-regular analytic unimodal maps by
topological classes is highly non-homogeneous, in the sense that it fails
drastically to be absolutely continuous}.  Indeed, that the lamination can
not be absolutely continuous is easily checked since the phenomena we
described imply the complete failure of Fubini's Theorem.  (Although the
setting is infinite dimensional, one can interpret those results in
parametrized families with at least two parameters.)

\subsubsection{On universality and the holonomy method} \label {universa}

The results of \cite {ALM} imply
that the parameter space of the quadratic
family do have a universal quasisymmetric structure
(due to the holonomy of the lamination).  Although quasisymmetric maps are
not necessarily absolutely continuous,
the metric universality was used in \cite {ALM} and \cite {AM2}
to transfer certain strong measure-theoretical
results (regular or stochastic dichotomy, Collet-Eckmann condition
and polynomial recurrence of the critical orbit)
from the quadratic family to other
analytic families of (quasiquadratic) unimodal maps.

This so called
{\it holonomy method}, consisting in the comparison between
parameter spaces of different families had to be applied to estimates which
are topological invariants.  More seriously, the set of
combinatorics concerned must have full measure {\it simultaneously} in all
non-trivial families of unimodal maps.

The lack of absolute continuity of the lamination established now
sets a limit to the metric universality
of the parameter space of unimodal families (as the quadratic family).  Our
Theorem A is particularly interesting in this respect since it gives an
example of a result which is definitely inaccessible by
the holonomy method (which clearly can not be used to prove that the
lamination itself is not absolutely continuous).

\comm{
However, the results of \cite {ALM} are still at the
base of our estimates in this paper.

We should point out that the Phase-Parameter relation of \cite {Av} used in
this work is still based on the lamination structure.  However, one uses
local holonomy maps to obtain relate phase and parameter of the same family,
instead of using the global holonomy to relate phase and parameter between
some non-trivial family and the quadratic family, which introduces serious
distortion and lack of sharpness in the estimates.

In \S \ref {universa}, we will comment
(after the necessary technical
preparation) on this subject further.
}

\subsubsection{Related matters}

Another consequence of our techniques is existence of a combinatorial
formula for the Lyapunov exponent of typical non-regular unimodal maps. 
This exponent coincides with the one of the critical value by Corollary \ref
{lyap}.  This formula is quite simple, but is formulated in terms of
the principal nest description of the combinatorics instead of
itineraries, so we postpone its formulation to \S \ref {lyapformula}.

\bigskip

In view of Theorem A, it is natural to ask how to effectively relate the
information about the exponents of periodic orbits to other properties of
interest of a typical non-regular unimodal map.  Although we will not
investigate this problem in this paper, we would like to call attention to
one situation where such a relation might be explicitly obtained.

It is common to organize periodic orbits in a {\it zeta function}.
The general formula for a zeta function is
\be
\zeta_\phi(z)=
\exp \left (\sum_{n=1}^\infty \frac {z^n} {n} \sum_{p \in \Fix(f^n)}
\prod_{k=0}^{n-1} \phi(f^k(p)) \right )
\ee
where $\Fix(f^n)$ is the set of fixed {\it points} of $f^n$ and
$\phi$ is a weight function which is to be chosen according to the
problem to be studied.

The relation of zeta functions and the thermodynamical formalism of
hyperbolic dynamical systems is well developed.
However it is reasonable to expect that this
relation might also hold for certain non-uniformly hyperbolic unimodal maps,
and in \cite {KN} some results in this direction were obtained in the
Collet-Eckmann case.

For the weight $\phi=|Df|^{-1}$, the zeta function can be written as
\be
\zeta_{|Df|^{-1}} (z)=
\exp \left (\sum_{n=1}^\infty \sum_{m=1}^\infty
\frac {z^{mn}} {m} \sum_{p \in \Per_n(f)}
\frac {1} {|Df^n(p)|^m} \right )
\ee
where $\Per_n(f)$ is the set of periodic {\it orbits} of (prime) period $n$.
Notice that in this case the zeta function only depends on the exponent of
periodic orbits, so by Theorem A it can be expressed combinatorially for
typical non-regular maps.  This choice of the weight is particularly
interesting: it is related to the physical measure $\mu_f$, and the results
of \cite {KN} show that the poles of $\zeta_{|Df|^{-1}}$ can be sometimes
related to parts of the spectrum of the Ruelle transfer operator,
which encodes (in some cases precise) information
about the rates of decay of correlations of the system
(for certain classes of observables).
It is a natural problem to show that the first pole of $\zeta_{|Df|^{-1}}$
gives indeed the exact rate of decay of correlations (for smooth enough
observables) of typical non-renormalizable unimodal maps.


\subsection{Complex techniques}

The succesful investigation of families of unimodal maps, specially
the quadratic family, was heavily tied to the
possibility of the intertwined use of real and complex techniques.
Our results are based on the coupling of two main methods.
For the analysis of the dynamics in phase space, we use
a statistical description of the critical orbit.  Techniques from complex
dynamics are used to obtain the Phase-Parameter relation, which allows to
compare the phase space and the parameter space of a
non-trivial family.  Those complex techniques are mainly based in the
theory of Lyubich (which in turn uses ideas from several different fields).

The Phase-Parameter relation was proved in \cite {AM1} in the case of the
quadratic family, and in \cite {Av} in all generality.
This last result can be
directly used in our context and will allow us to concentrate mostly on the
real dynamics of unimodal maps.

\begin{rem}

Although, as discussed in \ref {universa},
the holonomy method of \cite {ALM} can not be used
for this work, the lamination structure of the partition into topological
classes is still the key result from complex dynamics used in
the proof of the Phase-Parameter relation in \cite {Av}.
This is possible because the regularity of the holonomy map between two
transversals is related to their distance.
The original holonomy method has a global nature
and corresponds to relating phase and parameter between a given 
non-trivial family and the quadratic family, introducing serious 
distortion and lack of sharpness in the estimates.
To get over those limitations, one uses local holonomy maps (which are more
regular) to relate phase and parameter of the same family.

\end{rem}

\comm{
This last development makes
use of the lamination by topological classes of \cite {ALM} which we
mentioned before.

Although this paper will concentrate on the statistical
analysis, it is worth to understand a couple of the ideas involved:
this will allow to understand better the significance of
Theorem A in the study of universal properties of the
quadratic family.  We will comment on this in \S \ref {universa}, after
presenting formally our results on the pathological behavior of the
lamination by topological classes.

For the case of the quadratic family, the estimates of \cite {AM1} are
enough to conclude (together with the statistical analysis of this
paper) Theorems A, B and C.  To extend the result to any non-trivial
analytic family, we use the results
of \cite {Av}.  Although this paper will concentrate on the statistical
analysis, it is worth to present a couple of the ideas involved: this will
allow to understand better the significance of Theorem A in the study of
universal properties of the quadratic family.

\subsubsection{Lamination structure in spaces of analytic unimodal maps}

The key result to obtain measure-theoretical for general analytic families
of unimodal maps is the following structure theorem of \cite {ALM}.  Each
non-regular topological conjugacy class is a codimension-one analytic
submanifold in spaces of analytic unimodal maps.  Moreover, different
topological conjugacy
classes fit together to form, almost everywhere (indeed, outside of
parabolic leaves), a lamination.  The holonomy of this lamination is
quasisymmetric (and so the parameter space of the quadratic family
has a universal quasisymmetric structure).  This is used by \cite {ALM} as a
shortcut to extend the regular or stochastic dichotomy.
In \cite {AM2}, this shortcut is used to
deduce the Collet-Eckmann and polynomial recurrence conditions from
statements for the quadratic family.

Quasisymmetric maps are not
absolutely continuous, so this can be only applied to transfer results which
are robust in the sense that the corresponding set of combinatorics have
full measure simultaneously in all non-trivial analytic families. 

While this method allows to obtain strong results, it does not give access
to fine geometric information.  Roughly, instead of comparing the
phase and parameter space of a given family, we compare the phase space of
the family with the parameter space of the quadratic family.  This gives
less than optimal estimates for the phase-parameter relation.

To understand why the holonomy method is not appropriate in our setting, we
notice that one of our motivations to prove Theorem A is to obtain a
rigidity statement which translates in complete lack of
absolute continuity of the lamination.  It is clear that we can
not use such a holonomy method to prove
lack of regularity of the holonomy itself.
}

\comm{
The solution is to investigate the phase-parameter relation at any
non-trivial family.  In order to do so, we must get into the proof of the
lamination structure of \cite {ALM}, particularly the estimates related to
the existence of a transverse direction to the lamination.  Those estimates
will be used to create, near any family transverse to the lamination,
a special family which has (locally) similar properties to the
quadratic family itself.  Information is transferred from this special
family through a local holonomy map, which has better estimates then the
global holonomy map that relates to the quadratic family.

As a consequence of our improved estimates, we deduce an improvement on
\cite {AM2}, which is a direct consequence of the phase-parameter relation
and the statistical analysis of \cite {AM1}

\begin{cor}

Let $f_\lambda$ be a non-trivial analytic family of quasiquadratic maps. 
Then, for almost every non-regular parameter $\lambda$, the recurrence of
the critical point is polynomial with exponent $1$.

\end{cor}

In \cite {AM2}, the exponent depends on the map (essentially its distance
from the quadratic family).
}

\comm{
\subsection{Motivation and consequences}

The main motivation to study the exponents of periodic orbits is the
following result of Shub and Sullivan.  Let $f$ and $g$ be two smooth
uniformly expanding maps of the circle of the same degree and let $h$ be a
topological conjugacy between $f$ and $g$.  If $h$ preserves exponents
of periodic orbits, that is, for any periodic orbit $p$ of $f$, $p$ and
$h(p)$ have the same exponent, then $h$ is actually a smooth map.

In other words, the exponents of periodic orbits exhausts all smooth
invariants.  For the unimodal case, \cite {MM} obtained a corresponding
result.  Thus we have the following measure-theoretical rigidity statement.

The main motivation for Theorem A is, as described before, the possibility
to compute, from topological information, a complete set of smooth
invariants.  This geometric rigidity can be stated as follows:

\begin{prop}

Let $X$ denote the set of all non-regular
analytic unimodal maps which satisfy the conclusions of Theorems A, B and C
(in particular, $X$ intersects each non-trivial family of unimodal maps
in a full measure set in the complement of hyperbolic maps).
If $h:I \to I$ is a homeomorphism such that $h \circ f=g \circ h$ then
$h|A_f$.

\end{prop}

Notice that $h$ may fail to be globally analytic, see \S.

This allows us to obtain the following Corollary regarding the lamination of
\cite {ALM}.

\begin{prop}

The set $X$ intersects each leaf of the lamination by topological conjugacy
classes (a codimension-one submanifold) in a set of infinite codimension
(analytic conjugacy class on the non-trivial attractor).

\end{prop}

In particular, the holonomy maps of the lamination take in general a full
measure set of the complement of hyperbolic maps into zero Lebesgue measure
sets.  A corresponding statement for finitely parametrized families is made
in \S.

\subsubsection{Zeta functions}

Periodic orbits also show up in zeta functions, which can be related to the
thermodynamical formalism of dynamical systems.  Let $\Fix(f^n)$ be the set
of fixed points of $f^n$ and $\Per_n(f)$ be the set of periodic orbits of
(prime) period $n$.  The generic formula for a
zeta function is
\be
\zeta_\phi(z)=
\exp -\sum_{n=1}^\infty \frac {z^n} {n} \sum_{p \in \Fix(f^n)}
\prod_{k=0}^{n-1} \phi(f^k(p))
\ee
where $\phi$ is a weight function which is to be chosen according to the
problem at present.  In particular, the weight $|Df|^{-1}$ is related to the
physical measure $\mu_f$.  In this case, the zeta function can be written as
\be
\zeta_{|Df|^{-1}} (z)=
\exp -\sum_{n=1}^\infty \sum_{m=1}^\infty
\frac {z^{mn}} {mn} \sum_{p \in Per_n(f)}
\frac {1} {|Df^n(p)|^m}
\ee
which only depends on the exponent of periodic orbits, so it can be
expressed combinatorially for typical non-regular maps.  The poles of
$\zeta$ can be sometimes related to parts of the spectrum of the Ruelle
transfer operator.  This spectral information encodes for instance the
decay of correlations of the system (or rather, its last renormalization),
which is always exponential under
the Collet-Eckmann condition (but is not always in reach of the zeta
function).  It would be interesting to investigate this connection further.

\subsubsection{Formula for the exponent of $\mu$}

Another consequence of Theorems A and B
is the existence of a combinatorial formula
for the Lyapunov exponent of $\mu$, which is the same as the Lyapunov
exponent of the critical value.
}

\subsection{Outline}

In \S \ref {prelim}, we present some background on
the dynamics of unimodal maps.  In \S \ref {formula},
we state precisely the formula for periodic orbits.
We then prove Theorem A, assuming the validity of Theorems B and C.

In \S \ref {pp}, we discuss the combinatorics of the principal nest and
introduce our basic tool to make parameter estimates:
the Phase-Parameter relation, which was proved in \cite {AM1} (for the
quadratic family) and in \cite {Av}
in all generality.  We then present
some of the estimates obtained in \cite {AM1}.

\comm{
The proof of the
phase-parameter relation in the general case is very involved, and we
postpone the most technical part (which involve not only the main result of
\cite {ALM}, but its techniques) to \S.
}

In \S \ref {typical}, we prove Theorem B.  The proof is
technical but has a clear strategy, which we describe in \S \ref {thmB}.
In \S \ref {hyperbolic}, we reduce Theorem C to the so called Main estimate,
which we prove in \S \ref {est}.  The proof of the Main estimate is the
most technically involved part of this work.


\comm{
In Appendix A, we revisit M. Lyubich's theory of the quadratic family in
a slightly more general context.  We concentrate in describing precise
statements, but we do not attempt to reproduce Lyubich's proof in its
entirety in this new context (\cite {puzzle} is already quite long).  We
give a quick outline of the basic ideas on how to put the main estimates
together.

In Appendix B, we complete the proof of the phase-parameter relation.
The main
technique is a careful estimate of the ``explicit transverse'' to the
lamination obtained in \cite {ALM}.
}

\section{Preliminaries} \label {prelim}

\subsection{Notation}

As usual,
$\N=\{0,1,2,\dots\}$ stands for the set of natural numbers;
$\R$ stands for the real line;
$\C$ stands for the complex plane.

\comm{
Let $\D_r(x)= \{z \in \C:\ |z-x|<r\}$, $\D_r=\D_r(0)$, and let $\D=\D_1$.

We will reserve notation $I$ for the  interval  $[-1, 1]$. 
For  $a>0$ let $$ \Om_a= \{ z\in \C:\ \dist(z, I)< a \}.$$
A {\it topological disk} is a simply connected domain in $\C$;
a {\it Jordan disk} is a topological disk bounded by a Jordan curve.
}

The Lebesgue measure of a set $X\subset \R$ will be denoted by $|X|$.


Given a diffeomorphism $\phi: J\ra J'$ between two real intervals, its {\it
distortion}
or {\it non-linearity} is defined as 
\be
\dist(\phi)=\sup_{x,y\in J} \frac {|D\phi (x)|} {|D\phi(y)|}.
\ee
Its {\it Schwarzian derivative} is given by the formula:
\be
S\phi = \frac {D^3\phi}{D\phi} - \frac 32 \left(\frac
{D^2\phi}{D\phi}\right)^2.  
\ee 
The condition of negative Schwarzian derivative plays an important role in
one-dimensional
dynamics. This condition is preserved under composition.

\comm{
If $r>1$, let $A_r=\{z \in \C|1<|z|<r\}$.  A ring $A$ is a subset of
$\C$ such that there exists a conformal map from
$A$ to some $A_r$.  In this case, $r$ is uniquely defined
and we denote the moduli of $A$ as $\mod (A)=\ln(r)$.

Let $U\subset \C$ be a bounded open set. We say that a holomorphic
function $f \colon U \ra \C$ belongs to class $A^1(U)$ if $f$ and
its derivative $f'$ admit a continuous extension to the closure
$\overline U$. We will use the same notations
$f$ and $f'$ for the extensions.
We supply $A^1(U)$ with the seminorm 
\begin{equation}\label{A norm}
   \| f \|_1 = \max_{z\in \overline U} | f'(z) |.
\end{equation}
If $f \in A^1(U)$, $f|\overline U$ is a homeomorphism
onto its image and $f'$
does not vanish on $\overline U$, we say that
$f| \overline U$ is a diffeomorphism
(onto the image).
}

\comm{
\subsection{Quasiconformal maps} \label {quasiconformal maps}

Let $U \subset \C$ be a domain.  A map $h:U \to \C$ is
$K$-quasiconformal ($K$-qc) if it is a homeomorphism onto its image
and for any ring $A \subset U$,
$\mod (A)/K \leq \mod (h(A)) \leq K \mod (A)$.
The minimum such $K$ is called the dilatation of $h$.

Let $h:X \to \C$ be a homeomorphism and let $C,\epsilon>0$.  An extension
$H:U \to \C$ of $h$ to a Jordan disk $U$ is $(C,\epsilon)$-qc
if there exists a ring $A \subset U$ with $\mod (A)>C$
such that $X$ is contained in the bounded component of the complement of
$A$ and $H$ is $1+\epsilon$-qc.

}

\subsection{Quasisymmetric maps}

A quasisymmetric map is a homeomorphism $h:\R \to \R$ such that there exists
a constant $k$ such that for any $x \in \R$, $a>0$,
\be \label {k}
\frac {1} {k}<\frac {h(x+a)-h(x)} {h(x)-h(x-a)}<k.
\ee

Equivalently, $h$ is quasisymmetric if it has a symmetric quasiconformal
extension to the whole $\C$ (Ahlfors-Beurling).  We say that $h$ is $\g$-qs
if there exists such an extension with dilatation bounded by $\g$.
The quasisymmetric constant of a quasisymmetric map $h$ is the infimum
of the dilatations of all those extensions\footnote {It is
possible to work out upper bounds for the quasisymmetric
constant in terms of the $k$ in (\ref {k}) and inversely.}.
In particular, if $h_1$ is $\g_1$-qs and $h_2$ is $\g_2$-qs,
$h_2 \circ h_1$ is $\g_1\g_2$-qs.

If $h:X \to \R$ is a monotonic map defined on $X \subset \R$, we will also
say that $h$ is $\g$-qs if it has a $\g$-qs extension to $\R$.

\comm{
\begin{lemma} \label {qc to qs}

For all $\g>1$ there exists $C,\epsilon>0$ such that if
$h:X \to \R$ has a
$(C,\epsilon)$-qc extension $H:U \to \C$ which is symmetric then
$h$ is $\g$-qs.

\end{lemma}
}

One of the main concepts we will need in our paper was introduced in \cite
{AM1}.  The $\g$-qs capacity of a set $X \subset \R$ inside some interval $T
\subset \R$ is defined as
\be
p_\g(X|T)=\sup \frac {|h(X \cap T)|} {h(T)}
\ee
where $h:\R \to \R$ ranges over all $\g$-qs maps.  An important property of
$\g$-qs capacity is its behavior under tree decomposition: if $T^j
\subset T$ are disjoint intervals and $X \subset \cup T_j$ then
\be
p_\g(X|T) \leq p_\g(\cup T^j|T) \sup_j p_\g(X|T^j).
\ee

\comm{
\subsection{Holomorphic motions}

Given a domain $\VV$ in a complex Banach space  $E$ with a base point $*$
and a set  $X_*\subset\C$, a {\it holomorphic motion} of $X_*$ over $\VV$ is
a family
of injections $h_\la: X_*\ra\C$, $\la\in \VV$, such that $h_*=\id$ and
$h_\la(z)$ is
holomorphic in $\la$ for any $z\in X_*$. Let $X_\la = h_\la X_*$.

Holomorphic motions were studied in a series of works, specially \cite
{MSS}, \cite {BR} and \cite {Sl}.  The main properties of
holomorphic motions we will use are summarized in the following two lemmas.

Let $K:[0,1) \to \R$ be defined by $K(r)=(1+\rho)/(1-\rho)$ where $0 \leq
\rho <1$ is such that the hyperbolic distance between
$0$ and $\rho$ in $\DD$ is $r$.

\proclaim {$\lambda$-Lemma \cite {BR}.}
{
Let $h_\la: U_*\ra U_\la$ be a holomorphic motion of a domain $U_*\subset\C$
over a hyperbolic domain $D\subset\C$.  Then the maps $h_\la$ are
$K(r)$-qc, where $r$ is the hyperbolic distance between $*$ and $\la$ in
$D$.  Moreover, $K(r)=1+O(r)$ as $r\to 0$.
}

\proclaim {Extension Lemma \cite{Sl}.}
{
 A holomorphic motion $h_\la: X_*\ra X_\la$  of a set
$X_*\subset\C$ over a topological disk $D$ 
admits an extension to a holomorphic motion
$H_\la: \C\ra\C$ of the whole complex plane over $D$.  
}

Such an extension to the whole complex plane will also be called a
completion.

A symmetric holomorphic motion $h$ is a holomorphic motion such that
the image of each leaf of $h$ by complex conjugation (acting in $\C^2$) is
still a leaf of $h$.  The extension of Slodkovsky is not unique, but in  
\cite {AM1} it is remarked that the proof of Slodkovsky allows to obtain.

\proclaim{Real Extension Lemma}
{
Any symmetric
holomorphic motion over a simply connected domain can be completed to a
symmetric holomorphic motion.
}

Every time we consider a completion of
a symmetric holomorphic motion we will take a symmetric completion.
}

\subsection{Unimodal maps}

We refer to the book of de Melo and van Strien \cite{MS}
for the general background in one-dimensional dynamics. 

We will say that a smooth (at least $C^2$) map $f: I\ra I$ of the
interval $I=[-1,1]$ is {\it unimodal} if $f(-1)=-1$, $f(x)=f(-x)$
and $0$ is the only critical point of $f$ and is non-degenerate, so that
$D^2 f(0) \neq 0$.  The introduction of normalization and symmetry
in this definition is exclusively for the simplicity of the notation,
and is no loss of generality, see also Appendix C of \cite {ALM}.
The assumption of non-degeneracy of the critical point
is clearly typical.

\comm{
A smooth  map $f: I\ra I$ of the interval $I=[-1,1]$ is called {\it
unimodal}
if it has  a single critical point and this point is an extremum.
 We always assume that the critical point is located at the origin.
Let $\UU^3$ be the space of $C^3$ unimodal maps $f: I\ra I$ 
 with quadratic critical point, which are even (symmetric),
that is $f(x)=f(-x)$,
i.e., $Df(0)=0$ and $D^2(f(0)\neq 0$.  We normalize the maps so 
that $-1$ is a fixed point and $f(1)=-1$.  We endow $\UU^3$ with the $C^3$
topology.  If $Df(-1)<1$, either the
dynamics is trivial ($-1$ is the global attractor) or the map has a proper
unimodal restriction.  For this reason we will assume
further that $Df(-1) \geq 1$.
}

Basic examples of unimodal maps are given by quadratic maps
\begin{equation}\label{quadratic family}
q_\tau \colon I \to I,\quad  q_\tau(x)= \tau-1-\tau x^2,  
\end{equation}
where $\tau\in [1/2,2]$ is a real parameter.

Let $\U^k$, $k \geq 2$ be the space of $C^k$ unimodal maps.
We endow $\U^k$ with the $C^k$ topology.
A map $f \in \U^3$ is {\it quasiquadratic} if any nearby map $g \in
\U^3$ is topologically conjugate to some quadratic map.
We denote by $\U
\subset \U^3$ the space of quasiquadratic maps.  By the theory of
Milnor-Thurston and Guckenheimer \cite{MS},
a map $f\in \U^3$ with negative  Schwarzian derivative and $Df(-1)>1$
is quasiquadratic, so quadratic maps $q_\tau$, $\tau \in (1/2,2]$ belong
to $\U$.

A map $f \in \U^2$ is said to be {\it Kupka-Smale} if all
periodic orbits are
hyperbolic.  It is said to be hyperbolic if it is Kupka-Smale and
the critical point is attracted to a periodic attractor.  It is said to be
{\it regular} if it is hyperbolic and its critical point is not periodic or
preperiodic.  It is well known that regular maps are structurally stable.

In this paper, an analytic family of unimodal maps
will be understood as a one-parameter family
$\{f_\lambda \in \U^2\}_{\lambda \in \Lambda}$
(where $\Lambda \subset \R$ is an interval), such that
the correspondence $(\lambda,x) \mapsto f_\lambda(x)$ is analytic.
(The measure-theoretical
description of analytic families in several parameters follows from the
one-parameter case, see \cite {Av}.)


An analytic family of unimodal maps is called {\it non-trivial} if regular
parameters are dense.  If all maps in the family are quasiquadratic, it can
be shown that a family is non-trivial if it contains one regular
parameter (see Theorem A of \cite {ALM}).

\comm{
In this work we will only
consider quasiquadratic maps.

A quasiquadratic map $f$ is called {\it hyperbolic} or {\it regular}
if it has an attracting cycle $\overline q$.
In this case the orbit of the critical point 0 converges to
$\overline q$, hence a quasiquadratic map
can have at most one attracting cycle (Singer, see \cite{MS}).
Moreover, if it has one then almost all orbits converge to this cycle
(this follows from a result by Guckenheimer and Ma\~n\'e, see \cite{MS}).

A quasiquadratic map $f$ is called {\it parabolic} if it has a
parabolic cycle $\overline  q$.
Similarly to the hyperbolic case, in the parabolic case
the critical point $0$ belongs to the basin $D(\overline q)$
and this basin has full Lebesgue measure in $I$. Thus, a quasiquadratic map
can have at most one parabolic cycle (and the map cannot be simultaneously
hyperbolic and parabolic).
}

\comm{
\subsection{Spaces of unimodal maps}

Let $a>0$, and let  $\EE_a\subset \BB_{\Om_a}$ be the
complex Banach space of holomorphic maps
$v:\Om_a \to \C$ continuous up to the boundary
which are $0$-symmetric (that is, $v(z)=v(-z)$)
and such that $v(-1)=v(1)=0$,
endowed with the $\sup$-norm $\|v\|_a=\|v\|_\infty$.
It contains the real Banach space $\EE_a^{\R}$ of "real maps" $v$, 
i.e, holomorphic maps symmetric with respect to the real line: 
$v(\overline z) = \overline {v(z)}$. 

The  complex affine subspace $q_2+\EE_a$ will be denoted as $\AAA_a$.

If $f \in \AAA_a$, we denote the poscritical set
$\overline {\orb(f(0))}$ by $O_f$.

Let $\U_a=\U \cap \AAA_a$.  Note that $\U_a$ is the union of an open set 
in the affine subspace $\AAA_a^\R= q_2+ \EE_a^\R$ and a codimension-one
space of Ulam-Neumann maps.

\subsection{Laminations}

By definition, each quasiquadratic unimodal map is topologically equivalent
to some quadratic map.  In other words, each topological conjugacy class of
quasiquadratic maps intersect the quadratic family.  The following result
gives more information:

\begin{thm}

Each non-hyperbolic topological conjugacy class intersects the quadratic
family in one unique parameter.

\end{thm}

Of course the hyperbolic topological conjugacy classes intersect the
quadratic family in an open set.  For this result, it is sometimes
convenient to refine the hyperbolic topological conjugacy classes.  Let us
say that two unimodal maps are hybrid equivalent if they are topologically
equivalent and, moreover, if they are hyperbolic, the exponent of their
attracting periodic orbits is the same.  The hybrid class of a
quasiquadratic map $f$ is denoted $\HH_f$.
With this definition we have:

\begin{thm}

Each hybrid class intersects the quadratic
family in one unique parameter.

\end{thm}

One of the main results of \cite {ALM} is to describe the structure of
the partition in hybrid classes.  The main result is that each hybrid class
is a codimension-one analytic submanifold, which form an analytic
lamination outside countably many parabolic classes.

To be more precise, a codimension-one  {\it holomorphic
lamination} $\LL$ on an open subset $\VV \subset \AAA_a$
is a family of disjoint codimension-one  
Banach submanifolds of $\AAA_a$, called the {\it leaves} of the lamination
such that for any point $p\in \MM$, there exists a holomorphic local  chart
$\Phi: \WW \ra \VV \oplus \C$
(where $\VV$ is a neighborhood in some complex Banach space) 
such that for any leaf $L$ and any connected component $L_0$ of $L\cap \WW$,
the image $\Phi(L_0)$  is a graph of a holomorphic function $\VV \ra \C$.
  
The neighborhood $\WW$ in the above definition is called a {\it flow box},
and the connected components $L_0$  are called {\it local leaves} in this
flow box.

\begin{thm}[Theorem A of \cite {ALM}]

Every real hybrid class $\HH_f$, $f\in \U_a$,
  is an embedded codimension-one  real  analytic Banach submanifold of
$\U_a$.
 Furthermore, the hybrid classes laminate a neighborhood of any
non-parabolic map $f\in \U_a$.
 More precisely, such an $f$ has a neighborhood $\VV$
in  the complex affine space $\AAA_a$ endowed with
a codimension-one holomorphic lamination with the following properties:
\begin {itemize} 
\item the lamination is transverse to the real affine Banach subspace
$\AAA_a^{\
R}$; 
\item if $g\in \VV\cap \AAA_a^\R$ then the intersection of the leaf through 
$g$ with $\AAA_a^\R$ coincides with  $\HH_g^\R\cap \VV$. 
\end{itemize} 

\end{thm}

Moreover, the quadratic family is a universal transversal to the lamination.
}

\subsection{Renormalization}

Let $f \in \U^2$.  A symmetric (about $0$) interval $T \subset I$
is said to be {\it nice} if the iterates of $\partial T$ never return to
$\inter T$.  A nice interval $T \neq I$ is said to be a
restrictive (or periodic) interval of period $m$ for $f$ if
$f^m(T) \subset T$ and $m$ is
minimal with this property.  In this case,
the map $A \circ f^m \circ A^{-1}:I \to I$ is again unimodal for some affine
homeomorphism
$A:T \to I$ and is called a {\it renormalization}\footnote {A more usual
convention is to call $A \circ f^m \circ A^{-1}$ a unimodal restriction if
$m=1$, reserving the name renormalization for the case $m>1$, but we won't
make this distinction.} of $f$.  The map $f^m:T \to T$ will be called a
prerenormalization of $f$.

We say that $f$ is {\it infinitely renormalizable} if there exists
arbitrarily small restrictive intervals, and we say it is
{\it finitely renormalizable} otherwise.

Let $\FF \subset \U^2$ be the class of
Kupka-Smale finitely renormalizable maps whose critical point is
recurrent, but not periodic.

The following result shows that when investigating typical properties of
analytic unimodal maps, it is enough to deal with the quasiquadratic case.

\begin{thm}[Theorem B of \cite {Av}]

Let $f_\lambda$ be a non-trivial
analytic family of unimodal maps.  Then for almost
every non-regular parameter $\lambda$, $f_\lambda$ has a renormalization
which is quasiquadratic.

\end{thm}

It is easy to check that the conclusions of Theorems A, B, or C do not
depend on considering a map or its renormalization.  Due to this result,
in the arguments to follow, we will concentrate on the description of
quasiquadratic map and non-trivial analytic families of quasiquadratic maps.

\subsection{Some metric properties}

The condition of negative Schwarzian derivative plays an important role when
one needs to do distortion estimates.  One of the main tools is the
{\it Koebe Principle}:

\begin{lemma}[Koebe Principle, see \cite {MS}, page 258]

Let $f:T \to \R$ be a diffeomorphism with non-negative
Schwarzian derivative.  If $T' \subset T$ and both components
$L$ and $R$ of $T \setminus T'$ are bigger than $\epsilon |T'|$
then the distortion of $f|T'$ is bounded by $\frac {(1+\epsilon)^2}
{\epsilon^2}$.
In particular, we have $\min \{|f(L)|,|f(R)|\} \geq
\delta(\epsilon) \epsilon |f(T')|$, where $\delta(\epsilon)>0$ satisfies
$\lim_{\epsilon \to \infty} \delta(\epsilon)>\frac {9} {100}$.

\end{lemma}

Due to the recent results of Kozlovski, the condition of negative
Schwarzian is not needed for application of the Koebe Principle (for
unimodal maps in $\U^3$ which are Kupka-Smale),
see Theorem~B of \cite {K1} for instance.  We will thus apply the above
Koebe Principle without further comments in our setting.

\subsubsection{Hyperbolicity}

It was shown by Ma\~n\'e \cite {MS}
that (for one-dimensional maps of class $C^2$)
the obstruction to uniform expansion lies in critical points and
non-repelling periodic orbits.  Since quasiquadratic maps in $\FF$ do not
have non-repelling periodic orbits, this implies:

\begin{lemma} \label {hyperbolici}

Let $f \in \FF$ be a quasiquadratic map, and let $T$ be a
nice interval.  There
exists constants $C>0$, $\lambda>1$ such that
if $f^k(x) \in I \setminus T$,
$k=0,...,m-1$ then $|Df^m(x)|>C \lambda^m$.

\end{lemma}

\begin{cor}

Under the hypothesis of the previous lemma, if $K$ is a compact invariant
set which does not contain $0$, then $f|K$ is uniformly expanding.

\end{cor}

\comm{
\subsection{Combinatorics of unimodal maps}

There are many different presentations for the combinatorial theory of
unimodal maps.  The description in terms of kneading sequences of
Milnor-Thurston (which will be used in the formula claimed in Theorem A) is
not always convenient for our analysis.  We will need the description in
terms of renormalization and generalized renormalization that we give now.

\subsubsection{Nice intervals, return maps and landing maps}

Let $T \subset I$ be a symmetric interval which is nice, that is, for $k>0$,
$f^j(\partial T) \cap \inter T=\emptyset$.  We let $R_T$ denote the first
return map of $T$.  Due to the nice condition, the domain of $R_T$ is a
countable union of intervals $T^j \subset T$, indexed by integers,
where we let $T^0$ be the
central (that is, $0 \in T^0$), if it exists.  All branches, except the
central are diffeomorphisms.

Let $\Omega$ be the set of all finite sequences of non-zero integers
(possibly empty).  For
any element $\d \in \Omega$, $\d=(j_1,...,j_m)$ we associate a branch
$R^\d_T$ of $R^m_T$, whose domain is $\{x \in T|R^k(x) \in T^{j_{k+1}}, 0
\leq k < m\}$.


If $0$ belongs to the domain of $R_T$, we consider the first landing map
$L_T$ from $T$ to $T^0$.  The domain of $L_T$ is simply the union of
intervals $C^\d_T=(R^\d_T)^{-1}(T^0)$, and $L_T|C^\d_T=R^\d_T$.

Notice that $R_{T^0}=L_T \circ R_T$.

\subsubsection{Renormalization}

We say that $f$ is renormalizable
if there exists an interval $T \subset I$ and an integer $n>1$ such
that $f^n(T) \subset T$ and $f^k(\inter T) \cap \inter T=\emptyset$, $1 \leq
k <n$.  The integer $n$ is called the period of $T$.  The maximal interval
of period $m$ is called the renormalization interval of period $m$.

The set of renormalization periods of a unimodal map is an increasing
(possibly empty) sequence $m_i$, $i \geq 1$, each associated to a
renormalization interval $T^{(i)}$.  We set also $m_0=1$, $T^{(0)}=I$.

\subsubsection{The principal nest}

Let $\Delta_k$ be the set of all maps $f$ such that
$f$ is $k$-times renormalizable and
$f^{m_k}:T^{(k)} \to T^{(k)}$ has a fixed point $p_k$
in $\inter T^{(k)}$ and $Df^{m_k}(p_k) \leq -1$.  We define
$T^{(k)}_0=[-p_k,p_k]$.  If $0 \in T^{(k)}_0$ we let $T^{(k)}_1$ be the
central component of the domain of the first return map to
$T^{(k)}_0$.  Continuing in this way, we define a sequence (possibly finite)
of intervals $T^{k)}_i$, such that $T^{(k)}_{i+1}$ is the central domain of
the first return map to $T^{(k)}_i$.

It is possible to show that $T^{(k)}_i$ if either $f$ is $k+1$-times
renormalizable or $f$ is exactly $k$-times renormalizable and has a
recurrent critical point.  Those possibilities correspond to $\cap
T^{(k)}_i$ be an interval or $\{0\}$, in the former case, it is
$T^{(k+1)}$.

Assume that we are working in
some fixed level $k$ of renormalization.  In this
case we will simplify the notation dropping $k$ by letting $I_n=T^{(k)}_n$.
We let $R_n:\cup I^j_n \to I_n$ be the first return map and $L_n:\cup C^\d_n
\to I_{n+1}$ be the first landing map.
}

\subsection{Physical measures}

Let $\mu$ be a probability measure which is invariant under the dynamics of
$f$.  The {\it basin} of $\mu$ is the set of points $x \in I$ such that
\be
\lim \frac {1} {m} \sum_{k=0}^{m-1} \delta_{f^k(x)}=\mu
\ee
in the weak sense.  We say that $\mu$ is a {\it physical measure} if the
basin of $\mu$ has positive Lebesgue measure.
A quasiquadratic map can have at most one physical
measure \cite {BL},
which (if it exists) has always a basin of full Lebesgue measure.
If $f$ is hyperbolic, then the uniform distribution in the attracting
periodic orbit is the physical measure of $f$.  If $f$ is stochastic, that
is, it has an absolutely continuous invariant measure $\mu$,
then this measure is ergodic and, by Birkhoff's Theorem, it is a physical
measure.  Notice that there exist quadratic maps without any physical
measure, see \cite {MS}, Chapter V, Section 5.

If $f$ is stochastic, then it is finitely renormalizable.  Let
$f^k:T \to T$ be its last prerenormalization.  It turns out
that the support of $\mu$ is $A=T_0 \cup ... \cup T_{k-1}$
where $T_0=[f^{2k}(0),f^k(0)]$ and $T_j=f^j(T_0)$.  Notice that
$f^k(T_0)=T_0$.  We could have defined $A$ topologically in this
way without any reference to $\mu$.

The set $A$ has another remarkable
property: it is the smallest compact subset of $I$ such that

\begin{enumerate}

\item for almost every $x \in I$, $\omega(x) \subset A$;

\item for generic $x \in I$, $\omega(x) \subset A$.

\end{enumerate}

Those two conditions mean exactly that
$A$ is the topological and metric attractor of $f$ in the
sense of Milnor.

\begin{rem}

All quasiquadratic unimodal maps have a unique topological and a unique
metric attractor.  Both concepts of attractor coincide
by \cite {attractors}.

\end{rem}

A sufficient condition for $f$ to be stochastic is the Collet-Eckmann
condition: $|Df^n(f(0))|$ grows exponentially fast.

\begin{thm}[Corollary C of \cite {Av}]

Let $f_\lambda$ be a non-trivial family of analytic unimodal maps.  Then
almost every non-regular parameter belongs to $\FF$ and
satisfies the Collet-Eckmann condition.

\end{thm}

\comm{
This condition has been
improved a lot, it is known that it is enough to assume that
$\sum |Df^n(f(0))|^{-1/2}<\infty$.  Lyubich's
Regular or Stochastic dichotomy for the quadratic family was originally
obtained by showing that almost
every non-regular quadratic map satisfies this summability condition.
}

We will need the following result of Keller about general
stochastic unimodal maps:

\begin{thm}[see \cite {MS}, Theorem 3.2, Chapter V] \label {Below}

Let $f$ be a quasiquadratic stochastic map, and let $\mu$ be its physical
measure.  Then $d\mu$ is uniformly bounded from below on $A$.

\end{thm}

\begin{rem}

Keller's Theorem is stated in \cite {MS} for maps with negative Schwarzian
derivative.  The result for quasiquadratic maps can be obtained with the
same proof using the results of Kozlovski \cite {K1}.

\end{rem}

Notice that while $d\mu$ is always bounded from below, it
is definitely not bounded from above, and we will
need to work a lot to obtain in Theorem C a
reasonable estimate for $d\mu$.  Notice
also that our proof of Theorem C is not a general one for stochastic maps:
we have to exclude lots of them.  (It is easy to see that some exclusion has
to be done, for instance, one must exclude
stochastic maps with non-recurrent critical point.)

\comm{
\subsection{Statistical estimates} \label {values}

Several estimates obtained in \label {AM} will be used here, we collect the
main ones here:

\begin{lemma}

Almost every non-regular map satisfies
\be
\lim \frac {\ln v_{n+1}} {\ln c_n^{-1}}=\lim \frac {\ln s_n} {\ln c_n^{-1}}=
\lim \frac {\ln \ln c_{n+1}^{-1}} {\ln c_n^{-1}}=\lim \frac {\ln
r_n(\tau_n)} {\ln c_{n-1}^{-1}}=1,
\ee
\be
I(\{r_n(j)>c_{n-1}^4 k\},n) \leq e^{-k}, \quad k \geq 1.
\ee

\end{lemma}
}

\section{The formula} \label {formula}

\subsection{Combinatorics}

Let us have a symbol space $\Sigma$ with finitely many elements.  A (finite
or infinite) sequence of elements of $\Sigma$ will be called a word.  In the
space $\Sigma^\N$ of infinite words, we let the shift operator $\sigma$
act by $\sigma(\alpha_0\alpha_1...)=\alpha_1\alpha_2...$.

Given a finite word $\alpha$ and $r \in \N \cup \{\infty\}$, we let
$\alpha^r$ denote $r$ repetitions of $\alpha$.

A finite word $\alpha$ is said to be irreducible if $\alpha=\beta^r$ for
some $r$ implies $\alpha=\beta$.

If $\alpha$ is an infinite word which is periodic, there exists a unique
irreducible word $\overline \alpha$ such that $\alpha=\overline
\alpha^\infty$.

\subsubsection{Frequencies}

If $\alpha=\alpha_0...\alpha_{m-1}$ is a finite word
and $\beta=\beta_0\beta_1...$ is an infinite word, we define the
lower and upper frequencies of $\alpha$ in $\beta$ in the natural way:
\be
r^+(\alpha,\beta)=\limsup_{n \to \infty} \frac {1} {n} \#\{0 \leq k \leq
n-1|\alpha_i=\beta_{k+i}, 0 \leq i \leq m-1\},
\ee
\be
r^-(\alpha,\beta)=\liminf_{n \to \infty} \frac {1} {n} \#\{0 \leq k \leq
n-1|\alpha_i=\beta_{k+i}, 0 \leq i \leq m-1\}.
\ee

The frequency $r(\alpha,\beta)$ is defined as the common value of
$r^+(\alpha,\beta)$ and $r^-(\alpha,\beta)$ if they coincide.
We say that $\beta$ is normal if, for any $\alpha$,
$r^+(\alpha,\beta)=r^-(\alpha,\beta)$.

\subsubsection{Geometric frequencies}

Let $\alpha$ be a finite word and $\beta$ be a normal infinite word.
Let us consider the non-increasing sequence $r(\alpha^k,\beta)$.
We want to associate to $\alpha$ and $\beta$ a quantity related
to the decay of
$r(\alpha^k,\beta)$.  In the case of exponential decay,
it is natural to define the upper and lower geometric
frequencies:
\be
\rho^+(\alpha,\beta)=\limsup_{n \to \infty} r(\alpha^n,\beta)^{1/n},
\ee
\be
\rho^-(\alpha,\beta)=\liminf_{n \to \infty} r(\alpha^n,\beta)^{1/n}.
\ee

The geometric frequency $\rho(\alpha,\beta)$ is the common value of
$\rho^+(\alpha,\beta)$ and $\rho^-(\alpha,\beta)$ if they coincide.
We say that $\beta$ is geometrically normal if for any
$\alpha$, $\rho^+(\alpha,\beta)$ and $\rho^-(\alpha,\beta)$ coincide.

\subsection{Itineraries}

Let us associate to an unimodal map $f$ some symbolic dynamics.  We fix the
symbol space $\Sigma=\{0,c,1\}$.  Let
$\Theta:I \to \Sigma$ be defined by $\Theta|[-1,0)=0$, $\Theta|(0,1]=1$,
and $\Theta(0)=c$.

The {\it itinerary} of a point $x \in I$ is the infinite word
$\theta(x)=\theta_0\theta_1...$, where $\theta_k=\Theta(f^k(x))$.

The (discontinuous)
map $\theta:I \to \Sigma^\N$ satisfy $\theta \circ f=\sigma \circ
\theta$.  It is clear that if $p$ is a periodic point for $f$, then
$\theta(p)$ is a periodic word for $\sigma$.

Given a word $\alpha$, we let $I_\alpha \subset I$ be the set of points
whose itinerary starts with $\alpha$.  Depending on $\alpha$,
$I_\alpha$ can be either an interval, a point or empty.

\subsection{Proof of Theorem~A assuming Theorems B and C}

We will actually prove the following stronger:

\begin{thm}

Let $f$ be a quasiquadratic unimodal map such that

\begin{enumerate}

\item $f$ is Collet-Eckmann and has an absolutely continuous invariant
measure $\mu$ supported in a cycle of intervals $A$;

\item $0$ belongs to the basin of $\mu$;

\item For any invariant hyperbolic set $K$, and any $1 \leq p <\infty$,
$d\mu^K_f \in L^p$.

\end{enumerate}

Then $\theta(0)$ is geometrically normal and for any
$p \in A$ periodic (of period $m$),
\be \label {rhoformula}
|Df^m(p)|=\rho(\overline {\theta(p)},\theta(0))^{-1}.
\ee
Moreover, for any $\alpha$ such that
$\rho(\alpha,\theta(0))>0$, there exists a periodic orbit $p \in A$
such that $\theta(p)=\alpha^\infty$.

\end{thm}

\begin{pf}

Let $\theta(0)=\theta_0\theta_1...$.  Let $\alpha=\alpha_0...\alpha_{m-1}$
be an arbitrary finite word.  Notice that
$\theta_{k+j}=\alpha_j$, $0 \leq j
\leq m-1$, if and only if $f^k(0) \in I_\alpha$,
so by definition of basin of $\mu$,
$r(\alpha,\theta(0))=\mu(I_\alpha)$.  In particular, $\theta(0)$ is normal.

Let $p \in A$ be a periodic orbit, and let $\alpha=\overline
{\theta(p)}$.  By item 1, we conclude that $p$ is repelling, and
since $f$ is quasiquadratic, $\cap I_{\alpha^k}=p$, and the length $m$ of
$\alpha$ is the period of $p$.  Let
$q,q' \in I_\alpha$ be periodic orbits in opposite sides of $p$, and let
$q_k=(f^{km}|I_{\alpha^{k+1}})^{-1}(q)$ and
$q'_k=(f^{km}|I_{\alpha^{k+1}})^{-1}(q')$. 
Let $K$ be the hyperbolic set consisting of $p$, the forward orbit of $q$
and $q'$ and all $q_k$ and $q'_k$.  Let $T_k=[q'_k,q_k]$.
It is easy to see that there exists
$j>0$ such that for all $k>j$,
\be
T_{k+j} \subset I_{\alpha^k} \subset T_{k-j}.
\ee

In particular,
\be
\rho^+(\alpha,\theta(0))=\limsup_{n \to \infty} \mu(T_n)^{1/n},
\ee
\be
\rho^-(\alpha,\theta(0))=\liminf_{n \to \infty} \mu(T_n)^{1/n}.
\ee

By Theorem \ref {Below},
there exists a constant $C>0$ such that $d\mu|A \geq C$.
On the other hand, $T_k \subset A$ for $k$ big enough, so
$\mu(T_k) \geq C|[q'_k,q_k]|$.  It is clear that
\be
\lim_{n \to \infty}
|T_n|^{1/n}=|Df^m(p)|^{-1},
\ee
so $\rho^-(\alpha,\theta(0)) \geq
|Df^m(p)|^{-1}$.

By item 3, for all $1 \leq p <\infty$, there exists a constant
$C_p$ such that, for all $k \geq 0$,
\be
\left ( \int_{T_k} (d\mu^K)^p \right )^{1/p} \leq C_p.
\ee
In particular, by the H\"older inequality,
\be
\mu(T_k)=\int_{T_k} d\mu^K \leq
\left ( \int_{T_k} (d\mu^K)^p \right )^{\frac {1} {p}}
\left ( \int_{T_k} 1^{\frac {p} {p-1}} \right )^{1-\frac {1} {p}} \leq
C_p |T_k|^{1-\frac {1} {p}}.
\ee
Taking $k \to \infty$ we get
\be
\rho^+(\alpha,\theta(0)) \leq |Df^m(p)|^{-1+\frac {1} {p}}.
\ee
Since $1 \leq p < \infty$ is arbitrary, we get (\ref {rhoformula}).

If $\alpha$ is an arbitrary finite word, then
either $I_{\alpha^k} \cap A$ is eventually empty or
$\cap I_{\alpha^k}$ is a repelling periodic orbit $p$ in
$A$.  In the first case, obviously $\rho(\alpha,\theta(0))=0$. 
In the second case, by the previous discussion,
$\rho(\alpha,\theta(0))=|Df^m(p)|^{-1}>0$, where $m$ is the length of
$\alpha$.  In particular, $\theta(0)$ is geometrically normal.
\end{pf}

\begin{rem}

Let us note that the Collet-Eckmann condition already implies a number of
interesting properties (see \cite {NS}).  For instance, if $f$ is
a quasiquadratic Collet-Eckmann map, then there exists a constant
$\lambda>1$ such that if $p$ is a periodic orbit of period $n$ then
$|Df^n(p)| \geq \lambda^n$.

\end{rem}


\comm{
\subsection{Combinatorics of unimodal maps}

\subsubsection{Nice intervals, return maps and landing maps}

Let $T \subset I$ be a symmetric interval which is nice, that is, for $k>0$,
$f^j(\partial T) \cap \inter T=\emptyset$.  We let $R_T$ denote the first
return map of $T$.  Due to the nice condition, the domain of $R_T$ is a
countable union of intervals $T^j \subset T$, indexed by integers,
where we let $T^0$ be the
central (that is, $0 \in T^0$), if it exists.  All branches, except the
central are diffeomorphisms.

Let $\Omega$ be the set of all finite sequences of non-zero integers.  For
any element $\d \in \Omega$, $\d=(j_1,...,j_m)$ we associate a branch
$R^\d_T$ of $R^m_T$, whose domain is $\{x \in T|R^k(x) \in T^{j_{k+1}}, 0
\leq k < m\}$.

Let $\Omega_0$ denote $\Omega$ without the empty sequence.
We let $\sigma^+,\sigma^-:\Omega_0 \to \Omega$,
$\sigma^+(j_1,...,j_m)=(j_1,...,j_{m-1})$,
$\sigma^-(j_1,...,j_m)=(j_2,...,j_m)$.

If $0$ belongs to the domain of $R_T$, we consider the first landing map
$L_T$ from $T$ to $T^0$.  The domain of $L_T$ is simply the union of
intervals $C^\d_T=(R^\d_T)^{-1}(T^0)$, and $L_T|C^\d_T=R^\d_T$.

Notice that $R_{T^0}=L_T \circ R_T$.

\subsubsection{Renormalization}

We say that $f$ is renormalizable
if there exists an interval $T \subset I$ and an integer $n>1$ such
that $f^n(T) \subset T$ and $f^k(\inter T) \cap \inter T=\emptyset$, $1 \leq
k <n$.  The integer $n$ is called the period of $T$.  The maximal interval
of period $m$ is called the renormalization interval of period $m$.

The set of renormalization periods of a unimodal map is an increasing
(possibly empty) sequence $m_i$, $i \geq 1$, each associated to a
renormalization interval $T^{(i)}$.  We set also $m_0=1$, $T^{(0)}=I$.

\subsubsection{The principal nest}

Let $\Delta_k$ be the set of all maps $f$ such that
$f$ is $k$-times renormalizable and
$f^{m_k}:T^{(k)} \to T^{(k)}$ has a fixed point $p_k$
in $\inter T^{(k)}$ and $Df^{m_k}(p_k) \leq -1$.  We define
$T^{(k)}_0=[-p_k,p_k]$.  If $0 \in T^{(k)}_0$ we let $T^{(k)}_1$ be the
central component of the domain of the first return map to
$T^{(k)}_0$.  Continuing in this way, we define a sequence (possibly finite)
of intervals $T^{k)}_i$, such that $T^{(k)}_{i+1}$ is the central domain of
the first return map to $T^{(k)}_i$.

It is possible to show that $T^{(k)}_i$ if either $f$ is $k+1$-times
renormalizable or $f$ is exactly $k$-times renormalizable and has a
recurrent critical point.  Those possibilities correspond to $\cap
T^{(k)}_i$ be an interval or $\{0\}$, in the former case, it is
$T^{(k+1)}$.
}

\section{Phase-parameter relation and statistics of the principal nest}
\label {pp}

In this section we will discuss the principal nest combinatorics, and then
state the Phase-Parameter relation, which is our means to obtain parameter
results based on phase estimates.  We will then present some results on the
statistics of the principal nest.

\subsection{Principal nest combinatorics}

If $T \subset I$ is a nice interval, the domain of the
first return map $R_T$ to $T$ consists of a (at most)
countable union of intervals which we denote $T^j$.  We reserve the
index $0$ for the component of $0$: $0 \in T^0$, if $0$ returns to $T$.
From the nice condition, $R_T|T^j$ is a diffeomorphism if
$0 \notin T^j$, and is an even map if $0 \in T^j$.  The domain containing
$0$ will be called the central domain of $R_T$ and will be denoted $T^0$.
The return $R_T$ is said
to be central if $R_T(0) \in T^0$.  If $f$ is quasiquadratic with recurrent
but not periodic critical point, the domain of the first return map is dense
and its complement is a regular Cantor set.

Let $f \in \FF$ be quasiquadratic,
and let $T$ be its smallest restrictive interval (of period $m'$).
Define a sequence of nested nice intervals $I_n$ by induction as follows.
Let $I_0=[-p,p]$ where $p$ is the unique orientation reversing fixed
point of $f^{m'}:T \to T$.  Assuming $I_n$ defined, let
$R_n:I_n \to I_n$ be the first return map and $I_{n+1}=I^0_n$.
Since $f$ is finitely renormalizable, $\cap I_n=\{0\}$.

Let $\Omega$ be the set of all finite sequences of non-zero integers
(possibly empty).  For
any element $\d \in \Omega$, $\d=(j_1,...,j_m)$ we associate a branch
$R^\d_n$ of $R^m_n$, whose domain is $I^\d_n=\{x \in I_n|R^k(x) \in
I^{j_{k+1}}_n, 0 \leq k < m\}$.

Let $L_n:I_n \to I^0_n$ be the first landing map.  The domain of $L_n$ is
the union of intervals $C^\d_n=(R^\d_n)^{-1}(I^0_n)$.

\subsection{Phase-Parameter relation}

We will now quickly define formally the Phase-Parameter relation, and we
will discuss in the next section the way it is applied for
measure-theoretical problems.

\begin{definition}

Let us say that a family $f_\lambda$ of quasiquadratic maps satisfies the
Topological Phase-Parameter relation at a parameter $\lambda_0$ if
$f=f_{\lambda_0} \in \FF$, and there exists $i_0>0$ and
a sequence of nested intervals $J_i$, $i \geq i_0$ such that: 

\begin{enumerate}

\item $J_i$ is the maximal interval containing $\lambda_0$
such that for all
$\lambda \in J_i$ there exists a homeomorphism $H_i[\lambda]:I \to I$
such that
$f_\lambda \circ H_i[\lambda]|(I \setminus I_{i+1})=H_i[\lambda] \circ f$.

\item There exists a homeomorphism $\Xi_i:I_i \to J_i$ such that
$\Xi_i(C^\d_i)$ (respectively, $\Xi_i(I^\d_i)$)
is the set of $\lambda$ such that
the first return of $0$ to $H_i[\lambda](I_i)$
under iteration by $f_\lambda$ belongs to $H_i[\lambda](C^\d_i)$
(respectively, $H_i[\lambda](I^\d_i)$).

\end{enumerate}

\end{definition}

Let $K_i$ be the closure of the union of all
$\partial C^\d_i$ and $\partial I^\d_i$.
Notice that $H_i$ and $\Xi_i$ are only uniquely defined in $K_i$.
Condition (2) of the Topological Phase-Parameter relation can be
equivalently formulated as the existence of a homeomorphism $\Xi_i:I_i \to
J_i$ such that the first return of the critical point (under iteration by
$f_\lambda$) to $H_i[\lambda](I_i)$ belongs to
$H_i[\lambda](K_i)$ if and only if $\lambda \in \Xi_i(K_i)$.

Let us assume we have a non-trivial family of unimodal maps satisfying the
Phase-Parameter relation at a parameter $f=f_{\lambda_0}$.
It will be important to
estimate the metric properties of $H_i|K_i$ and $\Xi_i|K_i$.

Let $\tilde I_{i+2}=(R_i|I^0_i)^{-1}(I^\d_i)$, where $\d$ is such that
$(R_i|I^0_i)^{-1}(C^\d_i)=I_{i+2}$.

Let $\tau_i$ be such that $R_i(0) \in I^{\tau_i}_i$.  Let
$\tilde K_i=\overline {(\cup_j \partial I^j_i \cup \partial I_i)
\setminus \tilde I_{i+1}}$.

Let $J^j_i=\Xi_i(I^j_i)$.

Let us say that $f \in \FF$ is simple if only finitely many
$R_n$ have central returns.

\begin{definition}

Let $f_\lambda$ be a family of unimodal maps.  We say that $f_\lambda$
satisfies the Phase-Parameter relation at $\lambda_0$ if $f=f_{\lambda_0}$
is simple, $f_\lambda$ satisfies the
Topological Phase-Parameter relation at $\lambda_0$ and for every $\g>1$,
there exists $i_0>0$ such that for $i>i_0$ we have:

\begin{description}

\item[PhPa1] $\Xi_i|(K_i \cap I^{\tau_i}_i)$ is $\g$-qs,

\item[PhPa2] $\Xi_i|\tilde K_i$ is $\g$-qs,

\item[PhPh1] $H_i[\lambda]|K_i$ is $\g$-qs if $\lambda \in
J^{\tau_i}_i$,

\item[PhPh2] the map $H_i[\lambda]|\tilde K_i$ is
$\g$-qs if $\lambda \in J_i$.

\end{description}

\end{definition}

\begin{thm}[Theorem A of \cite {Av}]

Let $f_\lambda$ be a non-trivial analytic family of quasiquadratic maps. 
Then $f_\lambda$ satisfies the Phase-Parameter relation at
almost every non-regular parameter.

\end{thm}

(Theorem A of \cite {Av} actually
covers the non-quasiquadratic case as well.)

\comm{
\subsection{Using the Phase-Parameter relation}

Let $f_\lambda$ be a non-trivial analytic family of quasiquadratic
maps, and let $X$ be a set of parameters satisfying the Phase-Parameter
relation plus some set of properties which we were able to prove are
typical among non-regular analytic unimodal maps, so
that $X$ has full measure among non-regular parameters of $f_\lambda$.  For
instance, at this point we could ask that parameters in $X$ are
Collet-Eckmann.

We would like to prove that some other property is still typical among
non-trivial analytic unimodal maps.  Thus we must show that the set of
properties

Let us assume that we were able to prove that a certain property $A$
is typical among non-regular analytic unimodal maps.  In other words,
in any non-trivial family of
analytic unimodal maps, the set $X$ of parameters satisfying $A$
has full measure among non-regular maps.  We want to prove some other
property $B$is still typical, that is, if $Y$ is the set of parameters
satisfying $B$ then $|X \setminus Y|=0$.

In practice, property $A$ will always be a collection of properties
including being quasiquadratic, simple, plus whatever statistical properties
that have already been proved to be typical, and $B$ will be a certain new
statistical property that we want to prove

\begin{definition}

We say that $f_\lambda$ satisfies the $n_0$-topological
phase-parameter relation at $f$ if
there exists a sequence of intervals $\{J_n\}_{n \geq n_0}$ with the
following properties:

\begin{enumerate}

\item $J_n$ is the maximal interval containing $f$ such that for all $g \in
J_n$, the interval $T^{(k)}_{n+1}[g]$ is defined and changes in a continuous
way;

\item For every $n \geq n_0$, there exists a continuous
family of homeomorphisms $\{H_n[g]\}_{g \in J_n}$ such that for each $g \in
J_n$, $H_n[g]:I_n \to I_n[g]$ preserves the domains of landing maps:
$H_n[g](C^\d_n)=C^\d_n[g]$;

\item For every $n \geq n_0$, there exists a homeomorphism $\Xi_n:I_n \to
J_n$ such that $\Xi_n(C^\d_n)$ is the set of all $g \in J_n$ such that
$R_n[g](0) \in C^\d_n[g]$.

\end{enumerate}

\end{definition}

\subsection{Metric phase-parameter relation}

\begin{definition}

We say that $f_\lambda$ satisfies the phase-parameter relation at $f$ if
for all $\g>1$, there exists $n_0 \in \N$ such that
$f_\lambda$ satisfies the $n_0$-topological phase-parameter relation at
$f$ and moreover, for all $n \geq n_0$:

\begin{description}

\item[PhPa1] $\Xi_n|_{K^\tau_n}$ is $\g$-qs,

\item[PhPa2] $\Xi_n|_{\tilde K_n}$ is $\g$-qs,

\item[PhPh1] $H_n[g]|_{K_n}$ is $\g$-qs if $g \in J^{\tau_n}_n$,

\item[PhPh2] $H_n[g]|_{\tilde K_n}$ is
$\g$-qs if $g \in J_n$.

\end{description}

\end{definition}

\begin{thm}

Let $f$ be an analytic unimodal map which is finitely renormalizable and has
a recurrent critical point.
Let $f_\lambda$ be an analytic family of unimodal maps transverse at $f$.
Then there exists $n_0$ such that $f_\lambda$
satisfies the topological phase-parameter relation at $f$.
If $f$ has only finitely many central cascades, then it also satisfies the
phase-parameter relation.

\end{thm}

We will postpone the proof of Theorem until \S, in order to get directly
into the statistical analysis part of the paper.

\begin{thm}

Let $f_\lambda$ be a non-trivial family of unimodal maps.  For almost every
non-regular $\lambda$, $f_\lambda$ is finitely renormalizable and has a
recurrent critical point and finitely many central cascades.

\end{thm}
}

\subsection{Using the Phase-Parameter relation}

Let us now explain how the Phase-Parameter relation can be used to prove
that some property is typical among non-regular analytic unimodal maps.

\medskip
Notice that, due to the previous results, it is enough to prove that the
property is satisfied by almost every parameter in a non-trivial analytic
family of quasiquadratic maps.  {\it From now on we shall always work inside
such a fixed family.}  We can further restrict our scrutiny to the subset of
parameters which are simple and satisfy the
Phase-Parameter relation.  It is also clearly enough to
restrict ourselves to the analysis of
unimodal maps which are exactly $k$-times renormalizable for some fixed (but
arbitrary) $k$.  {\it We shall use  ``with total probability'' to denote
some property that is valid for a full measure set of parameters under those
restrictions.}

\medskip
We will now illustrate the basic principle with an example worked out in
\cite {AM1}.

For a simple map $f=f_{\lambda_0}$ which is quasiquadratic, simple and
satisfies the Phase-Parameter relation,
let us associate a sequence of ``statistical
parameters'' in the following way.  Let $s_n$ be the
number of times the critical point $0$ returns to $I_n$ before the first
return to $I_{n+1}$.  Let $c_n=|I_{n+1}|/|I_n|$.  Each of the points of the
sequence $R_n(0)$,...,$R_n^{s_n}(0)$ can be located anywhere
inside $I_n$.  Pretending that the distribution of those points is indeed
uniform with respect to Lebesgue measure, we may expect that $s_n$ is about
$c_n^{-1}$.

Let us try to make this rigorous.  Consider the set of points $A_k \subset
I_n$ which iterate exactly $k$ times in $I_n$ before entering $I_{n+1}$. 
Then most points $x \in I_n$ belong to some $A_k$ with $k$ in a neighborhood
of $c_n^{-1}$ (to be computed precisely using a statistical argument, in
this case, fixing some small $\epsilon>0$,
we can take the neighborhood to be
$c_n^{-1+2\epsilon}<k<c_n^{-1-\epsilon}$ for $n$ big).  By
most, we mean that, say, the complement has at most probability $\alpha_n$
which is some summable sequence.  In this case, it is not hard to see that
we can take $\alpha_n=c_n^\epsilon$, which indeed decays exponentially,
and so is summable, for all simple maps $f$ by \cite {attractors}.

If the phase-parameter relation were Lipschitz, we would now argue as
follows: the probability of a parameter be such that $R_n(0) \in A_k$ with
$k$ out of the ``good neighborhood'' of values of $k$ is also summable
(since we only multiply those probabilities by the Lipschitz constant) and
so, by Borel-Cantelli, for almost every parameter this only happens a finite
number of times.  More precisely, we would use the following version of
Borel-Cantelli:

\begin{lemma}[Lemma 3.1 of \cite {AM1}] \label {measure-theoretical}

Let $X \subset \R$ be a measurable set such that for each
$x \in X$ is defined a sequence
$D_n(x)$ of nested intervals converging to $x$
such that for all $x_1,x_2 \in X$
and any $n$, $D_n(x_1)$ is either equal or disjoint to $D_n(x_2)$.  Let
$Q_n$ be measurable subsets of
$\R$ and $q_n(x)=|Q_n \cap D_n(x)|/|D_n(x)|$.  Let $Y$ be
the set of $x$ in $X$ which belong to finitely many $Q_n$.
If $\sum q_n(x)$ is finite for almost any $x \in X$ then $|Y|=|X|$.

\end{lemma}

Unfortunately, the Phase-Parameter relation is not Lipschitz.  To make the
above argument work, we must have better control of the size of the ``bad
set'' of points which we want the critical value to not fall into.  In order
to do so, in the statistical analysis of the sets $A_k$
we control the quasisymmetric capacity (instead of Lebesgue measure)
of the complement of the set of points
whose entrance times belong to the
good neighborhood.  This makes the analysis sometimes much
more difficult: capacities are not probabilities (since they are not
additive), so we can have two disjoint sets with capacity close to $1$. 
This will usually introduce some error that was not present in the naive
analysis: this is the $\epsilon$ in the exponents present above.  If we were
not forced to deal with capacities, we could get much finer estimates.

Incidentally, to keep the error low, making $\epsilon$ close to $0$,
we need to use capacities with constant $\g$ close to $1$.  It will indeed
be very important for us that the
Phase-Parameter relation we use provides constants near $1$,
since this will allow us to partially get rid of those error terms.
This is also the reason that the estimates in \cite {AM2}
(which employed weaker Phase-Parameter estimates)
are worse than \cite {AM1}.

Coming back to our problem, we see that we should concentrate in proving
that for almost every parameter, certain bad sets have summable
$\g$-qs capacities for some constant $\g$ independent of $n$ (but which can
depend on $f$).

There is one final detail to make this idea work in this case: there are two
Phase-Parameter statements, and we should use the right one.  More
precisely, there will be situations where we are analyzing some sets which
are union of $I^j_n$ (return sets), and sometimes union of $C^\d_n$ (landing
sets).  In the first case, we should use the PhPa2 and in the second the
PhPa1.  Notice that
our Phase-Parameter relations only allow us to ``move the
critical point'' inside $I_n$ with respect to the partition by $I^j_n$, to
do the same with respect to the partition by $C^\d_n$, we must restrict
ourselves to $I^{\tau_n}_n$.  In all cases, however, the bad sets
considered should be either union of $I^j_n$ or $C^\d_n$.

For our specific example, the $A_k$ are union of $C^\d_n$, and we must use
PhPa1.  In particular we have to study the capacity of a bad set inside
$I^{\tau_n}_n$.  Here is the estimate that we should go after:

\begin{lemma}

For almost every parameter,
for every $\epsilon>0$, there exists $\g>1$ such that
$p_\g(X_n|I^\tau_n)$
is summable, where $X_n$ is the set of points $x \in I_n$ which enter
$I_{n+1}$ either before $c_n^{-1+\epsilon}$ or after $c_n^{-1-\epsilon}$
returns to $I_n$.

\end{lemma}

And as a consequence of PhPa1 we get:

\begin{lemma}

With total probability, for all $\epsilon>0$, for all $n$ sufficiently big,
\be
c_n^{-1+\epsilon}<s_n<c_n^{-1-\epsilon}.
\ee

\end{lemma}

In the language of Lemma \ref {measure-theoretical},
$X$ would be the set of
simple quasiquadratic parameters satisfying the
Phase-Parameter relation and which are exactly
$k$-times renormalizable, $D_n(\lambda)$, $\lambda \in X$ would be
$J^{\tau_n}_n(\lambda)$, and $Q_n \subset X$ would be the set of parameters
such that either $s_n<c_n^{-1+\epsilon}$ or $s_n>c_n^{-1-\epsilon}$.

\comm{
We note that we can use $s_n$ to estimate $c_{n+1}$ directly.
This last lemma implies that $c_n$ decays torrentially to $0$ for typical
parameters.  For general parameters (with finitely many central cascades),
the best information is
given by \cite {attractors}: $c_n$ decays exponentially
(this was actually used in the above argument).  This improvement from
exponential to torrential should give the reader an idea of
the power of this kind of statistical analysis.
}

\subsection{Some results on the statistics of the principal nest} \label
{values}

Let us collect here some results of \cite {AM1} on the dynamics of typical
non-regular analytic unimodal maps (the results were initially proved in the
quadratic setting, but hold in general due to \cite {Av}).

Let $r_n(j)$ be such that $R_n|I^j_n=f^{r_n(j)}$.  For $x \in I^j_n$, we let
$r_n(x)=r_n(j)$.  Let $l_n(\d)$ be such that $L_n|C^\d_n=f^{l_n(\d)}$, and
for $x \in C^\d_n$, let $l_n(x)=l_n(\d)$.  Let $v_n=r_n(0)$.  Recall that we
have defined $s_n=|\d|$ where $R_n(0) \in C^\d_n$, so that
$R_{n+1}(0)=R_n^{s_n+1}(0)$.
Let $c_n=|I_{n+1}|/|I_n|$.

We define the following convenient notation
\begin{align}
&I^X_n=\bigcup_{j \in X} I^j_n, &&I(X,n)=\frac {|I^X_n|} {|I_n|}=\sum_{j \in
X} \frac {|I^j_n|} {|I_n|}, && X \subset \Z,\\
&I^X_n=\bigcup_{\d \in X} I^\d_n, &&I(X,n)=\sum_{\d \in X} \frac {|I^\d_n|}
{|I_n|}, &&X \subset \Omega,\\
&C^X_n=\bigcup_{\d \in X} C^\d_n, &&C(X,n)=\frac {|C^X_n|} {|I_n|}=\sum_{\d
\in X} \frac {|C^\d_n|} {|I_n|}, &&X \subset \Omega.
\end{align}
(Thus $I^X_n$ and $I(X,n)$ are defined both for $X \subset \Z$ and $X
\subset \Omega$.)

The following summarizes Lemma~4.3,
Corollaries 6.8 and 6.10, and Remark~6.3 of \cite {AM1}.

\begin{lemma} \label {vn+1}

Almost every non-regular map satisfies
\be
\lim \frac {\ln v_{n+1}} {\ln c_n^{-1}}=\lim \frac {\ln s_n} {\ln c_n^{-1}}=
\lim \frac {\ln \ln c_{n+1}^{-1}} {\ln c_n^{-1}}=\lim \frac {\ln
r_n(\tau_n)} {\ln c_{n-1}^{-1}}=1,
\ee

\end{lemma}

In particular, $c_n$ decays very fast (this type of decay is called
torrential).

\subsubsection{Distortion estimates} \label {distortion}

Let us now discuss some estimates on the position of the critical value of
the return maps $R_n$,
which are relevant for distortion estimates.  The following summarizes
Lemmas 4.8 and 4.10 (and their proof) of \cite {AM1}.

\begin{lemma} \label {n1+delta}

For almost every non-regular map, for every $\delta>0$, for any $n$ big
enough, the following holds:

\begin{enumerate}

\item $|R_n(0)|>n^{-1-\delta} |I_n|$, and in particular, $R_n(0) \notin
\tilde I_{n+1}$,

\item The distance between $R_n(0)$ to $\partial I_n$ is at least
$n^{-1-\delta}|I_n|$,

\item For any $\d \in \Omega$, if $R_n(0) \notin C^\d_n$, then the distance
between $R_n(0)$ and $C^\d_n$ is at least $n^{-1-\delta} |C^\d_n|$,

\item For any $\d \in \Omega$, $\dist(R^\d_n) \leq n^{\frac {1}
{2}+\delta}$.

\end{enumerate}

\end{lemma}

The estimate above for distortion of branches $R^\d_n$ is relatively
pessimistic.  For most branches, we have much better bounds.  Indeed, if
$R_{n-1}(I^j_n) \subset
C^\d_{n-1}$ and $R_{n-1}(0) \notin I^\d_{n-1}$, then
$\dist (f|I^j_n)-1$ is at most of order of the quotient of $|I^j_n|$ by the
distance from $I^j_n$ to $0$ (this can be bounded from above by
$O(|C^\d_{n-1}|/|I^\d_{n-1}|)$ because $R_{n-1}(0) \notin I^\d_{n-1}$),
so $\dist (f|I^j_n)=1+O(c_{n-1})$.  Since
$R_n|I^j_n$ is the composition of $f|I^j_n$ and a diffeomorphism onto $I_n$
(which extends to $I_{n-1}$) with distortion bounded by $1+O(c_{n-1})$
(by the Koebe principle), we see that for all
those branches the distortion of $R_n$ is at most $1+O(c_{n-1})$.

Notice that for any $j$, both components of $I_n \setminus I^j_n$ have size
at least $|I^j_n| 2^n c_{n-1}^{-1/2}$.  Indeed, let
$R_{n-1}(I^j_n) \subset C^\d_{n-1}$.
Each connected component of $I_{n-1} \setminus C^\d_{n-1}$ must have size
at least of order $2^{4n} c_{n-1}^{-1} |C^\d_{n-1}|$ (which implies the desired
estimate), unless $|\d|=0$ (that is $C^\d_{n-1}=I_n$).
In this last case, the first item of the
previous lemma implies that each connected component of
$I_n \setminus I^j_n$ has size at least of order
$2^{-n} c_{n-1}^{-1} |I^j_n| \geq 2^n c_{n-1}^{-1/2} |I^j_n|$.

In particular, if $\dist (R_n|I^j_n)=1+O(c_{n-1})$ and the last entry of
$\d$ is $j$, we can also find better bounds for the distortion of
$R^\d_n$.  Indeed,
$R^\d_n$ is the composition of a map onto $I^j_n$ which extends to $I_n$,
and has distortion bounded by $1+o(c_{n-1}^{1/2})$ and $R_n|I^j_n$, so we
have $\dist(R_n^\d)=1+o(c_{n-1}^{1/2})$.

\subsubsection{Estimates on the capacity of some relevant sets}
\label {epsilon(gamma)}

In the course of proving the above estimates, one obtains several estimates
for the quasisymmetric capacities of certain sets, which will be important
here.  In order to be definite, let $\epsilon=\epsilon(\g)$ be the smallest
number such that, for $\kappa=1+\frac {\epsilon} {5}$ and for any
$\g$-qs map $h$ we have
\be
\frac {1} {\kappa} \left ( \frac {|J|} {|I|} \right )^\kappa
\leq \frac {|h(J)|}
{|h(I)|} \leq
\left ( \frac {\kappa|J|} {|I|} \right )^{1/\kappa},
\ee
so that $\epsilon(\g) \to 0$ as $\g \to 1$.

The following summarizes Corollaries 6.5 and 6.7 of \cite {AM1}.

\begin{lemma}

For almost every non-regular map, if $\epsilon_0=\epsilon(\g)<1/100$, then,
for $n$ large enough
\begin{align}
\label {r_n(j)}
&p_\g(r_n(x)>k c_n^{-4}|I_n) \leq e^{-k}, \quad k \geq 1,\\
\label {r_n(j) 1}
&p_\g(r_n(x)<c_{n-1}^{-1+2 \epsilon_0}|I_n) \leq
c_{n-1}^{\epsilon_0/10},\\
&p_\g(r_n(x)>c_{n-1}^{-1+2 \epsilon_0}|I_n) \leq
e^{-c_{n-1}^{-\epsilon_0/5}}.
\end{align}

\end{lemma}


\section{The critical orbit is typical} \label {typical}

\subsection{Outline} \label {thmB}
Let us summarize the main steps in the proof of Theorem~B.

\noindent(1)\,
We must show that (with total probability)
the proportion of time the critical
orbit spends in any given interval $T \subset I$
is given by $\mu(T)$.  It is of course enough to consider
a countable class of intervals which generates
all Borelians, and then prove
the distribution result (with total probability) for each interval in the
class.  Our choice of intervals will be domains $\xi$
of the first landing map from $I$ to $I_{n_0}$ (for arbitrary $n_0$).

\noindent(2)\,
We must be able to estimate $\mu(\xi)$ in terms of return branches.
Let $\psi^\xi_n(x)$ be the frequency of visits to $\xi$
of the iterates of a point $x \in I_n$ before $x$ returns to $I_n$
($\psi^\xi_n(x)$ only depends on the branch $I^j_n$ containing $x$).
We show that $\psi^\xi_n$ is concentrated around $\mu(\xi)$ and
indeed we show that $\mu(\xi)$ is the unique number $q$ such that, for every
$\epsilon>0$, we have $\lim_{n \to \infty}
p(|\psi^\xi_n(x)-q|>\epsilon|I_n)=0$.

\noindent(3)\,
We use an explicit Large Deviation Estimate to obtain a quantitative
estimate on the rate of decay of $p(|\psi^\xi_n(x)-\mu(\xi)|>\epsilon|I_n)$
(in $n$) using only the fact that it decays to $0$.  We obtain a
torrential estimate
($p(|\psi^\xi_n(x)-\mu(\xi)|>\epsilon|I_n)<c_{n-1}^{1/20}$).

\noindent(4)\,
We would like to show that returns $R_n(0)$ of
the critical point belong to branches of
$R_n$ with ``close to correct'' distribution on $\xi$, that is
$|\psi^\xi_n(R_n(0))-\mu(\xi)|<\epsilon$.
The previous estimate indicate that this should be the case, but the
Phase-Parameter relation is just quasisymmetric.  We show that the
torrential rate of decay still holds if instead of probabilities
$p(|\psi^\xi_n(x)-\mu(\xi)|>\epsilon|I_n)$ we consider qs-capacities
$p_{\g(n)}(|\psi^\xi_n(x)-\mu(\xi)|>\epsilon|I_n)$, provided we choose
$\gamma(n)$ very close to $1$.  This argument does not give any
reasonable bound on the rate of decay of
$\gamma(n)$ to $1$, it could be very fast.

\noindent(5)\,
We want to show that we may actually take $\gamma(n)$ as a constant $\g$
bigger than $1$.  For this we argue that a torrentially
small set of branches (in the $\g(n_s)$-qs)
of a fixed level $n_s$ has torrentially small effect
(in the $\g$-qs sense for some fixed $1<\g<\g(n_s)$) with respect to total
(and partial) time of branches in the subsequent levels.
This argument follows the proof of the
Collet-Eckmann condition in \cite {AM1}, where we used those ideas to
control the propagation of weakly hyperbolic branches.
A little bit of change is needed in order
to avoid a loss of the quasisymmetric
constant of level $n_s$, on which we do not have control.  For this reason,
we will work with modified quasisymmetric capacities in some arguments.

\noindent(6)\,
As a consequence, we see that except for a set with torrentially small
$\g$-qs capacity, return
branches of level $n$ are ``very good''
in the sense that they spend most of their time following
branches of level $n_s$ which satisfy $|\psi^\xi_{n_s}-\mu(\xi)|<\epsilon$.
As a consequence, those ``very good'' return branches of level $n$ satisfy
$|\psi^\xi_n-\mu(\xi)|<2\epsilon$.  As a bonus from the previous item
we get for free the estimates for
intermediate moments (not just full returns), which are needed also in the
proof of the Collet-Eckmann condition.

\noindent(7)\,
Using the Phase-Parameter relation we make the critical point falls in
``very good'' branches.  Thus the distribution of
the critical orbit on $\xi$ is $2\epsilon$ close to $\mu(\xi)$.
Making $\epsilon$ goes to $0$ we obtain Theorem B.

\subsection{Inductive estimates} \label {induc}

In this section we will show that a small (in the quasisymmetric sense)
set of branches of level $n_0$ has a small effect on most (in the
quasisymmetric sense) branches of level $n \geq n_0$.  This kind of argument
was already needed in the analysis of \cite {AM1}, so we will keep a similar
notation to that work, and will refer to it for some computations.

\subsubsection{Modified capacities}

For our application, we will need a modification
of the $\g$-qs capacities used by \cite {AM1}.  This is not the same
modification used by \cite {Av}.

We say that $h$ is a $(\g,C)$-homeomorphism if $h=h_2 \circ h_1$ where $h_2$
is $\g$-qs and $h_1$ is $C^1$ with distortion bounded by $C$.

If $X \subset I$ is a Borelian set, we let
\be
p_{\g,C}(X|I)=\sup \frac {|h(X \cap I)|} {|h(I)|}
\ee
where $h$ ranges over all $(\g,C)$-homeomorphisms.

Through the end of this section we will fix $\epsilon_0$ very small (say,
$1/1000$), but we won't need to make $\epsilon_0 \to 0$ later on.  Choose
$\hat \g$ very close to $1$ so that $\epsilon(\hat \g) \leq \epsilon_0$, in
the notation of \S \ref {epsilon(gamma)}
\comm{
if $h$ is a $(2\hat \g-1)$-qs
homeomorphism then for $k=1+\epsilon_0/5$ we have
\be
\frac {1} {k} \left ( \frac {|J|} {|I|} \right )^k \leq \frac {|h(J)|}
{|h(I)|} \leq
\left ( \frac {k|J|} {|I|} \right )^{1/k}.
\ee
}

Let us fix $C$ and $\g_0$ close to $1$
so that for $n$ big, any $(C \frac {n+1} {n},\g_0)$ homeomorphism is a
$\hat \g$-qs homeomorphism.
Let $C_n=C \cdot \frac {n+1} {n}$, $\tilde C_n=C \frac {2n+3} {2n+1}$.

In what follows, we will work with some fixed
$1 \leq \g \leq \g_0$, {\it but the estimates will be uniform for $\g$
in this range}, and with the sequences $C_n$ and $\tilde C_n$.

We will use $(\g,C_n)$ capacities to estimate the size of sets of return
branches of level $n$ and $(\g,\tilde C_n)$ for sets of landing branches of
level $n$.

The introduction of those constants is motivated by the following
result which can be proved using the methods of \cite {AM1}.

\begin{lemma}[Analogous to Remarks 5.1 and 5.2 of \cite {AM1}]

With total probability, there exists $n_0$ such that for $n>n_0$ and for all
$1 \leq \g \leq \g_0$, the following holds.
If $X \subset I_n$ then
\be \label {estind}
p_{\tilde C_n,\g}((R^\d_n)^{-1}(X)|I^\d_n) \leq 2^n p_{C_n,\g}(X|I_n).
\ee
And if $X \subset I_n$ and
\be
p_{\tilde C_n,\g}(X|I_n) \leq \delta \leq 2^{-n^2}.
\ee
then
\be \label {invimag}
p_{C_{n+1},\g}((R_n|I^0_n)^{-1}(X)|I_{n+1}) \leq \delta^{1/5}.
\ee

\end{lemma}

Induction applied to (\ref {estind}) gives:

\begin{lemma}[Analogous to Lemma 5.4 of \cite {AM1}]

With total probability, there exists $n_0$ such that for $n>n_0$ and all $1
\leq \g \leq \g_0$ the following holds.
Let $Q_n \subset \Z$ and let
$Q_n(m,r)$ be the set of all $\d$ with length $m$ and at least $r$ entries
on $Q_n$.  Let
\begin{align}
q_n&=p_{\g,C_n}(I^{Q_n}_n|I_n),\\
q_n(m,r)&=p_{\g,\tilde C_n}(I^{Q_n(m,r)}_n|I_n).
\end{align}
Then
\begin{equation}
q_n(m,r) \leq \binom {m} {r} (2^n q_n)^r.
\end{equation}

More generally, for any fixed $\d$, defining
\be
q^\d_n(m,r)=p_{\g,\tilde C_n}((R_n^\d)^{-1}(I^{Q_n(m,r)}_n)|I^\d_n),
\ee
we have
\begin{equation}
q^\d_n(m,r) \leq \binom {m} {r} (2^n q_n)^r.
\end{equation}

\end{lemma}

This estimate will be mainly used to estimate $q_n(m,r)$ for $m$ large and
$\frac {r} {m}$ larger than $(6 \cdot 2^n) q_n$.  Notice that if
$q^{-1} \geq 6 \cdot 2^n$ and $q \geq q_n$ then by Stirling formula,
\be \label {ld4}
q_n(m,(6 \cdot 2^n) q m) \leq 2^{-(6 \cdot 2^n) q m},
\ee
and
\be \label {ld5}
\sum_{k \geq q^{-2}} q_n(k,(6 \cdot 2^n) q k) \leq
2^{-n} q^{-1} 2^{-(6 \cdot 2^n) q^{-1}}.
\ee

\comm{
\subsubsection{Preliminary estimates on frequency of times}

The following estimates (valid with total probability, for $n$ big)
were obtained in \cite {AM1}.
\begin{align}
&p_{\hat \g}(r_n(x)>k c_n^{-4}|I_n) \leq e^{-k}, \quad k>1\\
&p_{\hat \g}(r_n(x)<c_{n-1}^{-1+2 \epsilon_0}|I_n) \leq
c_{n-1}^{\epsilon_0/10}\\
&p_{\hat \g}(r_n(x)>c_{n-1}^{-1+2 \epsilon_0}|I_n) \leq
e^{-c_{n-1}^{-\epsilon_0/5}}.
\end{align}

Notice that those estimates remain valid if instead of $\hat \g$ capacities
we consider $(\g,C_n)$ or $(\g,\tilde C_n)$ capacities with $1 \leq \g \leq
\g_0$, due to our choice of $\g_0$ and $C$.
}

\subsubsection{Estimates on time} \label {esttime}

Following \cite {AM1}, we define the set of standard landings at time $n$,
$LS(n) \subset \Omega$ as the set
of all $\d=(j_1,...,j_m)$ satisfying the following.
\begin{align}
\tag{LS1} &c^{-1/2}_n<m<c^{-1-2\epsilon_0}_n,\\
\tag{LS2} &r_n(j_i)<c^{-14}_{n-1}, \quad \text {for all } i,\\
\tag{LS3} &\#\{1 \leq i \leq k,\, r_n(j_i)<c^{-1+2\epsilon_0}_{n-1}\}<
(6 \cdot 2^n) c^{\epsilon_0/10}_{n-1} k, \quad \text {for }
c^{-2}_{n-1} \leq k \leq m,\\
\tag{LS4} &\#\{1 \leq i \leq k,\, r_n(j_i)>c^{-1-2\epsilon_0}_{n-1}\}<
(6 \cdot 2^n) e^{-c^{-\epsilon_0/5}_{n-1}} k, \quad
\text {for } c^{-1/n}_n \leq k \leq m.
\end{align}

The following estimate was obtained in \cite {AM1}, and is a consequence of
the estimates of \S \ref {epsilon(gamma)}:

\begin{lemma}[Analogous to Lemma 7.1 of \cite {AM1}] \label {standard}

Let $LS(n)$ denote the set of standard landings.  Then
\be \label {sta1}
p_{\hat \g}(C^{\Omega \setminus LS(n)}_n|I_n)<c_n^{1/3}.
\ee
\be \label {sta2}
p_{\hat \g}(C^{\Omega \setminus LS(n)}_n|I^{\tau_n}_n)<c_n^{1/3}.
\ee

\end{lemma}

Let $T \subset \Z$ be given.  Let us define $VG(T,n_0,n) \subset \Z$
and $LE(T,n_0,n) \subset \Omega$
inductively as follows.  Let $VG(T_{n_0},n)=\Z
\setminus T$. 
Assuming $VG(T,n_0,n)$ defined, let $LE(T,n_0,n)$ be
the set of all $\d \in LS(n)$ such that $\d=(j_1,...,j_m)$ and
\be
\tag{LE}
\#\{j_i \notin VG(T,n_0,n), 1 \leq i \leq k\}<(6 \cdot 2^n) c_n^{1/20} k,
\quad \text {for } c_{n-1}^{-2} \leq k \leq m.
\ee

Let $VG(T,n_0,n+1)$ be the set of all $j$ such that $R_n(I^j_{n+1})
\subset LE(T,n_0,n)$.

In what follows, we will put the condition that $T$ is a small set of
branches of some (deep) level $n_0$ in the sense that
\begin{equation} \label {small}
p_{\g,C_{n_0}}(I^T_{n_0}|I_{n_0})<c_{n_0-1}^{1/20}
\end{equation}
for some $n_0$ and some $1 \leq \g \leq \g_0$.

The definition of the class $VG(T,n_0,n)$ is designed so that such branches
do not pass very often by $T$ before returning.  The precise constants in
the definition were chosen so that they allow to show that
$VG(T,n_0,n)$ correspond to most branches of level $n$ (by induction).
Those two estimates are given below:

\begin{lemma}[see also Lemma 7.2 of \cite {AM1}] \label {ier 1}

With total probability, for all $n_0$ sufficiently big,
if $T$ satisfies (\ref {small})
for some $1 \leq \g \leq \g_0$ then for all $n \geq n_0$, we have
\begin{equation} \label {f}
p_{\g,\tilde C_n}(C^{\Omega \setminus LE(T,n_0,n)}_n|I_n)<c_n^{2/7}
\end{equation}
\begin{equation} \label {s}
p_{\g,C_n}(I^{\Z \setminus VG(T,n_0,n)}_n|I_n)<c_{n-1}^{1/20}.
\end{equation}
Furthermore,
\be \label {stau}
p_{\g,\tilde C_n}(C^{\Omega \setminus
LE(T,n_0,n)}_n|I^{\tau_n}_n)<c_n^{2/7}.
\ee
\end{lemma}

\begin{pf}

If (\ref {f}) is valid for $n$ then by (\ref {invimag}) we get
\be
p_{\g,C_{n+1}}(I^{\Z \setminus VG(T,n_0,n+1)}_{n+1}|I_{n+1})<
c_n^{2/35}<c_n^{1/20}
\ee
which gives (\ref {s}) for $n+1$.

Let us assume the validity of (\ref {s}) for $n$.  Then the
$(\tilde C_n,\g)$-capacity of the set of
standard landings which fail to satisfy LE is much less than $c_n$, by (\ref
{ld4}).  Using Lemma \ref {standard} we get
\be
p_{\g,\tilde C_n}(C^{\Omega \setminus LE(T,n_0,n)}_n|I_n)<c_n^{1/3}+c_n
\leq c_n^{2/7}.
\ee
This implies that (\ref {f}) is valid for $n$.  A similar computation gives
(\ref {stau}) for $n$.

Since (\ref {s}) is valid for $n_0$ by hypothesis, we get
(\ref {f}), (\ref {s}) and (\ref {stau}) for all $n$ by induction.
\end{pf}

\begin{lemma}[Analogous to Lemma 7.6 of \cite {AM1}] \label {ret time}

With total probability, for all $n_0$ big enough and for all $n \geq n_0$,
the following holds.

Let $j \in VG(T,n_0,n+1)$, and let $\d$ be such that
$R_n(I^j_{n+1}) \subset C^\d_n$ and
$\d=(j_1,...,j_m)$.  Let $c_n^{-2/n}<k \leq r_{n+1}(j)$.
Let $m_k$ be biggest possible with
\be
v_n+\sum_{j=1}^{m_k} r_n(j_i) \leq k
\ee
\be
\beta_k=\sum_{\ntop {1 \leq i \leq m_k,}
{j_i \in VG(T,n_0,n)}} r_n(j_i).
\ee
Then $1-\frac {\beta_k} {k}<c_{n-1}^{1/100}$.

\end{lemma}

\begin{lemma} \label {ier}

With total probability, for all $n_0$ big enough and for all $n \geq n_0$,
the following holds.
Let $j \in VG(T,n_0,n+1)$ and $x \in I^j_{n+1}$, and let
$c_n^{-2/n} \leq k \leq
r_{n+1}(x)$.  Then
\be
\sum_{\ntop {i<k,} {f^i(x) \in I^T_{n_0}}}
r_{n_0}(f^i(x))<c_{n_0-1}^{1/200} k,
\ee

\end{lemma}

\begin{pf}

Let $\alpha_n=\sum_{k=n_0}^{n-1} c_{k-1}^{1/110}<c_{n_0-1}^{1/200}$.
We show by induction that if
\be \label {sum}
\sum_{\ntop {i<r_n(x),} {f^i(x) \in I^T_{n_0}}}
r_{n_0}(f^i(x)) \leq \alpha_n r_n(x), \quad \text {for all }
x \in I^{VG(T,n_0,n)}_n,
\ee
then
\be
\sum_{\ntop {i<k,} {f^i(x) \in I^T_{n_0}}}
r_{n_0}(f^i(x))<\alpha_{n+1} k, \quad \text {for all }
x \in I^{VG(T,n_0,n+1)}_{n+1}, c_n^{-2/n}
\leq k \leq r_{n+1}(x).
\ee

Indeed (using the notation of Lemma \ref {ret time}),
\be
\sum_{\ntop {i<k,} {f^i(x) \in I^T_{n_0}}}
r_{n_0}(f^i(x)) \leq k-\beta_k+\alpha_n \beta_k+c_{n-1}^{-14}
\leq \left (1-\frac {\beta_k} {k}+\alpha_n+c_{n-1}^{-14} c_n^{2/n}
\right ) k \leq \alpha_{n+1} k.
\ee
This gives our result by induction, since for $n=n_0$, the left side of
(\ref {sum}) is $0$. 
\end{pf}

\subsubsection{Control of intermediate times}

Let us define $LC(T,n_0,n) \subset \Omega$, $n_0,n \in \N$,
$n \geq n_0$ as the set of all $\d=(j_1,...,j_m)$ in $LE(T,n_0,n)$
satisfying
\begin{align}
\tag {LC1} &j_i \in VG(T,n_0,n), \quad 1 \leq i \leq c_{n-1}^{-1/30},\\
\tag {LC2} &\#\{1 \leq i \leq k,\, r_n(j_i)<c^{-1+2\epsilon_0}_{n-1}\}<
(6 \cdot 2^n) c^{\epsilon_0/10}_{n-1} k, \quad
\text {for } c^{-\epsilon_0/5}_{n-1} \leq k \leq m,\\
\tag {LC3} &\#\{1 \leq i \leq k,\, j_i \notin VG(n_0,n)\}<
(6 \cdot 2^n) c^{1/60}_{n-1} k, \quad
\text {for } c^{-1/30}_{n-1} \leq k \leq m,\\
\tag {LC4} &\#\{1 \leq i \leq k,\, r_n(j_i)>c^{-1-2\epsilon_0}_{n-1}\}<
(6 \cdot 2^n) c^{100}_{n-1} k, \quad
\text {for } c_{n-1}^{-200} \leq k \leq m,\\
\tag {LC5} &r_n(j_i)<c_{n-1}^{-1-2\epsilon_0}, \quad
1 \leq i \leq e^{c_{n-1}^{-\epsilon_0/5}/2}.
\end{align}
\comm{
Notice that LC4 and LC5 overlap, since
$c_{n-1}^{-200}<e^{c_{n-1}^{-\epsilon_0/5}/2}$ as do LC1 and LC3,
from this we
can conclude that we can control the proportion of large times or non very
good times in all moments (and not only for not too small initial segments).
}

\begin{lemma}[Analogous to Lemma 7.7 of \cite {AM1}] \label {ier 2}

For all $n_0$ sufficiently big and all $n \geq n_0$, if $T$ satisfy (\ref
{small}), then
\be
p_{\g, \tilde C_n}(C^{\Omega \setminus LC(T,n_0,n)}_n|I_n)<c_{n-1}^{1/100}
\ee
and if $\tau_n \in VG(T,n_0,n)$,
\be
p_{\g, \tilde C_n}(C^{\Omega \setminus
LC(T,n_0,n)}_n|I^{\tau_n}_n)<c_{n-1}^{1/100}.
\ee

\end{lemma}

\begin{lemma} \label {ier 3}

For all $n_0$ is sufficiently big, for all $n \geq n_0+1$,
for all $T$, if $\d \in LC(T,n_0,n)$, then for all
$c_{n-1}^{-4/(n-1)}<k \leq l_n(\d)$, and for all $x \in C^\d_n$,
\be
\sum_{\ntop {i \leq k,} {f^i(x) \in I^T_{n_0}}}
r_{n_0}(f^i(x))<2 c_{n_0-1}^{1/200} k
\ee

\end{lemma}

\begin{pf}

The proof follows closely the argument of Lemma 7.11 of
\cite {AM1}, but since the claim is formally different, we will
repeat some steps here,
referring to the computations in \cite {AM1}.

Let $\d=(j_1,...,j_m)$.  Assume that $k \leq r_n(j_1)$.  Since $j_1 \in
VG(T,n_0,n)$, we get the result as a consequence of Lemma \ref {ier}.
This will still work if we take $k \leq r_n(j_1)+...+r_n(j_t)$,
where $j_i$
is very good for $1 \leq i \leq t$.

Let $m_k$ be the last return completed before $k$, that is
$\sum_{i=1}^{m_k} r_n(j_i) \leq k$.
We must analyze the case where $j_i$ is not very good
for some $i \leq m_k+1$. 
In this case, we must have, by LC1, $m_k \geq c_{n-1}^{-1/30}$.  Let
\be
\beta_k=\sum_{\ntop {i \leq m_k,} {j_i \in VG(T,n_0,n)}} r_n(j_i).
\ee
After some computations, we get (see \cite {AM1})
\be
\sum_{\ntop {i \leq m_k,} {j_i \notin VG(T,n_0,n)}} r_n(j_i) \leq 4
c_{n-1}^{1/35} k,
\ee
and
\be
r_n(j_{m_k+1}) \leq c_{n-1}^{1/80} k,
\ee
(if $m_k=|\d|$, that is, $k=l_n(\d)$, we will make the convention that
$r_n(j_{m_k+1})=0$).  We obtain
\be
\sum_{\ntop {i \leq k,} {f^i(x) \in I^T_{n_0}}}
r_{n_0}(f^i(x)) \leq
c_{n_0-1}^{1/200} \beta_k+r_n(j_{m_k+1})+
\sum_{\ntop {i \leq m_k,} {j_i \notin VG(T,n_0,n)}} r_n(j_i)
\leq 2 c_{n_0-1}^{1/200} k.
\ee
\end{pf}

\subsection{Large deviation estimate}

\comm{
We define the following convenient notation
\begin{align}
&I^X_n=\bigcup_{j \in X} I^j_n, &&I(X,n)=\frac {|I^X_n|} {|I_n|},
&&&X \subset \Z,\\
&I^X_n=\bigcup_{\d \in X} I^\d_n, &&I(X,n)=\frac {|I^X_n|} {|I_n|},
&&&X \subset \Omega,\\
&C^X_n=\bigcup_{\d \in X} C^\d_n, &&C(X,n)=\frac {|C^X_n|} {|I_n|},
&&&X \subset \Omega.
\end{align}
(Thus $I^X_n$ and $I(X,n)$ are defined both for $X \subset \Z$ and $X
\subset \Omega$.)
}

\comm{
Given a set $X \subset \Z$, we let $I^X_n$,
denote the union of $I^j_n$, $j \in X$, and $I(X,n)$ the density of
$I^X_n$ in $I_n$.  We define analogously $I^X_n$ and $I(X,n)$, for
$X \subset \Omega$.  We also define $C^X_n$,
$C(X,n)$ in the natural way.
}

\subsubsection{More precise estimates on return times}

We will need several times the following elementary result.

\begin{lemma} \label {binomial}

Let $m>100$.  If $q \in [0,1]$ and $\epsilon \geq m^{-1/4}$ then
\be
\sum_{|\frac {k} {m}-q|>\epsilon} \binom {m} {k} q^k (1-q)^{m-k} \leq
e^{-m^{1/4}}.
\ee

\end{lemma}

\begin{pf}

Let $y_k=\binom {m} {k} q^k (1-q)^k$, and recall that $y_k \leq 1$ for all
$k$ (since $\sum y_k=1$).  It is enough to show that
$\sum_{k/m>q+\epsilon} y_k \leq e^{-m^{1/4}}/2$,
since the estimate
corresponding to $\frac {k} {m}<q-\epsilon$ reduces to this
one after interchanging $q$ and $1-q$.  Let
$x_k=\frac {y_{k+1}} {y_k}=\frac
{1-k/m} {(k+1)/m} \cdot \frac {q} {1-q}$.
If $\frac {k} {m} \geq
q+\frac {\epsilon} {2}$ then $x_k<\frac {1-q-\epsilon/2}
{q+\epsilon/2} \cdot \frac {q} {1-q}<1-\frac {\epsilon} {2}$.
Notice that if $k$ is minimal with $\frac {k}
{m}>q-\epsilon$ then there are about $\frac {\epsilon m} {2}$ integers $j<k$
such that $q+\frac {\epsilon} {2} \leq \frac {j} {m}$.  We conclude that
\be
\sum_{\frac {k} {m}>q+\epsilon} y_k \leq
\sum_{i \geq
\frac {\epsilon m} {2}} \left (1-\frac {\epsilon} {2} \right )^i
\leq \frac {2} {\epsilon} e^{-m \epsilon^2/4} \leq \frac {e^{m^{1/4}}} {2},
\ee
and the result follows.
\end{pf}

{\it Notation warning:}
In what follows, we will work with a fixed typical unimodal map $f$.  We
will use $\delta^{(n)}_1,...,\delta^{(n)}_{14}$ to denote several small
constants (going to $0$ with $n$).  We shall always choose
$\delta^{(n)}_{i+1}$ after fixing $\delta^{(n)}_i$, and satisfying
(among other requirements) $\delta^{(n)}_{i+1} \geq 10 \delta^{(n)}_i$.  We
shall also take $\delta^{(n)}_1>n^{-1}$.

\bigskip
\paragraph{} \label {B_n}

Let $\d_{n-1}$
be such that $R_{n-1}(0) \in C^{\d_{n-1}}_{n-1}$, and
let $B_n$ be the set of all $j$ such that
$R_{n-1}(I^j_n)=C^\d_{n-1}$, where $|\d|<|\d_{n-1}|$ and
$\d=(j_1,...,j_k)$ is obtained by considering the first $k$ entries of
$\d_{n-1}$.  Let $A_n=\Z \setminus (\{0\} \cup B_n)$.
Recall the estimates of \S \ref {distortion}.
One readily sees that
$I(B_n,n) \leq 2^{-n} c_{n-1}^{1/2}$ and for
$j \in B_n$, $r_n(j) \leq v_n$.
Notice that if $j \in A_n$, the interval $I^j_n$ is far from the critical
point in the sense that $c_{n-1}^{-1}|I^j_n|$ is much smaller than
the distance from $I^j_n$ to $0$.  It follows that, for any
$\d$ whose last entry belongs to $A_n$, $\dist (R_n^\d)<1+\delta_1^{(n)}
c_{n-1}^{1/2}$.  If
the last entry belongs to $B_n$ we will use the general estimate $\dist
(R_n^\d) \leq n^{2/3}$.

\bigskip
\paragraph{}

Let $m(\d)$ be the number of entries of $\d$ which belong to
$A_n$.  The following
easy estimates follow from the previous discussion by induction:
\begin{align}
(1-(1+\delta^{(n)}_2 c_{n-1}^{1/2})c_n)^m &\leq I(\{m(\d)=m\},n) &&\leq
(1+\delta^{(n)}_2 c^{1/2}_{n-1})
(1-(1-\delta^{(n)}_2 c_{n-1}^{1/2})c_n)^m,\\
\label {m(d)}
c_n (1-(1+\delta^{(n)}_2 c_{n-1}^{1/2})c_n)^m &\leq C(\{m(\d)=m\},n) &&\leq
(1+\delta^{(n)}_2 c^{1/2}_{n-1})
c_n (1-(1-\delta^{(n)}_2 c_{n-1}^{1/2})c_n)^m.
\end{align}

Let $Q(m',m) \subset \Omega$
be the set of $\d$ of size $m'$ and with at least
$m$ entries on $B_n$, that is, $Q(m',m)=\{\d \in \Omega,\, |\d|=m',\,
|\d|-m(\d)| \geq m\}$.
Let $q(m',m)=I(Q(m',m),n)$.  From the definition and
the estimates on distortion we have
\begin{align}
q(m',m) &\leq q(m'-1,m)+\delta^{(n)}_3 c_{n-1}^{1/2}
(q(m'-1,m-1)-q(m'-1,m))\\
\nonumber
&=(1-\delta^{(n)}_3 c_{n-1}^{1/2}) q(m'-1,m)+\delta^{(n)}_3
c_{n-1}^{1/2} q(m'-1,m-1),
\end{align}
which implies by induction,
\be
q(m',m) \leq \sum_{k=m}^{m'}
\binom {m'} {k} (\delta^{(n)}_3 c_{n-1}^{1/2})^k
(1-\delta^{(n)}_3 c_{n-1}^{1/2})^{m'-k}.
\ee

Let us compute a few consequences of those estimates.  Let
$H$ be the set of all $\d \in \Omega$ such that {\it at least
one} of the following holds:
\begin{enumerate}
\item [(H1)] $|\d| \geq c_n^{-1/n}$ and $|\d|-m(\d) \geq 2
c_{n-1}^{1/2} |\d|$,
\item [(H2)] $\d$ has some entry $j_i$ satisfying
$r_n(j_i) \geq c_{n-1}^{-14}$.
\end{enumerate}
Using the present discussion to estimate (H1) and (\ref {r_n(j)}) to
estimate (H2) we get
\be \label {H}
C(H,n) \leq I(H,n) \leq e^{-c_n^{-1/(8n)}}+e^{-c_{n-1}^{-19/2}} \leq
e^{-c_{n-1}^{-9}}.
\ee

\comm{
Let $K(m')$ be the set of $\d$ with $|\d|=m'$ and $|\d|-m(\d)
\geq 2 c_{n-1}^{1/2} |\d|$.  Let $K=\cup_{m' \geq c_n^{-1/n}} K(m')$.  Then
\be \label {K}
C(K,n) \leq I(K,n) \leq e^{-c_n^{1/(8n)}}.
\ee
}
Let $V$ be the set of $\d$ with $m(\d) \leq c_n^{-1/n}$.  The present
discussion gives
\be \label {V}
C(V,n) \leq 2 c_n^{1-1/n}.
\ee
\comm{
Let $\hat V$ be the set of $\d$ with $m(\d) \leq c_n^{1/n}$ but with $|\d|
\geq 2 c_n^{1/n}$.  Then
\be
I(\hat V,n) \leq e^{-c_n^{1/(8n)}}.
\ee
}

\comm{
Using estimate (\ref {r_n(j)}) one also gets the following.
Let $H$ be the set of $\d$ which have some
entry with $r_n(j_i) \geq c_{n-1}^{-14}$.  Then
\be \label {H}
C(H,n) \leq I(H,n) \leq e^{-c_{n-1}^{-9}}.
\ee
}

\bigskip
\paragraph{}

We will also need the following easy estimate:

\begin{lemma} \label {bino}

Fix $P \subset A_n$, and let $p=I(P,n)$.
Let $P(m,r) \subset \Omega$ be the set of all $\d$ with $m(\d)=m$
and with exactly $r$ entries in $P$.
Let $\overline P(m,r)$ (respectively $\underline P (m,r)$)
denote the union of all $P(m,r)$ with $r' \geq r$
(respectively $r' \leq r$).
Let $p(m,r)=I(P(m,r),n)$, $\overline p(m,r)=I(\overline P (m,r),n)$
and $\underline p(m,r)=I(\underline P (m,r),n)$.

We have, with $\overline p=p (1+4\delta^{(n)}_1 c_{n-1}^{1/2})$ and
$\underline p=p (1-4\delta^{(n)}_1 c_{n-1}^{1/2})$
\begin{align}
\label {k1}
\overline p (m,r) &\leq (1-\overline p) \overline p (m-1,r)+
\overline p \cdot \overline p (m-1,r-1)\\
\label {k2}
\underline p (m,r) &\leq (1-\underline p) \underline p (m-1,r)+\underline p
\cdot \underline p (m-1,r-1)\\
\label {k3}
\overline p (m,r) &\leq (1+2\delta^{(n)}_2 c_{n-1}^{1/2})
\sum_{k=r}^m \binom {m} {k} \overline p^k
(1-\overline p)^{m-k}\\
\label {k4}
\underline p(m,r) &\leq (1+2\delta^{(n)}_2 c_{n-1}^{1/2})
\sum_{k=0}^r \binom {m} {k} \underline p^k
(1-\underline p)^{m-k}
\end{align}

\end{lemma}

\begin{pf}

We notice that $p(1,0) \leq 1-p$, $p \leq p(1,1) \leq
p (1+2\delta^{(n)}_1 c_{n-1}^{1/2})$.
Let us consider a connected component $E$
of $I^{\overline P (m,r)}_n$.  It is either
contained in a connected component of $I^{\overline P(m-1,r)}_n$ or it is
contained in a component $\hat E$ of $I^{P(m-1,r-1)}_n$.  In this last case,
the iterate of $R_n$ which takes $\hat E$ to $I_n$
(necessarily with distortion bounded by $1+\delta^{(n)}_1 c_{n-1}^{1/2}$)
must take $E$ to a component of $I^{P(1,1)}_n$.  It follows that
\be
\overline p(m,r) \leq \overline p (m-1,r)+(1+\delta^{(n)}_1 c_{n-1}^{1/2})
p(m-1,r-1) p(1,1).
\ee
Since $p(m-1,r-1)=\overline p(m-1,r-1)-\overline p(m-1,r)$, we get
(\ref {k1}), and (\ref {k3}) follows by induction.

Let us now consider a connected component $E$
of $I^{\underline P (m,r)}_n$.  It is either
contained in a connected component of $I^{\underline P (m-1,r-1)}_n$ or it is
contained in a component $\hat E$ of $I^{P(m-1,r)}_n$.  In this last case,
the iterate of $R_n$ which takes $\hat E$ to $I_n$
(necessarily with distortion bounded by $1+\delta^{(n)}_1 c_{n-1}^{1/2}$)
must take $E$ to a component of $I_n \setminus I^{P(1,1)}_n$.
It follows that
\be
\underline p(m,r) \leq \underline p(m-1,r-1)+p(m-1,r)
(1-(1-\delta^{(n)}_1 c_{n-1}^{1/2})p(1,1)).
\ee
Since $p(m-1,r)=\underline p (m-1,r)-\underline p (m-1,r-1)$, we get
(\ref {k2}), so (\ref {k4}) follows by induction.
\end{pf}

\subsubsection{Return times}

Let us fix $\Theta \subset \Z \setminus \{0\}$, $\theta=I(\Theta,n)$.
We would like to estimate
\be
\zeta=\sum_{j \in \Theta} r_n(j) I(j,n)
\ee
in terms of $\theta$ (specially for the case $\Theta=\Z \setminus \{0\}$).
In order to do so, it is convenient to write
$\zeta=\zeta^A+\zeta^B$, where
\be
\zeta^A=\sum_{j \in \Theta \cap A_n} r_n(j) I(j,n), \quad
\zeta^B=\sum_{j \in \Theta \cap B_n} r_n(j) I(j,n).
\ee
Notice that it is easy to estimate (using \S \ref {B_n})
\be
\zeta^B \leq v_n I(\Theta \cap B_n,n) \leq c_{n-1}^{-1-\delta^{(n)}_4}
\min \{\theta,2^{-n} c^{1/2}_{n-1}\}.
\ee
To estimate $\zeta^A$, we will consider
the level sets $M_s=\{j \in A_n \cap \Theta|\, r_n(j)=s\}$, so that
$\zeta^A=\sum s m_s$, where $m_s=I(M_s,n)$.
Let $L=\{s|m_s \geq c_n^{1/(8n)}\}$,
$S=\{s|m_s<c_n^{1/(8n)}\}$.  Define
\be
\zeta^L=\sum_{s \in L} s m_s, \quad \zeta^S=\sum_{s \in S} s m_s,
\ee
so that $\zeta^A=\zeta^L+\zeta^S$.
Notice that by (\ref {r_n(j)}),
\be
\zeta^S=\sum_{\ntop {s \in S,} {s \leq c_n^{-1/(32n)}}} s m_s+
\sum_{\ntop {s \in S,} {s>c_n^{-1/(32n)}}} s m_s
\leq c_n^{1/(16n)}+\sum_{t \geq c_n^{-1/(32n)}}
t e^{-c_{n-1}^{-4} t} \leq c_n^{1/(32n)}.
\ee

\bigskip
\paragraph{} \label {Dalpha}

Let $N$ be the set of all $\d \in \Omega$ such that $m(\d) \geq
c_n^{-1/n}$ and {\it at least one} of the following holds:
\begin{enumerate}
\item [(N1)] For some $s \in L$, the number $u$
of entries $j_i$ of $\d$ belonging to $M_s$ satisfies either
$\frac {u} {m}>(1+4\delta^{(n)}_1 c_{n-1}^{1/2}) m_s+c_n^{1/(8n)}$, or
$\frac {u} {m}<(1-4\delta^{(n)}_1 c_{n-1}^{1/2}) m_s-c_n^{1/(8n)}$,
\item [(N2)] For some $s \in S$, the number $u$
of entries $j_i$ of $\d$ belonging to $M_s$ satisfies $\frac {u} {m} \geq
2c_n^{1/(8n)}$.
\end{enumerate}
It follows from Lemmas \ref {binomial} and \ref {bino} that
\be \label {N}
I(N,n) \leq 2 e^{-c_n^{-1/(10n)}}.
\ee

\comm{
Let $N(m)$ be the set of $\d$ with $m(\d)=m$, and, for every $s \in L$,
the number of entries in $M_s$, divided by $m$,
is outside $[(1-\delta c_{n-1}^{1/2})m_s-c_n^{1/(8n)},
(1+\delta c_{n-1}^{1/2})m_s+c_n^{1/(8n)}]$.
Let $N=\cup_{m \geq c_n^{-1/n}} N(m)$.
We have $I(N,n) \leq e^{-c_n^{-1/(10 n)}}$, by Lemma \ref {binomial}.

Let $\hat N(m)$ be the set of $\d$ with $m(\d)=m$ and, for every
$s \in S$, the number of entries in
$M_s$, divided by $m$, is more than $2 c_n^{1/(8n)}$.
Let $\hat N=\cup_{m \geq c_n^{-1/n}} \hat N(m)$.
We have $I(\hat N,n) \leq e^{-c_n^{-1/(10 n)}}$ by Lemma \ref {binomial}.
}

Let $D=N \cup H \cup V$
and $\hat D=N \cup H$.
By (\ref {N}), (\ref {H}) and (\ref {V}) we have
\begin{align}
C(D,n)&\leq c_n^{1-2/n}\\
\label {hat D}
C(\hat D,n)&\leq e^{-c_{n-1}^{-17/2}}.
\end{align}

If $\d \notin D$, we have
\begin{align}
\label {D1}
\frac {1} {m(\d)} \sum_{j_i \in \Theta} r_n(j_i) &\geq
(1-\delta^{(n)}_5 c_{n-1}^{1/2}) \zeta^L,\\
\label {upper}
\frac {1} {m(\d)} \sum_{j_i \in \Theta} r_n(j_i)
&\leq (1+\delta^{(n)}_5 c_{n-1}^{1/2})
\zeta^L+(|\d|-m(\d))v_n+2 c_{n-1}^{-28} c_n^{1/(8n)}\\
\nonumber
&\leq (1+\delta^{(n)}_5 c_{n-1}^{1/2}) \zeta^L+c_{n-1}^{-\frac {1}
{2}-\delta^{(n)}_5}.
\end{align}
while, if $\d \notin \hat D$, we have either $\d \notin V$ (in which case
(\ref {D1}) and (\ref {upper}) hold)
or $\d \in V$ in which case we have
$l_n(\d) \leq c_n^{-2/n}$.

Notice that $D$, $\zeta$ and $\zeta^L$ depend on $\Theta$ (and on $n$).
If needed we will
stress this dependence by writing $D(\Theta)$, $\zeta(\Theta)$
and $\zeta^L(\Theta)$.

\bigskip
\paragraph{}

Let
\be
\alpha_n=\zeta(\Z \setminus \{0\})=\sum_{j \neq 0} r_n(j) I(j,n).
\ee
Notice that due to (\ref {r_n(j) 1}),
\be \label {alphandelta}
\alpha_n>c_{n-1}^{-1+\delta^{(n)}_6}.
\ee

\begin{lemma} \label {alphanlemma}

We have
\be \label {alpha_n}
\left |\frac {\alpha_n} {\alpha_{n-1} c_{n-1}^{-1}}-1 \right |
<c_{n-2}^{1/30},
\ee
and for any set $\Theta \subset \Z \setminus \{0\}$ with
$\theta=I(\Theta,n)$, we have
\be \label {G}
\zeta(\Theta)=\sum_{j \in \Theta} r_n(j) I(j,n) \leq
(3 \theta (1-\ln \theta)+c_{n-1}) \alpha_n
\ee

\end{lemma}

\begin{pf}

Letting $\Theta=\Z \setminus \{0\}$ and keeping the previous notation,
we have clearly
\be
\zeta^L \leq \alpha_n=\zeta^L+\zeta^B+\zeta^S \leq
\zeta^L+v_n I(B_n,n)+c_n^{1/(32n)} \leq \zeta^L+
c_{n-1}^{-\frac {1} {2}-\delta^{(n)}_7},
\ee
and since $\alpha_n \geq c_{n-1}^{-1+\delta^{(n)}_6}$
by (\ref {alphandelta}), we actually have
\be
1 \leq \frac {\alpha_n} {\zeta^L} \leq 1+c_{n-1}^{\frac {1}
{2}-\delta^{(n)}_8}.
\ee

The previous discussion in \S \ref {Dalpha} gives for $\d \notin D$,
\be \label {upper'}
(1-c_{n-1}^{\frac {1} {2}-\delta^{(n)}_9}) \alpha_n \leq
(1-\delta^{(n)}_5 c_{n-1}^{1/2})
\zeta^L \leq \frac {1} {m(\d)} \sum r_n(j_i)=\frac
{l_n(\d)} {m(\d)} \leq
(1+c_{n-1}^{\frac {1} {2}-\delta^{(n)}_9}) \alpha_n.
\ee

Using the estimate (\ref {m(d)})
on the distribution of $m(\d)$, we get
\be
(1-c_{n-1}^{\frac {1} {2}-\delta^{(n)}_{10}}) c_n^{-1} \alpha_n \leq
\sum_{\d \notin D} l_n(\d) C(\d,n)
\ee
which implies that for each $j \in A_n$ we have
\be
I(j,n) (1-c_{n-1}^{\frac {1} {2}-\delta^{(n)}_{11}})
c_n^{-1} \alpha_n \leq \sum_{C^\d_n \subset I^j_n} l_n(\d) C(\d,n).
\ee
Let us now consider the set $Z \subset A_n$ of all $j$ such that
$R_n(I^0_n)$ contains $I^j_n$, $r_n(j)<c_{n-1}^{-14}$, and
such that $R_n(0)$ is at least $c_{n-1}^{1/4} |I_n|$ away from $I^j_n$.
Let $\hat Z$ denote the set of $j \in \Z$ such that $R_n(I^j_{n+1}) \subset
I^Z_n$.  Then $I(\Z \setminus \hat Z,n+1)<c_{n-1}^{1/9}$.
Since $I(j,n) \leq \delta^{(n)}_{12} c_{n-1}^{1/2}$ for all
$j$, the distortion of
$(R_n|I^0_n)^{-1}$ restricted to any component of $I^Z_n$ is bounded by
$1+\delta^{(n)}_{14} c_{n-1}^{1/4}$.  We conclude
\be \label {below}
(1-c_{n-1}^{1/10}) \alpha_n c_n^{-1} \leq
\sum_{j \in \hat Z} r_{n+1}(j) I(j,n+1) \leq \alpha_{n+1}.
\ee

\comm{
Notice that if $\d \notin D$ then
\be
\frac {1} {m(\d)} \sum_{j_i} r_n(j_i) \leq (1+c_{n-1})
\alpha+c_{n-1}^{-1/2-\epsilon}+c_n^{1/(10 n)} \leq
(1+c_{n-1}^{1/2-4\epsilon}) \alpha_n
\ee

Moreover, if $\d \notin \hat D$, we either have the above estimate or
$l_n(\d) \leq c_n^{-2/n}$.  Thus we get
\be
\sum l_n(\d) C(\d,n) \leq (1+c_{n-1}^{1/3}) c_n^{-1} \alpha_n,
\ee
and, for every $j \in A_n$,
\be
\sum_{C^\d_n \subset I^j_n} l_n(\d) C(\d,n)
\leq I(j,n) (1+2 c_{n-1}^{1/3}) c_n^{-1} \alpha_n,
\ee
while, for $j \in B_n$ we have
\be
\sum_{C^\d_n \subset I^j_n} l_n(\d) C(\d,n)
\leq I(j,n) n c_n^{-1} \alpha_n.
\ee
}

\comm{
Let us consider $Z$ the set of all $j \in W$ with $r_n(j) \leq
c_{n-1}^{-14}$, and let $\hat Z$ be the set of $j$ such that $k(j) \in Z$. 
Notice that $I(\Z \setminus \hat Z,n+1) \leq c_{n-1}^{1/5}$.
}

Let $X_t$ be the set of $\d$ with $l_n(\d) \geq
t(1+c_{n-1}^{\frac {1} {2}-10 \delta^{(n)}_{14}})
c_n^{-1} \alpha_n$.  Notice that
\be
t \geq c_n^{1-2/n} \implies X_t \cap D=X_t \cap \hat D.
\ee
On the other hand, by (\ref {upper'}),
\be
\d \in X_t \setminus D \implies m(\d) \geq
t (1+c_{n-1}^{\frac {1} {2}-9 \delta^{(n)}_{14}}) c_n^{-1},
\ee
so, by (\ref {m(d)}), $C(X_t \setminus D,n) \leq
(1-c_{n-1}^{\frac {1} {2}-8 \delta^{(n)}_{14}})
e^{-t}$,
which gives by (\ref {hat D})
\be \label {563}
C(X_t,n) \leq (1-c_{n-1}^{\frac {1} {2}-8\delta^{(n)}_{14}})
(e^{-t}+e^{-c_{n-1}^{-25/3}}),
\quad t \geq c_n^{1-2/n}.
\ee
If $j \in Z$, we can estimate
\be \label {564}
C(X_t \cap \{\d \in \Omega,\, C^\d_n \subset I^j_n\},n)
\leq I(j,n) (1-c_{n-1}^{\frac {1} {2}-5\delta^{(n)}_{14}})
(e^{-t}+e^{-c_{n-1}^{-25/3}}), \quad t \geq c_n^{1-2/n}, j \in Z.
\ee

Let $Y_t$ be the set of $j \neq 0$ with
$R_n(I^j_{n+1})=C^\d_n$, $\d \in X_t$.
\comm{
Notice that
\be
\sum_{j \neq 0,
r_{n+1}(j)-v_n \geq (1+c_{n-2}) t \alpha_n c_n^{-1}} r_n(j)
I(j,n)=\int_t^\infty I(Y_t,n+1).
\ee
}
The following estimates are immediate from (\ref {563}), (\ref {564}):
\begin{align}
\label {569}
I(Y_t,n+1) &\leq 2^n (e^{-t}+e^{-c_{n-1}^{-8}})^{1/2}, &&t \geq
c_n^{1-2/n}\\
\label {570}
I(Y_t \cap \hat Z,n+1) &\leq e^{-t}+e^{-c_{n-1}^{-8}}, &&t \geq
c_n^{1-2/n}.
\end{align}
This last estimate implies in particular
\be \label {571}
I(Y_t, n+1) \leq e^{-t}+e^{-c_{n-1}^{-8}}+
c_{n-1}^{1/9} \leq (1+c_{n-1}^{1/20}) e^{-t},
\quad c_n^{1-2/n} \leq t \leq \ln c_{n-1}^{-1/20}.
\ee
Using additionally that by (\ref {r_n(j)}),
$I(Y_t,n+1) \leq e^{-t c_n^4}$ for $t \geq c_n^{-4}$,
and that obviously $I(Y_t,n+1) \leq 1$ for all $t$, we see that (\ref
{569}), (\ref {571}) imply
\be \label {bla3}
Y_t \leq s(t)=\left \{ \begin{array}{ll}
1 & \text {for } t <c_n^{1-2/n},
\\[5pt]
(1+c_{n-1}^{1/20}) e^{-t} & \text {for }
c_n^{1-2/n} \leq t<\ln c_{n-1}^{-1/20},
\\[5pt]
2^{n+1} e^{-t/2} & \text {for }
\ln c_{n-1}^{-1/20} \leq t<c_{n-1}^{-8},
\\[5pt] 2^{n+1} e^{-c_{n-1}^{-8}/2} & \text {for }
c_{n-1}^{-8} \leq t<c_n^{-5},
\\[5pt] e^{-t c_n^4} & \text {for } t \geq c_n^{-5}.
\end{array}
\right.
\ee
which gives
\be \label {bla1}
\int_0^\infty I(Y_t,n+1) dt \leq 1+c_{n-1}^{1/20}.
\ee
By definition of $X_t$ and $Y_t$, we have
\be \label {bla2}
0 \neq j \in Y_t \iff r_{n+1}(j) \geq v_n+
t (1+c_{n-1}^{\frac {1} {2}-10\delta^{(n)}_{14}}) c_n^{-1} \alpha_n,
\ee
so that (\ref {bla1}) implies
\be \label {above}
\alpha_{n+1} \leq v_n (1-c_{n+1})+(1+c_{n-1}^{\frac {1}
{2}-10\delta^{(n)}_{14}})
\alpha_n c_n^{-1}\int_0^\infty I(Y_t,n+1)
\leq (1+c_{n-1}^{1/30}) \alpha_n c_n^{-1}.
\ee

Estimates (\ref {below}) and (\ref {above}) imply (\ref {alpha_n}), shifting
$n$ to $n+1$.

Moreover, for any set $\Theta \subset \Z \setminus \{0\}$, with
$\theta=I(\Theta,n+1)$, (\ref {bla3}) implies
\be \label {Y_t}
\int_0^\infty I(Y_t \cap \Theta,n+1) dt \leq \int_0^\infty \min
\{\theta,s(t)\} dt \leq \frac {5} {2} \theta (1-\ln \theta)+\frac {c_n} {2},
\ee
which (together with (\ref {bla2}) and (\ref {below}))
implies (\ref {G}), shifting $n$ to $n+1$.
\end{pf}

\comm{

We have
\be
\int_0^{c_n^{1-2/n}} I(Y_t,n+1) dt \leq c_n^{1-2/n}
\ee
\be
\int_{c_n^{1-2/n}}^{c_{n-1}^{-1/20}} I(Y_t \cap Z,n+1) dt \leq
1+e^{-c_{n-1}^{-7}},
\ee
\be
\int_{c_n^{1-2/n}}^{c_{n-1}^{-1/20}} I(Y_t \setminus Z,n+1) dt \leq
c_{n-1}^{1/10},
\ee
\be
\int_{c_{n-1}^{-1/20}}^{c_{n-1}^{-8}} I(Y_t,n+1) dt \leq
e^{-c_{n-1}^{-1/50}},
\ee
\be
\int_{c_{n-1}^{-8}}^{e^{-c_{n-1}^{-7}}} I(Y_t,n+1) dt \leq
e^{-c_{n-1}^{-7}},
\ee
\be
\int_{e^{-c_{n-1}^{-7}}}^\infty I(Y_t,n+1) dt \leq
\int_{e^{-c_{n-1}^{-7}}}^\infty e^{-t c_n^4} dt \leq e^{-c_n^{-1}}.
\ee
It follows that
\be
\alpha_{n+1}=v_n (1-c_n)+(1+c_{n-2})
\alpha_n c_n^{-1}\int_0^\infty I(Y_t,n+1) \leq (1+2 c_{n-2}) \alpha_n
c_{n-1}^{-1}.
\ee

If $G \subset \Z \setminus \{0\}$ satisfies $I(G,n+1)=g$, we have
\be
\int_0^\infty I(Y_t \cap G,n+1) \leq g \ln g
\ee
}

We can now conclude:

\begin{lemma}[Large Deviation Estimate] \label {large}

Let $\Theta \subset \Z \setminus \{0\}$, $\theta=I(\Theta,n)$.
Let $F$ be the set of $\d$ such that
\be
\frac {1} {l_n(\d)} \sum_{j_i \in \Theta} r_n(j_i) \geq 4 (\theta (1-\ln
\theta)+c_{n-1}^{1/4}).
\ee
Then $C(F,n) \leq 2 c_n^{1-2/n}$.

\end{lemma}

\begin{pf}

By the previous considerations \S \ref {Dalpha}, except for $\d$ in
an exceptional set $D(\Theta)$
satisfying $C(D(\Theta),n) \leq c_n^{1-2/n}$, (\ref {upper}) holds, that is
\be \label {upper 1}
\frac {1} {m(\d)} \sum_{j_i \in \Theta} r_n(j_i) \leq
(1+\delta^{(n)}_5 c_{n-1}^{1/2}) \zeta^L(\Theta)+c_{n-1}^{-\frac {1}
{2}-\delta^{(n)}_5},
\ee
where, by Lemma \ref {alphanlemma},
\be
\zeta^L(\Theta) \leq (3 \theta (1-\ln \theta)+c_{n-1}) \alpha_n.
\ee
By the proof of Lemma \ref {alphanlemma}, except for $\d$ in
an exceptional set $D(\Z \setminus \{0\}) \subset \Omega$
satisfying $C(D(\Z \setminus \{0\}),n) \leq c_n^{1-2/n}$, (\ref
{upper'}) holds, that is
\be \label {upper' 1}
\frac {l_n(\d)} {m(\d)} \geq (1-c_{n-1}^{\frac {1} {2}-\delta^{(n)}_9})
\alpha_n.
\ee
Estimates (\ref {upper 1}) and (\ref {upper' 1}) imply that for
$\d \notin D(\Theta) \cup D(\Z \setminus \{0\})$,
\be
\frac {1} {l_n(\d)} \sum_{j_i \in \Theta} r_n(j_i) \leq \frac
{(1+\delta^{(n)}_5 c_{n-1}^{1/2})
\zeta^L(\Theta)+c_{n-1}^{-\frac {1} {2}-\delta^{(n)}_5}}
{(1-c_{n-1}^{\frac {1} {2}-\delta^{(n)}_9})\alpha_n}
\leq 4 \theta (1-\ln \theta)+c_{n-1}^{1/3}.
\ee
thus $F \subset D(\Theta) \cup D(\Z \setminus \{0\})$.  The result follows.
\end{pf}

\comm{

Let $\hat \alpha=\sum_{s \in L} s m_s$.  We notice that clearly
$\alpha_n \leq \hat \alpha$, but also $\alpha_n-\hat \alpha<
c_{n-1}^{-1/2-\epsilon}$
and so
\be
1 \leq \frac {\alpha_n} {\hat \alpha} \leq 1+c_{n-1}^{1/2-2\epsilon}.
\ee

Let $N(m)$ be the set of $\d$ with $m(\d)=m$ and, for every $s \in L$,
the number of entries in $M_s$, divided by $m$,
is outside $[m_s-c_n^{-1/4}, m_s+c_n^{1/4}]$.  Let $\hat N(m)$ be the set of
$\d$ with $m(\d)=m$ and, for every $s \in S$, the number of entries in
$M_s$, divided by $m$, is more than $2 c_n^{1/8}$.
If $m \geq c_n^{-1/2}$, we have
\be
I(N(m) \cup \hat N(m),n) \leq 2 e^{-c_n^{-1/8}}.
\ee

Let $H(m)$ be the set of $\d$ with $m(\d)=m$ and with no entries in $H$.

If $c_n^{-1/2} \leq m(\d) \leq c_n^{-2}$, $\d \notin N(m(\d)) \cup \hat
N(m(\d)) \cup K(m(\d) \cup H(m(\d)))$, then
\be
(1-c_n^{1/8}) \hat \alpha \leq
\sum_{s \in L} s (1-c_n^{1/8}) m_s \leq \frac {l_n(\d)} {m(\d)}
\ee
and
\be
\sum_{s \in L} s (1+c_n^{1/8}) m_s+\sum_{s \in S} 2 c_n^{1/8}
s+c_{n-1}^{-1/2+\epsilon} \geq \frac {l_n(\d)} {m(\d)},
\ee
which implies
\be
(1+c_{n-1}^{1/2-3\epsilon} \hat \alpha \geq \frac {l_n(\d)} {m(\d)}.
\ee

\begin{lemma}

We have
\begin{align}
C(\{|(l_n(\d)/m(\d))-\alpha_n|<c_{n-1}^{1/3}\},n)<c_n^{1/2}\\
C(\{|(l_n(\d)/m(\d))-\alpha_n|<c_{n-1}^{1/3},l_n(\d) \geq
c_n^{-2/3}\},n)<e^{-c_{n-1}^{-2}}.
\end{align}

\end{lemma}

\begin{pf}

Let $M_s=\{j \in A|r_n(j)=s\}$ and $m_s=I(M_s,n)$.
Let $L=\{s \leq c_{n-1}^{-14}|m_s \geq c_n^{1/8}\}$,
$S=\{s \leq c_{n-1}^{-14}|m_s<c_n^{1/8}\}$, $H=\{s>c_{n-1}^{-14}\}$.

Let $\hat \alpha=\sum_{s \in L} s m_s$.  We notice that clearly
$\alpha_n \leq \hat \alpha$, but also $\alpha_n-\hat \alpha<
c_{n-1}^{-1/2-\epsilon}$
and so
\be
1 \leq \frac {\alpha_n} {\hat \alpha} \leq 1+c_{n-1}^{1/2-2\epsilon}.
\ee

Let $N(m)$ be the set of $\d$ with $m(\d)=m$ and, for every $s \in L$,
the number of entries in $M_s$, divided by $m$,
is outside $[m_s-c_n^{-1/4}, m_s+c_n^{1/4}]$.  Let $\hat N(m)$ be the set of
$\d$ with $m(\d)=m$ and, for every $s \in S$, the number of entries in
$M_s$, divided by $m$, is more than $2 c_n^{1/8}$.
If $m \geq c_n^{-1/2}$, we have
\be
I(N(m) \cup \hat N(m),n) \leq 2 e^{-c_n^{-1/8}}.
\ee

Let $H(m)$ be the set of $\d$ with $m(\d)=m$ and with no entries in $H$.

If $c_n^{-1/2} \leq m(\d) \leq c_n^{-2}$, $\d \notin N(m(\d)) \cup \hat
N(m(\d)) \cup K(m(\d) \cup H(m(\d)))$, then
\be
(1-c_n^{1/8}) \hat \alpha \leq
\sum_{s \in L} s (1-c_n^{1/8}) m_s \leq \frac {l_n(\d)} {m(\d)}
\ee
and
\be
\sum_{s \in L} s (1+c_n^{1/8}) m_s+\sum_{s \in S} 2 c_n^{1/8}
s+c_{n-1}^{-1/2+\epsilon} \geq \frac {l_n(\d)} {m(\d)},
\ee
which implies
\be
(1+c_{n-1}^{1/2-3\epsilon} \hat \alpha \geq \frac {l_n(\d)} {m(\d)}.
\ee
\end{pf}

We now obtain easily:
\be
\left |\frac {\alpha_{n+1}} {\alpha_n c_n^{-1}}-1 \right | \leq
c_{n-1}^{1/4}.
\ee

\begin{lemma}

We have, for $t \geq 1$
\be
I(\{r_{n+1}(j) \geq 1+c_{n-1}^{1/4} t c_n^{-1}
\alpha_n,n+1\},n+1) \leq e^{-t}
\ee

\end{lemma}

\begin{pf}

Let $M_s=\{j \in A \cap P|r_n(j)=s\}$ and $m_s=I(M_s,n)$.
Let $L=\{s \leq c_{n-1}^{-14}|m_s \geq c_n^{1/8}\}$,
$S=\{s \leq c_{n-1}^{-14}|m_s<c_n^{1/8}\}$, $H=\{s>c_{n-1}^{-14}\}$.

Let $\hat \alpha=\sum_{s \in L} s m_s$.  We notice that clearly
$\alpha_n \leq \hat \alpha$, but also $\alpha_n-\hat \alpha<
c_{n-1}^{-1/2-\epsilon}$
and so
\be
1 \leq \frac {\alpha_n} {\hat \alpha} \leq 1+c_{n-1}^{1/2-2\epsilon}.
\ee

Let $N(m)$ be the set of $\d$ with $m(\d)=m$ and, for every $s \in L$,
the number of entries in $M_s$, divided by $m$,
is outside $[m_s-c_n^{-1/4}, m_s+c_n^{1/4}]$.  Let $\hat N(m)$ be the set of
$\d$ with $m(\d)=m$ and, for every $s \in S$, the number of entries in
$M_s$, divided by $m$, is more than $2 c_n^{1/8}$.
If $m \geq c_n^{-1/2}$, we have
\be
I(N(m) \cup \hat N(m),n) \leq 2 e^{-c_n^{-1/8}}.
\ee

Let $H(m)$ be the set of $\d$ with $m(\d)=m$ and with no entries in $H$.

If $c_n^{-1/2} \leq m(\d) \leq c_n^{-2}$, $\d \notin N(m(\d)) \cup \hat
N(m(\d)) \cup K(m(\d) \cup H(m(\d)))$, then
\be
(1-c_n^{1/8}) \hat \alpha \leq
\sum_{s \in L} s (1-c_n^{1/8}) m_s \leq \frac {l_n(\d)} {m(\d)}
\ee
and
\be
\sum_{s \in L} s (1+c_n^{1/8}) m_s+\sum_{s \in S} 2 c_n^{1/8}
s+c_{n-1}^{-1/2+\epsilon} \geq \frac {l_n(\d)} {m(\d)},
\ee
which implies
\be
(1+c_{n-1}^{1/2-3\epsilon} \hat \alpha \geq \frac {l_n(\d)} {m(\d)}.
\ee
\end{pf}

Let $P \subset \Z \setminus \{0\}$, and let $p=I(P,n)$.  Then
\be
C(\{\d=(j_i),\sum_{j_i \in P} r_n(j_i) \geq -p \ln p\},n) \leq
c_{n-1}^{1/10}+
\ee

We now conclude that for any $P \subset \Z \setminus \{0\}$,
\be
\sum_{j \in P} r_{n+1}(j) I(j,n+1) \leq -\alpha_{n+1} I(P,n+1) \ln I(P,n+1)
\ee
}

\subsection{Proof of Theorem~B}

\subsubsection{Series of reductions}

We will argue by contradiction.  If Theorem B is false,
there exists a positive measure set $\JJ_1$ of
non-regular parameters $\lambda$ such that the critical point is not
in the basin of the physical measure $\mu_{f_\lambda}$.
Since almost all parameters in $\JJ_1$ are
finitely renormalizable, there exists a subset $\JJ_2 \subset \JJ_1$
of positive measure of parameters
which are exactly $k$ times renormalizable, with some fixed $k$.

For each parameter $\lambda$ in $\JJ_2$, let us consider the sequence of
partitions $\Upsilon_n$ of the interval $I$ in connected components of
the domain
of the first landing map from $I$ to $I_n$.  Those partitions get more
refined as $n$ increases, the size of the largest component (of order at
most $c_{n-1}$) decreasing to $0$ with $n$.  Thus, there exists some $\n>0$
and a positive measure set of parameters $\JJ_3 \subset \JJ_2$,
such that for all parameters in $\JJ_3$ there exists at least
one component $\xi^\lambda \in \Upsilon_\n$ (that may be chosen to depend
measurably on $\lambda$) such that the asymptotic frequency of
the critical orbit in $\xi^\lambda$ either does not exists or is different
from $\mu_{f_\lambda}(\xi^\lambda)$.
Proceeding further, there exists $\epsilon>0$ and
a positive measure set $\JJ_4 \subset \JJ_3$
such that for all parameters in
$\JJ_4$.
\be \label {delta4}
\limsup \left | \frac {1} {k} \#\{i \leq k, f_\lambda^i(0) \in
\xi^\lambda\}-\mu_{f_\lambda}(\xi^\lambda) \right |>\epsilon.
\ee

The set $\JJ_4$ is contained in the union of parameter
intervals $J_\n$ ($\n$ fixed) associated to the principal nest (of $k$-th
renormalization).  It follows that at least one such interval $J_\n$
intersects $\JJ_4$ in a positive measure set $\JJ_5$.  For any
$\lambda_1,\lambda_2 \in J_\n$, there is a homeomorphism
$h[\lambda_1,\lambda_2]:I \to I$ such that $h[\lambda_1,\lambda_2]
\circ f_{\lambda_1}|(I \setminus I_{\n+1}[\lambda_1])=
f_{\lambda_2} \circ h[\lambda_1,\lambda_2]$.  Thus,
there exists a positive measure subset $\JJ \subset \JJ_5$ such that for
$\lambda_1,\lambda_2 \in \JJ$, $h[\lambda_1,\lambda_2]$ takes
$\xi^{\lambda_1}$ to $\xi^{\lambda_2}$.  In other words, the combinatorics of
$\xi^\lambda$ does not depend on $\lambda \in \JJ$.

In order to get a contradiction and
prove Theorem B, we will show that for almost every parameter in $\JJ$,
\be
\limsup \left | \frac {1} {k} \#\{i \leq k, f_\lambda^i(0) \in
\xi^\lambda\}-\mu_{f_\lambda}(\xi^\lambda) \right |<\epsilon.
\ee

To simplify the notation, we will write $\xi$ for $\xi^\lambda$.  We will
also write $\mu$ for $\mu_{f_\lambda}$.
For $x \in I$, and a measurable set $\Lambda
\subset I$, let
\be
\Psi(\Lambda,x,k)=\frac {1} {k} \#\{i \leq k,f^i(x) \in \Lambda\}.
\ee

Notice that if $\Lambda=\xi$, and $n>\n$ then
$x \mapsto \Psi(\Lambda,x,k)$ is
constant in each interval $I^j_n$ for $k \leq r_n(x)$,
while for $k \leq l_n(x)$, $x \mapsto \Psi(\Lambda,x,k)$ is constant
in each $C^\d_n$.  Those quantities
stay unchanged if we vary the parameter $\lambda$ inside some
$J_n$, if we keep the combinatorics constant, that is, if we choose
a varying point $x_\lambda$ inside $I^j_n[\lambda]$ or $C^\d_n[\lambda]$,
$j$ or $\d$ fixed.
\comm{
\subsubsection{Reduction to the study of fixed gaps}

Let us say that $P$ is an $n$-gap if it is a gap of the Cantor set of points
that never enter $I_n$.  The following is a convenient reduction of the
problem:

\begin{lemma}

If for all $n$-gap $P$,
$$
\lim_{k \to \infty} \frac {\#\{j < k|f^j(0) \in P\}} {k}
$$
exists and is equal to $\mu(P)$, then $0$ belongs to the basin of $\mu$.

\end{lemma}

So we just have to prove that for almost every non-regular parameter, for
all $n$-gaps, the distribution of the critical orbit is correct.

Since the set of $n$-gaps is countable, we just have to show that for any
fixed $n$-gap $P$, the set of non-regular
parameters with correct distribution of the critical orbit in $P$ has full
measure.  We must however define what is a fixed $n$-gap, since it has, of
course, to depend (as an interval) on the parameter.

Fix some level $\Delta_k$ of renormalization and some window $J_{n'}$
of the parapuzzle (the set of all such windows, for all $k$ and $n'$
is still countable).  Let us also fix $f \in J_{n'}$.
In $J_{n'}$, the Cantor set of points which do not enter (under iteration by
$f$) on $T^{(k)}_{n'}$ is persistent:
there exists a homeomorphism $h_g$ (with continuous
dependence on $g$) such that $h_g$ conjugates
$f|I \setminus T^{(k)}_{n'}$ to
$g|I \setminus T^{(k)}_{n'}[g]$.  Fixing a $n'$-gap $P_f$ for $f$,
we have defined a ``fixed gap along $J_{n'}$'' by taking $P_g=h_g(P_f)$.

Let $\JJ \subset J_{n'}$ be the set of non-regular
exactly $k$-times renormalizable maps.

We will now show that for almost every $g \in \JJ$,
the distribution of the critical orbit is correct in $P$, and thus conclude
the proof.
}

\subsubsection{Computing $\mu$ in the principal nest} \label {comp}

For $x \in I$, let $\vs_n(x)=\inf \{k,\, f^k(x)
\in I_n\}$, so that $f^{\vs_n(x)}(x)$ is the
first landing of $x$ in $I_n$.

For $x \in I_n$, $\vs_n(x)=0$, and in general we have
$\vs_{n+1}(x)-\vs_n(x)=l_n(f^{\vs_n}(x))$.  Notice that Lemma~\ref
{standard} implies that
\be
\frac {|\{x \in I_n,c_n^{-1/2}<l_n(x)<c_n^{-2}\}|} {|I_n|}<c_n^{1/3}.
\ee
Since each branch of the first landing map from $I$ to $I_n$ has
distortion bounded by $1+O(c_n)$ (see \cite {ALM}, Theorem 2.14),
we obtain the estimate
\be
\frac {|\{x \in I,c_n^{-1/2}<l_n(x)<c_n^{-2}\}|} {|I|}<2 c_n^{1/3}.
\ee
By Borel-Cantelli, for almost every $x$, for $n$ sufficiently big,
\be
c_n^{-1/2}<\vs_{n+1}(x)-\vs_n(x)<c_n^{-2}.
\ee
In particular
\be
\lim \frac {\vs_{n+1}(x)} {\vs_n(x)}=\infty, \quad \text {for almost
every } x \in I.
\ee

Thus, for all $\Lambda \subset I$ measurable, for almost every $x \in I$,
\be \label {landcomp}
\lim \Psi(\Lambda,f^{\vs_n(x)}(x),\vs_{n+1}(x)-\vs_n(x))=\mu(\Lambda).
\ee

Given a measurable subset $\Lambda \subset I$, we let $M_l(\Lambda,n,\delta)
\subset I_n$ be the set of all $x$ such that
\be
\left |\Psi(\Lambda,x,l_n(x))-\mu(\Lambda) \right | >\delta.
\ee
We let $M_r(\Lambda,n,\delta) \subset I_n$ be the set of all
$x$ such that
\be
\left | \Psi(\Lambda,x,r_n(x))-\mu(\Lambda) \right | >\delta.
\ee

\begin{lemma}

For any measurable set $\Lambda \subset I$, for any $\delta>0$,
\be \label {lcomp}
\lim \frac {|M_l(\Lambda,n,\delta)|} {|I_n|}=0,
\ee
\be \label {rcomp}
\lim \frac {|M_r(\Lambda,n,\delta)|} {|I_n|}=0.
\ee

\end{lemma}

\begin{pf}

Let $H_n$ be the set of $x \in I$, such that the first landing of $x$
on $I_n$ belongs to $M_l(\Lambda,n,\delta)$.  If (\ref {lcomp}) is not true,
using the small distortion of the first landing map we conclude that
$\limsup |H_n|>0$,
so there exists a positive measure set of $x$ which belong
to infinitely many $H_n$.  But this is incompatible with (\ref
{landcomp}).

Let $T_n \subset I_n$ be the union of $I^j_n$ with the following
properties:
\begin{align}
&I^j_n \subset R_n(I_{n+1}),\\
&\dist ((R_n|I_{n+1})^{-1}|I^j_n)<2,\\
&\dist (R_n|I^j_n)<2,\\
&r_n(j)<c_{n-1}^{-14}.
\end{align}

It follows that
\be
1-\frac {|(R_n|I_{n+1})^{-1}(T_n)|} {|I_{n+1}|}<c_{n-1}^{1/10}.
\ee
Let
\be
Y_{n+1}=\{x \in (R_n|I_{n+1})^{-1}(T_n), R_n^2(x) \in
M_l(\Lambda,n,\delta/2), l_n(R_n^2(x))>c_n^{-1/2}\}.
\ee
Then
\be
\frac {|Y_{n+1}|}
{|(R_n|I_{n+1})^{-1}(T_n)|} \leq 4 \left (\frac {|M_l(\Lambda,n,\delta/2)|}
{|I_n|}+c_n^{1/3} \right ),
\ee
thus
\be
1-\frac {|Y_{n+1}|} {|I_{n+1}|}
\leq 4 \left (\frac {|M_l(\Lambda,n,\delta/2)|}
{|I_n|}+c_n^{1/3}+c_{n-1}^{1/10} \right ),
\ee
so that, by (\ref {lcomp}),
$\lim \frac {|Y_{n+1}|} {|I_{n+1}|}=1$.
On the other hand, if $x \in Y_{n+1}$,
\begin{align}
|\Psi(\Lambda,x,r_{n+1}(x))-\mu(\Lambda)| &\leq
\frac {1} {r_{n+1}(x)} \left ( \frac {\delta} {2}
l_n(R_n^2(x))+r_{n+1}(x)-l_n(R_n^2(x)) \right )\\
\nonumber
&\leq \left ( \frac {\delta} {2} +
\frac {2 c_{n-1}^{-14}} {c_n^{-1/2}} \right ),
\end{align}
so that $Y_{n+1} \subset M_r(\Lambda,n,\delta)$ for $n$ sufficiently big.
\end{pf}

\subsubsection{Distribution of the critical orbit}

\begin{lemma} \label {gamma}

Let $S \subset I_n$ be a union of $I^j_n$.  For any $C \geq 1$:
\be
\lim_{\g \to 1} p_{\g,C}(S|I_n)=p_{1,C}(S|I_n).
\ee

\end{lemma}

\begin{pf}

If $X \subset I_n$ is any finite union of intervals, by compactness of
quasisymmetric maps we get
\be
\lim_{\g \to 1} p_{\g,C}(X|I_n)=p_{1,C}(X|I_n).
\ee

It is clear that for any $\g \geq 1$, 
$p_{\g,C}(S|I_n) \geq p_{1,C}(S|I_n)$.  On the other hand, since $I_n
\setminus \cup I^j_n$ is a regular Cantor set,
\be
\lim_{k \to \infty} p_{2,C}(x \in I^j_n, |j|>k|I_n)=0,
\ee
since the qs-capacity of gaps of generation $t$ decays exponentially
with $t$ (see Lemma~6.1 of \cite {AM1} for a related estimate).

Given $\delta>0$, we can fix a subset $S' \subset S$ which is a union of
finitely many $I^j_n$ such that
\be
p_{2,C}(S \setminus S'|I_n)<\delta.
\ee

Hence
\begin{align}
\limsup_{\g \to 1} p_{\g,C}(S|I_n) &\leq
\limsup_{\g \to 1} p_{\g,C}(S \setminus S'|I_n)
+\limsup_{\g \to 1} p_{\g,C}(S'|I_n) \leq \delta+p_{1,C}(S'|I_n)\\
\nonumber
& \leq \delta+p_{1,C}(S|I_n).
\end{align}
The result follows.
\end{pf}

We now specify this discussion to $\Lambda=\xi$ (the gap fixed at the
beginning).  We are now in situation to apply the Large Deviations
Estimate to obtain:

\begin{lemma}

For all $\delta>0$, for all $n_s$ sufficiently big, there exists $\g>1$
such that
\be
p_{\g,10}(M_r(\xi,n_s,\delta)|I_{n_s})<c_{n_s-1}^{1/20},
\ee

\end{lemma}

\begin{pf}

Let $n_0$ be very big and $\delta'>0$ be such that
\be
4 (\delta' (1-\ln \delta')+c_{n_0-1}^{1/4})<\frac {\delta} {2},
\ee
and $c_{n_0}^{1/3} \ll \delta$.

Let $n>n_0$ be such that
\be
\frac {|M_r(\xi,n,\delta/3)|} {|I_n|}<\delta'.
\ee
Notice that $M_r(\xi,n,\delta/3)=I^\Theta_n$ for some set $\Theta \subset \Z$,
and $I(\Theta,n)<\delta'$.
Let $F \subset \Omega$ be the set of $\d=(j_1,...,j_m)$ such that
\be
\frac {1} {l_n(\d)} \sum_{j_i \in \Theta} r_n(j_i)>\frac {\delta} {2}>
4 (\delta'(1-\ln \delta')+c_{n-1}^{1/4}).
\ee
Then, by the Large Deviation Estimate (Lemma \ref {large})
we get $C(F,n) \leq c_n^{1-3/n}$.
Let $F'=F \cup (\Omega \setminus LS(n))$.  It follows that
$C(F',n) \leq c_n^{2/7}$.
Let $E' \subset \Z$ be the set of $j$ such that $R_n(I^j_{n+1}) \subset
C^{F'}_n$.  Then $I(E',n+1) \leq c_n^{2/35}$.

Notice that if $x \in I^{\Z \setminus E'}_{n+1}$ and $\d=(j_1,...,j_m)$ is
such that $R_n(x) \in C^\d_n$, then
\be
\frac {v_n} {r_{n+1}(x)}<c_n^{1/3} \ll \delta.
\ee
(since $\d \in LS(n)$ and $r_{n+1}(x)>|\d|$) so we can conclude
\begin{align}
|\Psi(\xi,x,r_{n+1}(x))-\mu(\xi)|&\leq
\frac {1} {r_{n+1}(x)} \left ( v_n+\sum_{i \leq m}
r_n(j_i) |\Psi(\xi,R_n^i(x),r_n(j_i))-\mu(\xi)| \right )\\
\nonumber
&\leq \frac {1} {r_{n+1}(x)} \left ( v_n+
\sum_{\ntop {i \leq m,} {j_i \in \Theta}} r_n(j_i)+
\frac {\delta} {3}
\sum_{\ntop {i \leq m,} {j_i \notin \Theta}} r_n(j_i) \right )\\
\nonumber
&\leq \frac {1} {r_{n+1}(x)} \left (
v_n+\left ( \frac {\delta} {2}+\frac {\delta} {3} \right )
l_n(\d) \right ) <\delta.
\end{align}

So $I^{E'}_{n+1} \supset M_r(\xi,n+1,\delta)$.  But $I(E',n+1) \leq
c^{2/35}_n$ implies that
\be
p_{1,10}(I^{E'}_{n+1}|I_{n+1})<c_n^{2/39},
\ee
so the result now follows by Lemma \ref {gamma} with $n_s=n+1$.
\end{pf}

Let us select $\delta=\epsilon/3$, and using the previous lemma, we select
$n_s$ very large and such that $c_{n_s-1}^{1/400}<\delta$.
Let $T$ be such that $I^T_{n_s}=M_r(\xi,n_s,\delta)$.
Using Lemma \ref {ier 1} we get
\be
p_{\g,C_n} (I^{\Z \setminus VG(T,n_s,n)}_n|I_n) \leq c_{n-1}^{1/20}.
\ee
Using PhPa2 we get:

\begin{lemma}

For almost every parameter in $\JJ$,
for all $n$ sufficiently big, we have $\tau_n \in
VG(T,n_s,n)$.

\end{lemma}

Using Lemma \ref {ier 2}, we get, for $n$
sufficiently big,
\be
p_{\g,\tilde C_n} (C^{\Omega \setminus LC(T,n_s,n)}_n|I^{\tau_n}_n)<
c_{n-1}^{1/100}.
\ee

Using PhPa1 we get:

\begin{lemma}

For almost every parameter in $\JJ$,
for all $n$ sufficiently big, $R_n(0) \in C^\d_n$
with $\d \in LC(T,n_s,n)$.

\end{lemma}

Let us now consider a parameter which satisfies the conclusion of the two
previous lemmas.  Let us show that for $k$ big enough,
\be \label {typ}
\left | \Psi(\xi,0,k)-\mu(\xi) \right |<2 \delta<\epsilon.
\ee

Indeed, if $v_n+c_{n-1}^{-4/(n-1)}<k \leq v_{n+1}$, by Lemma \ref {ier 3}
\be \label {delta2}
|\Psi(\xi,f^{v_n}(0),k-v_n)-\mu(\xi)|<\delta+2 c_{n_s-1}^{1/200},
\ee
in particular, for $n$ big enough
\be \label {delta/2}
|\Psi(\xi,0,v_n)-\mu(\xi)|<3\delta/2.
\ee
Notice that (\ref {delta2}) and (\ref {delta/2}) imply (\ref {typ}) for
$n$ big enough and for
$v_n+c_{n-1}^{-4/(n-1)}<k \leq v_{n+1}$.  For $v_n \leq k \leq
v_n+c_{n-1}^{-4/(n-1)}$, (\ref {typ}) follows from (\ref {delta/2}) since
$v_n>c_{n-1}^{-1/2} \gg \delta^{-1} c_{n-1}^{-4/(n-1)}$ for $n$ big enough.

Thus, for almost every parameter in $\JJ$, (\ref {typ}) holds, which
contradicts (\ref {delta4}) and completes the proof of Theorem B.

\subsection{Proof of Corollary \ref {lyap}}

We want to show that
\be
\int \ln |Df| d\mu=\lim \frac {1} {k} \ln |Df^k(f(0))|=\lim \frac {1} {n}
\sum_{k=1}^n \ln |Df(f^k(0))|.
\ee
The fact that $0$ belongs to the basin of $\mu$ means that for all
continuous $\phi$,
\be
\frac {1} {n} \sum_{k=1}^n \phi(f^k(0))=\int \phi d\mu.
\ee
Since $\mu$ has no atoms, this formula still holds if $\phi$ is
a bounded function with at most finitely many discontinuities.
Unfortunately, $\ln |Df|$ is not bounded, so we only have, for every
$\delta>0$ small
\be
\int_{I \setminus (-\delta,\delta)} \ln |Df| d\mu=
\lim_{n \to \infty} \frac {1} {n}
\sum_{\ntop {1 \leq k \leq n} {f^k(0) \in I \setminus (-\delta,\delta)}}
\ln |Df(f^k(0))|.
\ee
Since
\be
\lim_{\delta \to 0} \limsup_{n \to \infty} \frac {1} {n}
\sum_{\ntop {1 \leq k \leq n} {f^k(0) \in (-\delta,\delta)}}
\ln |Df(f^k(0))| \leq 0,
\ee
we have to prove that for almost every non-regular parameter,
\begin{equation} \label {weakly regular}
\lim_{\delta \to 0} \liminf_{n \to \infty} \frac {1} {n}
\sum_{\ntop {1 \leq k \leq n} {f^k(0) \in (-\delta,\delta)}}
\ln |Df(f^k(0))|=0.
\end{equation}

Condition (\ref {weakly regular}) is called Weak Regularity by Tsujii.  In
\cite {Av}, Theorem 10.2 (see also Remark 10.3 in that paper),
it was shown that almost every non-regular
parameter in non-trivial analytic families of unimodal maps
satisfies (\ref {weakly
regular}), so, together with Theorem B, it implies Corollary \ref {lyap}.

\comm{
Condition (\ref {weakly regular}) is called by Tsujii Weak Regularity.
The proof that almost every non-regular parameter in a non-trivial family of
unimodal maps is given in Theorem~B of \cite {AM2}.
Let us just outline the main steps of the proof.

Tsujii's Theorem \cite {T2} implies that a Collet-Eckmann parameter
with subexponential recurrence of the critical point in an analytic family
(he actually needs much weaker requirements) satisfying a transversality
condition (for that specific parameter)
is a density point of Weak Regular parameters.

The transversality condition is shown in \cite {AM2} to be implied (for
Collet-Eckmann parameters) by transversality to the lamination by
topological conjugacy classes.  By \cite {ALM}, almost every non-regular
parameter is transversal to the lamination.  So almost every non-regular
parameter is a density point of Weak Regular parameters and
the result follows.
}

\section{Regularity of the physical measure} \label {hyperbolic}

\subsection{Outline}

Theorem C is a statement of regularity of $\mu$.  We can think of
$d\mu^K$ as a regularization of $\mu$, designed to allow an
understanding of the relation between ergodic and geometric properties of
hyperbolic Cantor sets.
Before tackling the problem of studying the regularization of
$\mu$, it is
important to understand the limitations on the regularity of $\mu$
and identify the source of the difficulties.

According to Theorem \ref {Below},
$d\mu$ is bounded from below on $A$ (by some constant
$C>0$).  As a consequence, if $T$ is an interval of radius $\epsilon$
centered on $0$ then
$\mu(f(T))=\mu(T)\geq 2C\epsilon\geq 2C\sqrt{|f(T)|}$.  This shows
that $d\mu$ has a ``pole'' at the critical value and, due to invariance of
$\mu$, there are also poles all along the orbit of the critical value.

In particular, for a general measurable set contained in $A$, the
measure-theoretical quantity
$\ln(\mu(A))$ only gives information about the geometric quantity
$\ln |A|$ up to a factor of $2$ (for $|A|$ small enough).  This is the
main reason why we have to introduce the regularization procedure.  We would
not be able to prove Theorem A just with general information on $\mu$.

According to \cite {MS}, this estimate on the non-regularity of $\mu$
is optimal: it implies that
$d\mu \notin L^2(I)$, but it is known that for maps satisfying the
Collet-Eckmann condition $d\mu \in L^p(I)$, $p<2$.  This is better
explained by Benedicks and Carleson \cite {BC}, who, for a smaller set
of parameters (contained in the set of
good Benedicks-Carleson parameters) described $d\mu$ as a sum of a
bounded distribution and infinitely many poles (called square-root
singularities by them) along the orbit of the critical value.  Although this
was not proved in general, this is the picture to imagine as a guide.

Since the critical orbit is distributed according to $\mu$, those poles are
everywhere (they are dense in the attractor).  However, not all is lost:

\begin{enumerate}

\item The strength of the poles decreases exponentially fast along the
critical orbit (because of the Collet-Eckmann condition),

\item The regularized
$d\mu^K$ averages over the gap, and dissipate the pole with a strength
proportional to the size of the gap.

\end{enumerate}

Thus, a naive argument to prove Theorem C would be to obtain, with total
probability, some ``quantitative transversality'' of the critical orbit
with respect to $K$ which would guarantee that strong poles are located in
big gaps.  For instance, we could expect that the time
of the first visit of the critical point to some gap of $K$ is inversely
proportional to the size of the gap.  Such a situation would imply that
strong poles fall in (very) big gaps and should help\footnote{One also
needs to guarantee that strong poles fall well inside a gap in order to
cantrol the effect on small nearby gaps.} us to conclude that
$d\mu^K \in L^p(I)$ for $1 \leq p<\infty$.

This would be much easier to deal with if the location of the successive
poles was independent and uniformly distributed with respect to Lebesgue
measure.  However, there
is quite a bit of interaction between different poles.  In particular, new
poles tend to show up more frequently near earlier poles than elsewhere
(since the critical orbits distributes according to $\mu$ which in turn
is more concentrated near the poles).

Our strategy will be to hierarchize the gaps according to the principal
nest.  To estimate the measure of a given gap, we will study their frequency
in return branches.  To estimate the possible increase in frequency
between levels (caused by the distortion originated on the poles),
we introduce a transversality condition (which we call
``Strong poles fall in big gaps''), which means that $R_n(0)$ falls
transversely enough with respect to the Cantor set of points that never land
on $I_n$ (the concept of transverse involves the hierarchy).
This analysis (which will be carried out in the next section) will
allow us to conclude the ``Main estimate'', which gives
bounds on the $\mu$-measure of gaps.

In this section we state the ``Strong poles fall in big gaps'' condition,
prove that it is a total probability one, and conclude Theorem C assuming
the Main estimate.

\subsection{The ``Strong poles fall in big gaps condition''}

We say that $f$ satisfies the ``Strong poles fall in big gaps'' condition if

\begin{description}

\item [SP1]
For all $\d \in \Omega$, $|\d| \geq 1$,
the distance between $R_n(0)$ and
$\partial I^\d_n$ is bounded by
\be
\frac {|I^\d_n|} {2^n |\d|^2},
\ee

\item [SP2]
$R_n(0) \in C^\d_n$, where $\d=(j_1,...,j_m)$ satisfies
\begin{align}
&r_n(j_i) \leq c_{n-1}^{-11}, \quad 1 \leq i \leq m,\\
&I(j_i,n) \geq e^{-c_{n-1}^{-12}}, \quad 1 \leq i \leq m,
\end{align}

\item [SP3]
$R_n(0) \in C^\d_n$, where $\d=(j_1,...,j_m)$, and for each $1 \leq i
\leq e^{c_{n-2}^{-4}}$ we have $R_{n-1}(I^{j_i}_n) \subset
C^{\d_i}_{n-1}$ where
$\d_i=(j^i_1,...,j^i_{s(i)})$ and
\begin{align}
&r_{n-1}(j^i_k) \leq c_{n-2}^{-11}, \quad 1 \leq k \leq s(i),\\
&I(j^i_k,n-1) \geq e^{-c_{n-2}^{-12}}, \quad 1 \leq k \leq s(i).
\end{align}

\end{description}

\begin{lemma} \label {SP}

Almost every non-regular parameter satisfies the ``Strong poles fall in big
gaps'' condition.

\end{lemma}

\begin{pf}

Let $\g$ be such that $\epsilon(\g)<\delta_0$, in the notation of \S \ref
{epsilon(gamma)}, with $\delta_0>0$ very small (say, $1/1000$).

Let us first deal with SP1.
We will consider two cases $|\d|=1$ and $|\d|>1$.  Let $d(X,Y)$
denote the distance between $X$ and $Y$.

In the first case, let $A_n$ be the set of $k$ such that
\be
d(I^k_n,\partial I^j_n) \leq \frac {|I^j_n|} {2^{3n/4}}, \quad \text {for
some } j \neq k.
\ee
Then
\be
p_\g(I^{A_n}_n|I_n)<2^{-n/2}.
\ee
Applying PhPa2, we see that with total probability, $\tau_n \notin A_n$ for
$n$ large enough.  This implies that,
with total probability, for $n$ big, if
\be
d(R_n(0),\partial I^j_n) \leq \frac {|I^j_n|} {2^n}
\ee
then $j=\tau_n$.  Let $B_n$ be the set of $\d$ such that $C^\d_n \subset
I^{\tau_n}_n$ and
\be
d(C^\d_n,\partial I^{\tau_n}_n) \leq \frac {|I^{\tau_n}_n|} {2^{3n/4}}.
\ee
Then
\be
p_\g(C^{B_n}_n|I^{\tau_n}_n)<2^{-n/2},
\ee
and by PhPa1 we see that $R_n(0) \notin C^{B_n}_n$ for $n$ large enough.
In particular, we conclude the result for $\d=1$.

In the second case, let $E(n)$ be the set of $\d$ such that there exists
some $\td$ with $|\td| \geq 2$ and
\be \label {bla}
d(C^\d_n, \partial I^{\td}_n) \leq \frac {|I^{\td}_n|}
{|\td|^{3/2} 2^{3n/4}}.
\ee
Let us show that
\be \label {vra}
p_\g(E(n)|I^{\tau_n}_n) \leq \frac {1} {2^{n/2}} \sum_{k \geq 2} k^{-4/3}.
\ee

Notice that if $I^{\td}_n \subset I^j_n$ with $j \neq \tau_n$, then no
$C^\d_n \subset I^{\tau_n}_n$ satisfies (\ref {bla}), since
\be
d(C^\d_n,\partial I^{\td}_n) \geq d(I^{\td}_n,\partial I^j_n) \gg
|I^{\td}_n|.
\ee
On the other hand, for each $I^{\td}_n \subset I^{\tau_n}_n$, the set
$E(\td)$ of $\d$ satisfying (\ref {bla})
has the property that, for any $h$ $\g$-qs,
\be
\frac {|h(C^{E(\td)}_n)|} {|h(I^{\td}_n)|}<\frac {1}
{2^{n/2} |\td|^{4/3}},
\ee
and since all $I^{\td}_n$ with $|\td|=k$ are disjoint, letting
$E(k,n)=\cup_{|\td|=k} E(\td)$, we get
\be
p_\g(C^{E(k,n)}_n|I^{\tau_n}_n) \leq \frac {1}
{2^{n/2} k^{4/3}}.
\ee
This implies (\ref {vra}).

Applying PhPa1, we see that with total probability, for $n$ big enough,
$R_n(0) \notin C^{E(n)}_n$, which gives SP1 for $|\d|>1$.

Let us consider SP2.  Notice that $r_n(j)<c_{n-1}^{-11}$ implies
$I(j,n)>e^{c_{n-1}^{-12}}$ for $n$ big,
since the derivative of $f$ is bounded.
Let $F(n)$ be the set of $j$ satisfying $r_n(j)>c_{n-1}^{-11}$, and $F'(n)$
be the set of $\d$ with at least one entry in $F(n)$.  We get
\be
p_\g(I^{F(n)}_n|I_n) \leq e^{-c_{n-1}^{-7}},
\ee
which implies using PhPa2 that $\tau_n \notin F(n)$ with total probability,
and thus
\be
p_\g(C^{F'(n)}_n|I_n) \leq e^{-c_{n-1}^{-7+\delta}},
\ee
\be
p_\g(C^{F'(n)}_n|I^{\tau_n}_n) \leq e^{-c_{n-1}^{-7+\delta}},
\ee
where $\delta$ goes to $0$ when $n$ grows.
Using PhPa1 we get $R_n(0) \notin C^{F'(n)}_n$ with total probability,
which implies SP2.

Let us consider SP3.  Keeping the notation of the previous discussion, let
$G(n+1)$ be the set of $j$ such that $R_n(I^j_{n+1}) \subset C^\d_n$ with
$\d \in F'(n)$.  Let $G'(n+1)$ be the set of $\d$ with at least one entry in
$G(n+1)$ among its first $e^{-c_{n-1}^{-4}}$ entries.  It follows that
\be
p_\g(I^{G(n+1)}_{n+1}|I_{n+1}) \leq e^{-c_{n-1}^{-7+\delta}}.
\ee
which by PhPa2 implies that $\tau_{n+1} \notin G(n+1)$ with total
probability and thus
\be
p_\g(C^{G'(n+1)}_{n+1}|I_{n+1}) \leq e^{-c_{n-1}^{-7+\delta}},
\ee
\be
p_\g(C^{G'(n+1)}_{n+1}|I^{\tau_{n+1}}_{n+1}) \leq e^{-c_{n-1}^{-7+\delta}}.
\ee
This implies, using PhPa1, that $R_{n+1}(0) \notin C^{G'(n+1)}_{n+1}$ with
total probability, which implies SP3.
\end{pf}

\comm{
\subsubsection{Some kinds of branches and landings}

Let us define the set of fine returns $F(n) \subset \Z$ as the set of all
$j$ such that $|I^j_{n}| \geq e^{-c_{n-1}^{-12}}$.

The set of legal landings $LL(n) \subset \Omega$ is the set of all
$\d=(j_1,...,j_m)$
satisfying $m<e^{c_n^{-3}}$ and $I^{j_i}_n$ is fine for all $i$.

A legal return $L(n) \subset \Z$ is the set of all $j$ such that
$R_{n-1}(I^j_n) \subset C^\d_n$ with $\d \in LL(n-1)$.

An acceptable landing $AL(n) \subset \Omega$ is the set of all $C^\d_n$ with
$\d=(j_1,...,j_r)$ and such that $j_i \in L(n)$, $1 \leq i \leq
e^{c_{n-2}^{-4}}$.

\subsubsection{The condition}

We say that $f$ satisfies the ``Strong poles fall in big gaps condition'' if
for all $n$ sufficiently big:

\begin{enumerate}

\item For all $\d \in \Omega$, the distance between $R_n(0)$ and
$\partial I^\d_n$ is bounded by
\be
\frac {|I^\d_n|} {2^n |\d|^2};
\ee

\item $R_n(0)$ belongs to a legal landing;

\item $R_n(0)$ belongs to an acceptable landing.

\end{enumerate}
}

\subsection{Main estimate}

\begin{thm}[Main estimate] \label {mest}

Let $f$ be a unimodal map with the following properties:

\begin{enumerate}

\item $f$ is Collet-Eckmann and has an absolutely continuous invariant
measure $\mu$;

\item The several asymptotic limits and
estimates described in \S \ref {values} hold;

\item $f$ satisfies the ``Strong poles fall in big gaps'' condition.

\end{enumerate}

Then, there exists $n_0>0$ such that for every
$\delta>0$, and all $n \geq n_0$,
there exists $C_n$ such that for any $I^j_n$,
$\mu(I^j_n)<C_n |I^j_n|^{1-\delta}$.

\end{thm}

It turns out that Theorem C implies that
we can take $n_0=1$ in the Main estimate.

\begin{rem}

We think that it is possible to refine the conditions of the Main estimate
(keeping total probability) in order to obtain better estimates for
$\mu(I^j_n)$ (of the type $-C|I^j_n|\ln |I^j_n|$ or even better).
It is an interesting problem whether
a bound of the type $C|I^j_n|$ is valid with total
probability.  Such a bound is equivalent to obtaining
$d\mu^K_f \in L^\infty$ in Theorem~C.

\end{rem}

\comm{

\begin{lemma}

With total probability, $f$ satisfies the hypothesis of the main estimate.

\end{lemma}

\begin{pf}

We let $\g=1+\epsilon$ ($\epsilon<1/1000$).

Strong away from the boundary.  There are two cases: $|\d|=1$ and $|\d|>1$.
In the first case, the set $A$ of $I^k_n \subset I_n$ which are at
distance at most $|I^j_n|/2^{3n/4}$ from some
$\partial I^j_n$ (with $j \neq k$) satisfy
$p_\g(A|I_n)<2^{n/2}$.  We can use PhPa2 to get the
parameter result.  This shows that $R_n(0)$ can only be $|I^j_n|/2^n$ close
to $I^j_n$ if $R_n(0) \in I^j_n$, that is $j=\tau_n$.  On the other hand,
the set $B$ of $C^\d_n \subset I^{\tau_n}_n$ which are
$|I^{\tau_n}_n|/2^{3n/4}$
close to $\partial I^{\tau_n}_n$ satisfy $p_\g(B|I^{\tau_n}_n)<2^{n/2}$, and
the parameter estimate follows by PhPa1.

In the second case, the set $C(k)$ of $C^\d_n \subset I^{\tau_n}_n$ with
$|\d|=k$
which are at distance at most
$|I^{\td}_n|/|\td|^{3/2} 2^{3n/4}$ from some $I^{\td}
\subset I_n$ satisfy $p_\g(C|I^{\tau_n}_n)<1/2^{n/2} |\d|^{4/3}$.

Indeed, if $I^{\td}_n \not \subset I^{\tau_n}_n$, then $C^\d_n$
is not close to $\partial I^{\td}_n$ by the previous case
(this distance is
bigger than the distance from $I^{\tau_n}_n$
to $I^j_n$ with $I^{\td}_n \subset I^j_n$,
which is bigger than $|I^j_n|/2^n>|I^{\td}_n|/2^{3n/2} k^{4/3}$).

On the other hand, for each $I^{\td}_n \subset I^\tau_n$, the set
$C(\td)$ of
$C^\d_n$ which intersect a $|I^{\td}_n|/2^{3n/2} |\d|^{3/2}$
neighborhood of $I^{\td}_n$ satisfy $|h(C(\td))|/|I^{\td}_n|<
2^{n/2} k^{4/3}$, 
and the result follows since $I^\d_n \subset I^{\tau_n}_n$ with
$|\d|=k$ are disjoint.

The result follows using PhPa2
since $\sum_n 2^{n/2} \sum_k 1/k^{4/3}<\infty$.

Legal.  Notice that a return is fine provided $r_n(j)<e^{c_{n-1}^{-11}}$.
Using the estimate on
$v_n$, $\tau_n \in F(n)$.  So $p_\g(x \in C^\d_n, \d \notin LL(n)|x \in
I^{\tau_n}_n)<e^{-c_{n-1}^{-5}}$, and the result follows from PhPa1.

Acceptable.  As before, we may estimate $p_\g(x \in I^j_n|\d \notin L(n)|x
\in I_n)<e^{-c_{n-2}^{-5}}$.  By PhPa2, $\tau_n \in L(n)$.  It follows that
$p_\g(x \in C^\d_n|x\in I^{\tau_n}_n)<e^{-c_{n-2}^{-4}}$.  The result
follows by PhPa1.
\end{pf}
}

\subsection{Proof of Theorem~C assuming the Main estimate}

By Lemma~\ref {SP} and the results of \cite {AM1},
we get that, with total probability, $f$ satisfies the
hypothesis of the Main estimate.  Let us now fix such an $f$.

If $K$ is a hyperbolic set for $f$, then it avoids a neighborhood of
the critical point.  On the other hand, if $K \subset K'$ and
$d\mu^{K'} \in L^p$, then necessarily
$d\mu^K \in L^p$.
So we just have to
consider the case of $K_n$, the maximal invariant of $I \setminus I_n$
for $n$ big.  The gaps of $K_n$ are connected components of the domain of
the first landing map from $I$ to $I_n$.

We will use the following:

\begin{lemma} \label {M^k}

For all $n>0$, there exists a finite partition of
$I \setminus I_n$ on intervals $M_i$, such that for each $M_i$, $f|M_i$ is a
diffeomorphism onto the union of some $M_i$ and, possibly, $I_n$.
Moreover, there exists constants $C>0$ and $t<1$
such that for any $x$ such that $x,...,f^k(x) \in I \setminus I_n$,
we can associate an interval $M^k(x)$ such that

\begin{enumerate}

\item $f^k:M^k(x) \to I$ is a diffeomorphism over some $M_i$;

\item Two intervals $M^k(x)$ and $M^k(y)$ are
either disjoint or coincide.

\item $|M^k(x)|<C t^k$;

\item $\sum_{k \geq 1} |M^k(x)|<C$;

\item The distortion of $f^k|M^k(x)$ is bounded by $C$;

\item For each $k$, $|\cup M^k(x)|<C t^k$.

\end{enumerate}

\end{lemma}

\begin{pf}

Let $Q$ be the finite set consisting of all points in the forward orbit of
$\partial I_n$.  Let $M_i$ be the connected components of $I \setminus (Q
\cup I_n)$.  It is clear that the image of $M_i$ consists of a union of
$M_j$, possibly together with $I_n$.  The $M_i$ form a Markov partition of
$I \setminus I_n$, and so the first and second item follow.
The third item follows from hyperbolicity of $f|I
\setminus I_n$ (see Lemma~\ref {hyperbolici}),
and the fourth follows from the third.  The fifth follows from
the fourth by a classical argument (it is enough to use that
$\ln |Df|$ is H\"older in $I \setminus I_n$).  Notice that for each $i$,
there exists $j_i>0$ such that $f^{j_i}(M_i)$ contains $I_n$.  This and the
fifth item show that $|\cup M^{k+j}(x)| \leq t |\cup M^k(x)|$ for some
$t<1$ and for $j=\max j_i$, and this gives the sixth item.
\end{pf}

\begin{cor} \label {gap of K_n}

For all $n>0$, there exists $\theta_n>0$ such that
\be
\sum_{\Lambda \mathrm {\,\, gap\,\, of\,\,} K_n} |\Lambda|^{1-\theta_n}<\infty.
\ee

\end{cor}

\begin{pf}

Let us say that a gap $\Lambda$ of $K_n$ is of generation $k$ if
$f^k(\Lambda)=I_n$.  Let $k(\Lambda)$ be the generation of $\Lambda$.
Notice that each $M^k(x)$ contains at most
one gap of generation $k+1$ (and no gaps of generation $\leq k$).  On the
other hand, each gap $\Lambda$
of generation $k+1$ is contained on some $M^k(x)$, which we denote
$M(\Lambda)$.  Notice that since the derivative of $f$ is bounded by some
constant $\kappa$,
\be
|\Lambda| \geq |I_n| \kappa^{-k(\Lambda)}.
\ee
We can estimate
\begin{align}
\sum_{\Lambda \text { gap of } K_n} |\Lambda|^{1-\theta_n} &\leq
\sum_{k \geq 0}
\sum_{k(\Lambda)=k} |M(\Lambda)|^{1-\theta_n}
\leq \sum_{k \geq 0} \sum_{k(\Lambda)=k} |\Lambda|^{-\theta_n}
|M(\Lambda)|\\
\nonumber
&\leq |I_n|^{-\theta_n} \left (\sum_{k \geq 0} \kappa^{\theta_n k}
\sum_{k(\Lambda)=k} |M(\Lambda)| \right )
\leq C |I_n|^{-\theta_n} \sum_{k \geq 0} \kappa^{\theta_n k} t^k,
\end{align}
where $C>0$, $t<1$ comes from item (5) of Lemma \ref {M^k}.  The result
follows with $\theta_n>0$ such that $\kappa^{\theta_n} t<1$.
\end{pf}

\subsubsection{}
Let $K^r_n$ be the Cantor set $I_n \setminus \cup I^j_n$, and let
$d\mu^{K^r_n}$ be the function which takes, in each gap of $K^r_n$
the average value of $d\mu$ in that gap, and let
$d\mu^{K^r_n}=0$ outside $I_n$.

\comm{
Since $f|I \setminus I_{n+1}$ is hyperbolic, there exists $A>0$,
$t<1$ such that
\be
I(\{r_n(j)=k\},n)<A t^k.
\ee
On the other hand, since the derivative of $f$ is bounded by $4$,
\be
\min_{r_n(j)=k} I(j,n)>4^{-k}.
\ee
This implies that, if $\theta_n$ satisfies $t 4^{\theta_n}<1$,
there exists $C>0$ such that
}
Notice that by Corollary \ref {gap of K_n},
\be
\sum_j |I^j_n|^{1-\theta_{n+1}} \leq \sum_\d |C^\d_n|^{1-\theta_{n+1}} \leq
C<\infty.
\ee

Fix $1 \leq p<\infty$.  Using the Main Estimate, let $C'$ be such that
\be
\mu(I^j_n)<C' |I^j_n|^{1-\frac {\theta_{n+1}} {p}}.
\ee
We estimate
\be
\int_{I_n} (d\mu^{K^r_n})^p=
\sum |I^j_n|\left (\frac {\mu(I^j_n)} {|I^j_n|} \right )^p \leq
C'^p \sum |I^j_n|^{1-\theta_{n+1}}<C C'^p.
\ee
In particular, $d\mu^{K^r_n} \in L^p$.

\comm{
\begin{lemma}

For every $n$ there exists a $\theta_n>0$ and $C>0$ such that
\be
\sum_j |I^j_n|^{1-\theta_n}<C.
\ee
In particular, for all $p<\infty$ we have
$d\mu^{K^r_n} \in L^p$.

\end{lemma}

\begin{pf}

Since $f|I \setminus I_{n+1}$ is hyperbolic,
$|\cup_{r_n(j)=k} I^j_n|<A t^k$, for some $t<1$.  On the other hand,
$\min_{r_n(j)=k} |I^j_n|>B \tilde t^k$ for some $\tilde>0$.  This gives the
first result for $\theta_n<\tilde t/t$.

For the second, using the Main Estimate, let $C'$ be such that
\be
\mu(I^j_n)<C' |I^j_n|^{1-\theta_n/p}.
\ee
We estimate
\be
\int E(d\mu|\BB^r_n)^p=
\sum |I^j_n|\left (\frac {\mu(I^j_n)} {|I^j_n|} \right )^p \leq
C'^p \sum |I^j_n|^{1-\theta_n}<C C'^p.
\ee
\end{pf}
}

\subsubsection{}

Given $\Lambda \subset I \setminus I_n$ measurable, let
$\hat \Lambda \subset I$ be the set of $x$ such that $\min\{k
\geq 1, f^k(x) \in \Lambda\} \leq \min\{k \geq 1,f^k(x) \in I_n\}$ (that is,
the orbit of $f(x)$ intersects $\Lambda$ before intersecting $I_n$).  Let
$\Lambda_l=\hat \Lambda \cap (I \setminus I_n)$ and $\Lambda_r=\hat \Lambda
\cap I_n$.

Notice that $f^{-1}(\Lambda_l \cup \Lambda)=\Lambda_r \cup \Lambda_l=\hat
\Lambda$, thus
\be \label {Lambda_l}
\mu(\Lambda)=\mu(\Lambda_r) \quad \text {provided } \Lambda_l \cap
\Lambda=\emptyset.
\ee

\comm{
\begin{lemma}

If $\Lambda_l \cap \Lambda=\emptyset$ then $\mu(\Lambda_r)=\mu(\Lambda)$.

\end{lemma}
}

\comm{
\begin{lemma}

There exists a finite partition of $I \setminus I_n$ on intervals $M_i$
and $C,C',C''>0$ and $t<1$
such that for any $x$ such that $x,...,f^k(x) \in I \setminus I_n$,
we can associate an interval $M^k(x)$ such that

\begin{enumerate}

\item $f^k:M^k(x) \to I$ is a diffeomorphism over some $M_i$;

\item Two intervals $M^k(x)$ and $M^k(y)$ are
either disjoint or coincide.

\item $|M^k(x)|<C t^k$;

\item $\sum_{k \geq 1} \sum |M^k(x)|<C$;

\item The distortion of $f^k|M^k(x)$ is bounded by $C$;

\item For each $k$, $|\cup M^k(x)|<C t^k$.

\end{enumerate}

\end{lemma}

\begin{pf}

Let $Q$ be the finite set consisting of all points in the forward orbit of
$\partial I_n$.  Let $M_i$ be the connected components of $I \setminus I_n
\setminus Q$.  It is clear that the image of $M_i$ consists of a union of
$M_j$, possibly together with $I_n$.  The $M_i$ form a Markov partition of
$I \setminus I_n$, and so the first and second item follow.
The third item follows from hyperbolicity of $f|I
\setminus I_n$, and the fourth follows from the third.
The fifth follows from
the fourth classically (it is enough to use that
$\ln |Df|$ is H\"older in $I \setminus I_n$).  Notice that for each $i$,
there exists $j_i>0$ such that $f^{j_i}(M_i)$ contains $I_n$.  This and the
fifth item show that $|\cup M^{k+j}(x)| \leq t |\cup M^k(x)|$ for some
$t'<1$ and for $j=\max j_i$, and this gives the sixth item.
\end{pf}
}

Let $\Lambda(k) \subset \Lambda_l$ be the set of points $x$ with
$k=\min\{i>0,f^i(x) \in \Lambda\}$.  Then, by Lemma \ref {M^k},
$\Lambda(k)$ is covered by disjoint intervals $M^k(y_i)$.  By items $5$ and
$6$ of Lemma \ref {M^k}, there exists $C>0$, $t<1$ such that
\be
|\Lambda(k)|<C t^k \max_i \frac {|\Lambda(k) \cap M^k(y_i)|} {|M^k(y_i)|}<
C' t^k|\Lambda|,
\ee
since the density of $\Lambda(k)$ inside some $M^k(y_i)$ is comparable
with the density of $\Lambda$ in $f^k(M^k(y_i))$ (by bounded distortion)
which is at most $\max_j |\Lambda| |M_j|^{-1}$.
Thus, there exists a constant $C>0$ such that
\be \label {C>0}
|\Lambda_l|<C |\Lambda|
\ee
independently of $\Lambda$.

\comm{
\begin{lemma}

If $\Lambda \subset I \setminus I_n$ is measurable then
\be
|\Lambda_l|<C |\Lambda|.
\ee

\end{lemma}

\begin{pf}

Let $\Lambda(k) \subset \Lambda_l$ be the set of points $x$ such that
$f^k(x) \in \Lambda$ but for $0 < i < k$, $f^i(x) \notin \Lambda$.  Then, by
the previous Lemma, $\Lambda_k \subset \cup M^k$.  By the fifth and sixth
item of the previous Lemma, there exists $C>0$, $t<1$ such that
$$
|\Lambda(k)|<C t^k \max_i \frac {|\Lambda \cap M_i|} {|M_i|}<C't^k|\Lambda|.
$$
This gives the result.
\end{pf}
}

\subsubsection{}

Let now $\Lambda$ be a gap of $K_n$.  Assume first that
$\Lambda$ is a gap of $K_n$ which does not intersect $\{f^i(0),0
\leq i<v_n\}$.  In particular, $\Lambda \neq I_n$ and
$\Lambda_r$ does not contain the critical point.
Since $\Lambda$ is a connected component of the domain of the first
landing map from $I$ to $I_n$, we have that
$\Lambda_l \cap \Lambda=\emptyset$, and thus
(\ref {Lambda_l}) holds.
By (\ref {C>0}),
\be
|\Lambda_r| \leq 2 |\Lambda_l| \sup_{x \in I_n \setminus I^0_n} |Df(x)|^{-1}
\leq C |\Lambda|
\ee
for some constant $C>0$ independent of $\Lambda$.
Using the H\"older inequality we get
\be \label {lambdar}
\mu(\Lambda)=\mu(\Lambda_r) \leq
\left ( \int_{\Lambda_r} (d\mu^{K^r_n})^p \right)^{\frac {1} {p}}
|\Lambda_r|^{1-\frac {1} {p}} <
C_p |\Lambda_r|^{1-\frac {1} {p}} \leq C'_p |\Lambda|^{1-\frac {1} {p}},
\ee
where $C'_p$ depends on $p$ but not on $\Lambda$.

The set of gaps $\Lambda$ of $K_n$ which intersect
$\{f^i(0),0 \leq i<v_n\}$ are in finite number, so there exists $C>0$ such
that for any such $\Lambda$,
\be \label {lambdar1}
\mu(\Lambda) \leq C|\Lambda|.
\ee

Putting together (\ref {lambdar}) and (\ref {lambdar1}), and varying $p$, we
see that for any $\delta>0$ there exists a constant $C(\delta)$ such that
for any $\Lambda$ gap of $K_n$ we have
\be
\mu(\Lambda) \leq C(\delta) |\Lambda|^{1-\delta}.
\ee

By Corollary \ref {gap of K_n}, for $\delta<\theta_n/p$ we have
\be
\int (d\mu^{K_n})^p=
\sum_{\Lambda \text { gap of } K_n}
|\Lambda|\left (\frac {\mu(\Lambda)} {|\Lambda|} \right )^p \leq
C(\delta)^p \sum |\Lambda|^{1-\theta_n}<C C(\delta)^p.
\ee

\comm{
For the second we argue as in Lemma,
$$
\sum_{\Lambda \text{ gap of } K_n} |\Lambda|^{1-\theta_n}<\infty
$$

\begin{lemma}

For all $n$ sufficiently big, for all $\delta>0$ there exists $C>0$ such
that if $\Lambda$ is a gap of $K_n$ then
$$
\mu(\Lambda)<C |\Lambda|^{1-\delta}.
$$
In particular, for all $p<\infty$ we have $d\mu^{K_n} \in L^p$.

\end{lemma}

\begin{pf}

We may assume that $\Lambda_r$ does not contain the critical point,
since there exists only finitely many gaps $\Lambda$ with this property. 
Let $U$ and $V$ be the two components of $I_n \setminus I^0_n$.  Then $f|U$
and $f|V$ have bounded distortion.  This shows that there exists a constant
$C'$ such that $|\Lambda_r|<C' |\Lambda_l|$,  Using the H\"older inequality
we get
$$
\mu(\Lambda)=\mu(\Lambda_r) \leq
\left ( \int_{\Lambda_r} E(d\mu|\BB^r_n)^p \right)^{1/p}
|\Lambda_r|^{1-1/p} <
C_p |\Lambda_r|^{1-1/p}.
$$
Since $p<\infty$ is arbitrary, we get the first result using the previous
Lemma.  For the second we argue as in Lemma,
$$
\sum_{\Lambda \text{ gap of } K_n} |\Lambda|^{1-\theta_n}<\infty
$$
for some $\theta_n$, and together with the first result this gives the
estimate.
\end{pf}
}

\section{Proof of the Main Estimate} \label {est}

\subsection{Outline}

Our problem is to analyze the asymptotics of the physical
measure of $I^j_n$ as the Lebesgue measure of $I^j_n$ decreases, $n$ fixed.
Let us fix some advanced level $\n$.
Fix a small interval $I^\j_\n$.
To the end of this section, $a=|I^\j_\n|$.  The
critical step $l$ is defined as the unique number with
\be
c_{l+1} \leq a < c_l.
\ee
Since our estimate is only relevant if $I^\j_\n$ is small,
we can assume that $l$ is very big.

The proof will be based on the analysis, for each level $n \geq \n$, of the
frequency of visits to $I^\j_\n$ before a return to $I_n$.  Those estimates
can be passed from level to level if one can control the distortion
introduced by the critical orbit.  The argument will take distinct steps.

In the early stages (before $l-1$), very few branches (measure of order $a$)
pass at all in $I^\j_\n$ before returning.
The critical orbit falls in big holes away from the hierarchic
structure of those branches and does not distort much the measure.

In the later stages (after $l+2$) most branches have total time much bigger
then $a^{-1}$, and they spend a proportion of time of order $a$ in
$I^\j_\n$.  The exceptional branches have measure much smaller than $a$, and
we use the inductive estimate of \S \ref {induc}
to show that they do not contribute much for the next levels.

In the intermediate stage, there is a delicate
transition between those two situations.
To complicate further, at this moment
the position of the critical point could introduce distortion of strength
comparable with $a$. 
We will need to use the hierarchical structure of the set of branches
passing through $I^\j_\n$
combined with our conditions on the critical orbit to control the
distortion of pullbacks.

{\it In what follows, $\delta$ will denote several constants which go to $0$
uniformly as the critical step $l$ goes to infinity}.

\subsection{Preliminaries}

Let us define
\be
X_n(j)=\#\{k<r_n(j), f^k(I^j_n) \subset I^\j_\n\},
\quad n \geq \n,
\ee
\be
\X_n(\d)=\#\{k<l_n(\d), f^k(C^\d_n) \subset I^\j_\n\},
\quad n \geq \n,
\ee
so that
\be
\X_n(\d)=\sum_{i=1}^m X_n(j_i), \quad \d=(j_1,...,j_m).
\ee
\be
X_{n+1}(j)=X_n(0)+\X_n(\d), \quad R_n(I^j_{n+1}) \subset C^\d_n.
\ee

Define
\be
x_n(r)=I(\{X_n(j) \geq r\},n), \quad x_n=x_n(1),
\ee
\be
\x_n(r)=C(\{\X_n(\d) \geq r\},n), \quad \x_n=\x_n(1).
\ee
One immediately gets (see \S \ref {ld1} for a derivation)
\be
\x_n \leq n c_n^{-1} x_n.
\ee

\comm{
\begin{align}
\x_n &\leq \sum_{k \geq 1} I(\{\d=(j_1,...,j_k),\, X_n(j_k) \geq 1\},n)\\
\nonumber
&\leq 2^{n/2} \sum_{k \geq 1} x_n I(\{|\d|=k-1\},n)\\
\nonumber
&\leq 2^n \sum_{k \geq 1}
x_n c_n^{-1} C(\{|\d|=k-1\},n)\\
\nonumber
&=2^n c_n^{-1} x_n.
\end{align}
}

Before getting into the more complicate intermediate steps, let us deal with
the initial steps and discuss our strategy for the later steps.

\subsubsection{Initial steps}

Notice that if $x_n<c_n$ then $X_n(0)=0$, and if additionally
$\x_n<c_{n+1}^5$ then we conclude that $x_{n+1}<c_{n+1}$ and
$X_{n+1}(0)=0$ as well.  In this case we can estimate
\be
x_{n+1} \leq \dist (R_n|(I_{n+1} \setminus I_{n+2})) 2^n \x_n \leq
2^{2n} c_{n+1}^{-1} x_n.
\ee
Since $X_\n(0)=0$ and $x_\n(0)=a$, we conclude by induction
that for $\n \leq n \leq l-1$ we have
\be \label {X_n}
X_n(0)=0, \quad x_n \leq a c_n^{-1} c_{n-1}^{-5/2} \ll a^{1-\delta}<c_n,
\quad \x_n \leq a c_n^{-9/4}
\ee
(using that $a c_n^{-9/4}<c_{n+1}^5$ for $n \leq l-2$).
Thus we have, just before the critical time:
\be \label {precritical}
x_{l-1}, \x_{l-1} \leq a c_{l-1}^{-5/2} \ll a^{1-\delta}, \quad
X_{l-1}(0)=0.
\ee

\subsubsection{Later steps} \label {l+4}

Let us fix, for the end of this section, some very small $\epsilon$.
Our aim is to estimate
\be
\mu(I^\j_\n) \leq a^{1-10 \epsilon^{1/2}}.
\ee
To attack this problem,
we will need to compute $\mu(I^\j_\n)$ somehow.  Using the idea of \S \ref
{comp}, one sees that we only have to show that
\be \label {VG12}
\lim_{n \to \infty} I(\{j,X_n(j) \leq a^{1-10 \epsilon^{1/2}} r_n(j)\},n)=1.
\ee

This will be done in the following way.  We will show that there exists
$n_0$ such that, defining
\be
T=\{j,X_{n_0}(j)>a^{1-8\epsilon^{1/2}} r_{n_0}(j)\},
\ee
we have
\be
I(T,n_0)<c_{n_0-1}^{1/20}.
\ee

Our actual choice of $n_0$ will be $n_0=l+3$ if $a<c_{l+1}^{\epsilon^{1/2}}$
and $n_0=l+2$ otherwise.  Notice that in both cases
$c_{n_0-1}^{1/200}<a$.

By the inductive estimate, see \S \ref {esttime} (we only need the case
$\g=1$ corresponding to Lebesgue measure),
we get, for $n \geq n_0$,
\be \label {VG1}
I(\Z \setminus VG(T,n_0,n),n) \ll c_{n-1}^{1/20},
\ee
and for $j \in VG(T,n_0,n)$
\be \label {VG2}
X_n(j) \leq (a^{1-8 \epsilon^{1/2}}+
c_{n_0-1}^{1/200})r_n(j)<a^{1-10 \epsilon^{1/2}} r_n(j).
\ee
Estimates (\ref {VG1}) and (\ref {VG2}) imply (\ref {VG12}).

\subsection{Transition from $l-1$ to $l$} \label {l-1}

The analysis of the transition from the $l-1$ level to the $l$ level is more
complicated.  We will need to consider a special sequence $E_t$ of nested
intervals in level $I_{l-1}$
around the critical value $R_{l-1}(0)$,
where we can analyze (using SP2)
the density of the set of points visiting $I^\j_\n$.
We then pullback this information to a sequence $F_t$ of nested intervals
in level $l$ around the
critical point, with a control of the distortion by SP1.

\subsubsection{}
Let $\d=(j_1,...,j_s)$ be such that $R_{l-1}(0) \in C^\d_{l-1}$ and for each
$0 \leq t \leq s$, let $\d(t)=(j_1,...,j_t)$.  By condition SP2
we can estimate
\be \label {aa1}
I(j_i,l-1)>e^{-c_{l-2}^{-12}} \gg a^{\delta}, \quad 1 \leq i \leq s.
\ee

Define
\be
E_t=I^{\d(t)}_{l-1}, \quad 0 \leq i \leq s.
\ee
From (\ref {aa1}) we get
\be \label {aa3}
\frac {|E_{t+1}|} {|E_t|}>2^{-l} e^{-c_{l-2}^{-12}} \gg a^{\delta}.
\ee
Denote
\be
F_t=(R_{l-1}|I_l)^{-1}(E_t).
\ee

Notice that (\ref {aa1}), (\ref {precritical}) imply that $X_{l-1}(j_i)=0$.
In particular,
\be
C^{\{\X_{l-1}>0\}}_{l-1} \cap E_t=(R^t_{l-1}|E_t)^{-1}
C^{\{\X_{l-1}>0\}}_{l-1}.
\ee
This forces the density estimate
\be \label {est 1}
\frac {|C^{\{\X_{l-1}>0\}}_{l-1} \cap E_t|} {|E_t|}<a c_{l-1}^{-3} \leq
a^{1-\delta}.
\ee

Moreover, since $\X_{l-1}(\d)=\sum_{i=1}^s X_{l-1}(j_i)=0$ and
$X_{l-1}(0)=0$ (by \ref {precritical}), we have that
\be \label {X_l}
X_l(0)=0.
\ee


\comm{
Using the SP1, we see that
\be
|E_t| \geq |R_{l-1}(I_l) \cap E_{t+1}| \geq \frac {|E_{t+1}|} {2^l (t+1)^2}
\ee
so that
\begin{equation} \label {F}
\frac {|F_t|} {|F_{t+1}|} \leq 2^l (t+1)^2 \frac {|E_t|} {|E_{t+1}|}.
\end{equation}
}

\subsubsection{}
Let us define $a_t$ and $b_t$ by
$E_t=(R_{l-1}(0)-a_t,R_{l-1}(0)+b_t)$.
Let ${\bf m}_t=\min\{a_t,b_t\}$ and ${\bf M}_t=\max\{a_t,b_t\}$, so that
$\frac {|E_t|} {2} \leq {\bf M}_t<|E_t|$.  By SP1, we have
${\bf m}_t>\frac {|E_t|} {2^{l-1} (t+1)^2}$,
so we have
\be
\frac {{\bf M}_t} {{\bf m}_{t+1}} \leq \frac {|E_t|} {{\bf m}_{t+1}}
\leq 2^{l-1} (t+2)^2
\frac {|E_t|} {|E_{t+1}|} \leq e^{c_{l-2}^{-14}} \leq
a^{-\delta}.
\ee
This allows us to estimate
\be \label {dist'}
\dist (R_{l-1}|(F_t \setminus F_{t+1})) \leq e^{c_{l-2}^{-15}} \leq
a^{-\delta}, \quad 0 \leq t<s.
\ee

If $0 \leq t \leq s-1$, let $E^1_t$, $E^2_t$ be the
connected components of $E_t \setminus E_{t+1}$, and if $t=s$, let
$E^1_t$ and $E^2_t$ be the connected components of
$I^{\d(s)}_{l-1} \setminus C^{\d(s)}_{l-1}$, (recall that
$E_s=I^{\d(s)}_{l-1}$).
Using SP1 and (\ref {est 1}) we see that
\be \label {est 1'}
\frac {|C^{\{\X_{l-1}>0\}}_{l-1} \cap E^i_t|}
{|E^i_t|}<a c_{l-1}^{-3} 2^l (t+1)^2 \leq a^{1-\delta}, \quad 0 \leq t
\leq s, \quad i=1,2.
\ee

Using (\ref {dist'}), (\ref {est 1'}) and SP1, we obtain,
\be \label {est 21}
\frac {|(F_t \setminus F_{t+1}) \cap I^{\{X_l>0\}}_l|}
{|F_t|} \leq a e^{c_{l-2}^{-16}} \ll a^{1-\delta}, \quad 0 \leq t<s.
\ee

Notice that $X_l(0)=0$, so $F_s \cap I^{\{X_l>0\}}_l=(F_s \setminus I_{l+1})
\cap I^{\{X_l>0\}}_l$.  Notice that $R_{l-1}$ takes each component of $F_s
\setminus I_{l+1}$ to either $E^1_s$ or $E^2_s$, and we have the obvious
estimate $\dist(R_{l-1}|(F_s \setminus I_{l+1}))<2^l c_{l-1}^{-1/2}$.  By
(\ref {est 1'})
\be \label {est 22}
\frac {|F_s \cap I^{\{X_l>0\}}_l|}
{|F_s|} \leq a 2^{2l} c_{l-1}^{-7/2}(s+1)^2 \leq
a e^{-c_{l-2}^{-16}} \ll a^{1-\delta},
\ee
as well.

We have
\be
I_l=F_s \cup \bigcup_{t=0}^{s-1} (F_t \setminus F_{t+1}),
\ee
so, by (\ref {est 21}) and (\ref {est 22}),
\be \label {l-1 1}
x_l \leq a e^{c_{l-2}^{-16}} \ll a^{1-\delta},
\ee
and as a consequence,
\be \label {l-1 2}
\x_l \leq a 2^l c_l^{-1} e^{c_{l-2}^{-16}} \ll a^{1-\delta} c_l^{-1}.
\ee

\subsection{The critical step}

We will consider two cases: $a<e^{-c_{l-1}^{-20}}$ (Case 1)
and $a \geq e^{-c_{l-1}^{-20}}$ (Case 2).

\subsubsection{Case 1}

The first case can be dealt by an argument which is analogous to the
analysis in \S \ref {l-1}.  We consider a sequence of nested
intervals $E_t$ in level $l$ around $R_l(0)$, and also
their pullback $F_t$ in level $l+1$ (the definitions are the same of \S \ref
{l-1} up to a shift in the indexes).
Using SP2 we analyze the density of $C^{\{X_l>0\}}_l$ in
$E_t$, which we bound by $a^{1-\delta}$ (this only works in Case 1),
and since $a^{1-\delta} \ll c_l$ (this only works in Case 1 also)
we conclude that $X_{l+1}(0)=0$.  We use SP1 to
control the pullback to $F_t$.  The reader can check the estimate
\be \label {xl+1}
x_{l+1}<a^{1-\delta}, \quad \text {in Case 1}.
\ee

\subsubsection{Case 2}

Let $r=e^{c_{l-2}^{-4}}$ and $q=(\ln a)^2$, so that $q<r$ (since we are in
Case 2).

Let $\d=(j_1,...,j_s)$ be such that $R_l(0) \in C^\d_l$.
For $1 \leq u \leq q$, we let $\d_u$ be defined by
\be
I_l^{j_u}=(R_{l-1}|I_l)^{-1}(C^{\d_u}_{l-1}).
\ee
and we let $\d_u=(j^u_1,...,j^u_{s(u)})$.

Notice that by SP2,
\begin{equation} \label {v}
s(u) \leq r_l(j_u) \leq c_{l-1}^{-14}.
\end{equation}

For $1 \leq u \leq q$ and $0 \leq v \leq s(u)$, let us define a
sequence of nested intervals $S_{u,v}$ containing $C^{\d_u}_{l-1}$ by
$S_{u,v}=I^{(j^u_1,...,j^u_v)}_{l-1}$.
Let us define nested intervals $T_{u,v}$ containing $I^{j_u}_l$ by
taking $T_{u,v}$ as the connected
component of $(R_{l-1}|I_l)^{-1}(S_{u,v})$ containing $I^{j_u}_l$.  Let
\be
W_{u,v}=
(R_l^{u-1}|I_l^{(j_1,...,j_{u-1})})^{-1}(T_{u,v}),
\ee
which is some interval containing $R_l(0)$.
Notice that if $(u_1,v_1) \leq (u_2,v_2)$ in the lexicographic order we have
$W_{u_1,v_1} \supset W_{u_2,v_2}$.

Note also that $|W_{q,s(q)}|<e^{-q} |I_l|<a^{10}$ so that
\be
\frac {|(R_l|I_{l+1})^{-1}(W_{q,s(q)})|} {|I_{l+1}|} \ll a.
\ee

By our choice of $q$, we can apply SP3 and conclude
\be \label {j^u_v}
I(j^u_v,l-1)>e^{-c_{l-2}^{-12}}>
a^{1/2}> x_{l-1}, \quad 1 \leq u \leq q, \quad v \leq s(u).
\ee
This gives
\be \label {j_u 1}
X_{l-1}(j^u_v)=0, \quad 1 \leq u \leq q, \quad v \leq s(u).
\ee
If $X_l(j_u)>0$ then there exists some $v \leq s(u)$ such that
$X_{l-1}(j^u_v)>0$
(since $X_{l-1}(0)=0$), so (\ref {j_u 1}) implies
\be \label {j_u 2}
X_l(j_u)=0, \quad 1 \leq u \leq q.
\ee
Notice that (\ref {j_u 1}), (\ref {j_u 2}) imply
\be \label {S_u,v}
S_{u,v} \cap C^{\{\X_{l-1}>0\}}_{l-1}=(R_{l-1}^v|S_{u,v})^{-1}
(C^{\{\X_{l-1}>0\}}_{l-1}),
\ee
\be \label {j_u 3}
C^{\{\X_l>0\}}_l \cap
W_{u,v}=(R_l^{u-1}|I_l^{(j_1,...,j_{u-1})})^{-1}(T_{u,v} \cap
C^{\{\X_l>0\}}_l).
\ee
In particular, (\ref {S_u,v}) gives
\be \label {S u,v}
\frac {|S_{u,v} \cap C^{\{\X_{l-1}>0\}}_{l-1}|} {|S_{u,v}|}<2^l \x_{l-1}
\leq a^{1-\delta}.
\ee

\bigskip
\paragraph{}

Let us now show that
\be \label {est'}
\frac {|T_{u,v} \cap I^{\{X_l>0\}}_l|} {|T_{u,v}|} \leq a^{1-\delta}.
\ee
There are two cases: $T_{u,v} \supset I_{l+1}$ and otherwise.

In the first case, $T_{u,v} \supset I_{l+1}$, we have
$S_{u,v}=E_v$ and $T_{u,v}=F_v$ in the notation of
\S \ref {l-1}, and (\ref {est'}) follows from (\ref {est 21}) and (\ref {est
22}).

In the second case, $T_{u,v} \not \supset I_{l+1}$, using SP1
and that $v \leq s(u) \leq c_{l-1}^{-14}$ (see (\ref {v})), we get
\be \label {Tuv}
\dist(R_{l-1}|T_{u,v})<v^2 2^{l-1} \ll c_{l-1}^{-50} \leq a^{-\delta}.
\ee
Since $X_{l-1}(0)=0$, we have
\be
T_{u,v} \cap I^{\{X_l>0\}}_l=
(R_{l-1}|T_{u,v})^{-1} (C^{\{\X_{l-1}>0\}}_{l-1} \cap S_{u,v}),
\ee
and (\ref {est'}) follows from (\ref {S u,v}) and (\ref {Tuv}).

\bigskip
\paragraph{}

Since $X_l(0)=0$ (see (\ref {X_l})),
\be
T_{u,v} \cap C^{\{\X_l>0\}}_l=T_{u,v} \cap \left (I^{\{X_l>0\}}_l
\cup \bigcup_{\ntop {I^j_l \subset T_{u,v},} {j\neq 0}}
(R_l|I^j_l)^{-1} C^{\{\X_l>0\}}_l \right ).
\ee
And by (\ref {l-1 1}), (\ref {l-1 2}) and (\ref {est'}) we get
\be \label {est 4}
\frac {|T_{u,v} \cap C^{\{\X_l>0\}}_l|} {|T_{u,v}|} \leq
a^{1-\delta} c_l^{-1}.
\ee

By (\ref {est 4}) and (\ref {j_u 3}) we get
\be
\frac {|C^{\{\X_l>0\}}_l \cap W_{u,v}|} {|W_{u,v}|}<a^{1-\delta} c_l^{-1},
\quad 1 \leq u \leq q, \quad v \leq s(u).
\ee

\bigskip
\paragraph{}

Let $Z_0 \supset Z_1 \supset ... \supset Z_t$ be an enumeration of
the $W_{u,v}$.  Let us show that
\be \label {S_k}
\frac {|Z_k|} {|Z_{k+1}|} \ll a^{-\delta}, \quad 0 \leq k \leq t-1.
\ee

Notice that if $v<s(u)$, and $T_{u,v+1}$ does not contain $0$,
\be
\frac {|T_{u,v}|} {|T_{u,v+1}|} \leq 10 \frac {|S_{u,v}|} {|S_{u,v+1}|}
\ee
(using that $R_{l-1}|I_l$ is almost purely quadratic\footnote {This follows
from the following estimate: if $0 \leq a<b<c<d$ then $2\frac
{d^2-a^2} {c^2-b^2} \geq 2 \frac {d+a} {c+b} \frac {d-a} {c-b} \geq
\frac {d-a} {c-b}$.  Thus the quadratic part of the pullback can not
decrease the relative size of $S_{u,v+1} \subset S_{u,v}$ by a factor worse
than $2$ if $0 \notin S_{u,v}$, or $4$ if $0 \in S_{u,v} \setminus
S_{u,v+1}$.}).
So in this case,
\be
\frac {|W_{u,v}|} {|W_{u,v+1}|} \leq
2^l \frac {|T_{u,v}|} {|T_{u,v+1}|} \leq 2^{2 l} \frac {|S_{u,v}|}
{|S_{u,v+1}|} \leq  2^{4l} I(j^u_{v+1},l-1)^{-1} \leq
e^{c_{l-2}^{-14}} \ll a^{-\delta}
\ee
using (\ref {j^u_v}) to estimate
$I(j^u_{v+1},l-1)$.

Let us now assume that $v<s_u$ but that
$T_{u,v+1}$ contains $0$ (in this case
$T_{u,v+1}=F_{v+1}$, $T_{u,v}=F_v$, $S_{u,v+1}=E_{v+1}$, $S_{u,v}=E_v$
in the notation of \S \ref {l-1}).  Notice that
$R_{l-1}(0) \in S_{u,v+1}$, so we can apply SP1 to see that
\be
|R_{l-1}(T_{u,v+1})| \geq \frac {|S_{u,{v+1}}|} {2^l (v+1)^2}.
\ee
Thus,
\be
\frac {|W_{u,v}|} {|W_{u,v+1}|} \leq 2^l \frac {|T_{u,v}|} {|T_{u,v+1}|} \leq
2^l \cdot 10 \left (\frac {R_{l-1}(T_{u,v})}
{R_{l-1}(T_{u,v+1})} \right )^{1/2} \leq 2^l \cdot 10 \left
(\frac {2^l (v+1)^2 |S_{u,v}|} {|S_{u,v+1}|} \right )^{1/2}.
\ee
By (\ref {v}), $v+1 \leq c_{l-1}^{-14}$, while SP2 implies
that $|S_{u,v+1}| \geq 2^{-l} e^{-c_{l-2}^{-12}} |S_{u,v}|$ (this is the
same estimate as (\ref {aa3})), so we get
\be
\frac {|W_{u,v}|} {|W_{u,v+1}|} \leq 10 \cdot 2^{2l} c_{l-1}^{-14}
e^{c_{l-2}^{-12}/2} \leq e^{c_{l-2}^{-14}} \ll a^{-\delta}.
\ee

\comm{
so using (\ref {v}) and (\ref {F}),
\be
\frac {|T_{u,v}|} {|T_{u,v+1}|} \leq 10 \cdot 2^l c_{l-1}^{-28}
\frac {|S_{u,v}|} {|S_{u,v+1}|}.
\ee
and
\be
\frac {|W_{u,v}|} {|W_{u,v+1}|} \leq
2^l \frac {|T_{u,v}|} {|T_{u,v+1}|} \leq 10 \cdot 2^{2 l} c_{l-1}^{-28}
\frac {|S_{u,v}|} {|S_{u,v+1}|} \leq e^{c_{l-2}^{-14}} \ll a^{-\delta}.
\ee
}

Consider now the case of $v=s(u)$.  Notice that
$W_{u+1,0}=(R_l^{u-1}|I^{(j_1,...,j_{u-1})}_l)^{-1}(I^{j_u}_l)$, so
\be
\frac {|W_{u,s(u)}|} {|W_{u+1,0}|} \leq 2^l \frac {|T_{u,s(u)}|}
{|I^{j_u}_l|}.
\ee
Moreover, $I^{j_u}_l$ is a connected component of
$(R_{l-1}|I_l)^{-1}(C^{\d_u}_{l-1})$ and $T_{u,s(u)}$ is a connected
component of $(R_{l-1}|I_l)^{-1}(I^{\d_u}_{l-1})$.  Since
$C^{\d_u}_{l-1}$ does not contain $R_{l-1}(0)$ (since $j_u \neq 0$), we
conclude
\be
\frac {|W_{u,s(u)}|} {|W_{u+1,0}|} \leq  2^l \frac {|T_{u,s(u)}|}
{|I^{j_u}_l|} \leq 2^{2l} \frac {|I^{\d_u}_{l-1}|}
{|C^{\d_u}_{l-1}|} \leq c_{l-1}^{-2} \ll a^{-\delta}.
\ee

\bigskip
\paragraph{}

As in \S \ref {l-1}, define $a_i,b_i>0$ by $Z_i=(R_l(0)-a_i,R_l(0)+b_i)$,
and let ${\bf m}_i=\min\{a_i,b_i\}$, ${\bf M}_i=\max\{a_i,b_i\}$.  Notice
that for $i<t$
we have $m_i>2 M_{i+1}$ (both components of $Z_i \setminus Z_{i+1}$ are
much bigger than $Z_{i+1}$).
It follows that for $i<t-1$, we have (using (\ref {S_k}))
\be
\frac {{\bf M}_i} {{\bf m}_{i+1}} \leq \frac {{\bf M}_i} {{\bf M}_{i+1}}
\frac {{\bf M}_{i+1}} {{\bf M}_{i+2}} \leq
a^{-\delta}.
\ee

Let $V_i=(R_l|I_{l+1})^{-1}(Z_i)$.  Repeating the argument used to obtain
(\ref {est 21}) and (\ref {est 22}) we get
\be
\frac {|(V_i \setminus V_{i+1}) \cap I^{\{X_{l+1}>0\}}_{l+1}|} {|V_i|} \leq
a^{1-\delta} c_l^{-1}.
\ee

On the other hand,
\be
|Z_{t-2}| \ll e^{-q} |I_l|\ll a^{10} |I_l|,
\ee
so $|V_{t-2}| \ll a |I_{l+1}|$.  Repeating the argument of (\ref {l-1 1})
and (\ref {l-1 2}) we get
\be \label {high}
x_{l+1} \leq a^{1-\delta} c_l^{-1}, \quad \text {in Case 2}.
\ee

\begin{rem}

The above estimate from above could be bigger than one if $a$ is near $c_l$. 
This means that $X_{l+1}$ could be supported on most branches of $R_{l+1}$.
\end{rem}

\subsection{Dealing with the later steps}

We are now in position to work out the later steps, aiming at the estimates
outlined in \S \ref {l+4}.  Before doing so, let us present a couple of
tools that will be repeatedly used.

\subsubsection{Useful estimates}

We will need several times the following easy estimates.

\bigskip
\paragraph{} \label {ld1}

Let $T \subset \Z$ and let $q=I(T,n)$.  Let $\hat T$ be the set of $\d$ with
at least one entry in $T$.  Then
\comm{
\be
C(\hat T,n) \leq c_n^{-1} nq.
\ee
Indeed
}
\begin{align}
C(\hat T,n) &\leq \sum_{k \geq 1} I(\{\d=(j_1,...,j_k),\, j_k \in T\},n)
\leq n^{1/2+\delta} \sum_{k \geq 1} q I(\{|\d|=k-1\},n)\\
\nonumber
&\leq n \sum_{k \geq 1} q c_n^{-1} C(\{|\d|=k-1\},n)=n c_n^{-1} q.
\end{align}

\bigskip
\paragraph{} \label {ld2}

Let $T \subset \Z$, let $q=I(T,n)$, and assume $nq<1/2$.
Assume that $k>k_0>c_{n-1}^{-2}$ satisfy $k k_0^{-1} c_n>n^2 q$. 
Let $\hat T$ be the set of $\d$ with at least $k$ entries in $T$.  Then
\be
C(\hat T,n)<e^{-k_0/8}.
\ee

To see this, let $\hat T(m)$ be the set of $\d \in \hat T$
of length $m$.  Then
\be
C(\hat T(m),n) \leq \sum_{j=k}^m \binom {m} {j} (nq)^j (1-nq)^{m-j}.
\ee

Let $m_0=c_n^{-1} k_0$.  For $m \leq m_0$,
\be
C(\hat T(m),n) \leq \sum_{j=k}^{m_0} \binom {m_0} {j} (nq)^j (1-nq)^{m_0-j}.
\ee
Notice that for $j \geq k$,
\be
\binom {m_0} {j} (nq)^j (1-nq)^{m_0-j} \leq
\binom {m_0} {k} (nq)^k (1-nq)^{m_0-k}
\ee
so for $k \leq m \leq m_0$ we have
\begin{align}
C(\hat T(m),n) &\leq m_0 \binom {m_0} {k} (nq)^k (1-nq)^{m_0-k}
\leq m_0 \left (\frac {e m_0} {k} \right )^k (nq)^{k} (1-nq)^{m_0-k}
\leq m_0 e^{-2k}\\
\nonumber
&\leq e^{-k},
\end{align}
since $(m_0 k^{-1}) nq \ll 1$.

Thus
\be
C(\hat T,n) \leq C(\{|\d|>m_0\},n)+\sum_{m \leq m_0} C(\hat T(m),n) \leq
e^{-k_0/4}+m_0 e^{-k} \leq e^{-k_0/8}.
\ee

\bigskip
\paragraph{} \label {ld3}

Let $T \subset \Z$ and let $q=I(T,n)$.  Let $q_0>n^2 q$ be such that
$q_0>c_n^{1/5}$.  Let $\hat T$ be
the set of $\d$ with $|\d|>c_n^{-1/2}$ and with at least
$q_0 |\d|$ entries in $T$.  Then
\be
C(\hat T,n) \leq e^{-c_n^{-1/4}}.
\ee
Indeed, by similar considerations as in \S \ref {ld2},
\begin{align}
C(\hat T,n) &\leq \sum_{m>c_n^{-1/2}} \sum_{j=q_0 m}^m \binom {m} {j} (nq)^j
(1-nq)^{m-j}\\
\nonumber
&\leq \sum_{m>c_n^{-1/2}} \sum_{j=q_0 m}^m \binom {m} {q_0 m}
(nq)^{q_0 m} (1-nq)^{m-q_0 m}\\
\nonumber
&\leq \sum_{m>c_n^{-1/2}} m \left (\frac {e} {q_0} \right )^{q_0 m}
(nq)^{q_0 m} (1-nq)^{m-q_0 m}\\
\nonumber
&\leq \sum_{m>c_n^{-1/2}}m e^{-q_0 m} \leq e^{-c_n^{-1/4}}.
\end{align}

\subsubsection{Standard landings}

We will need only
a couple of properties of standard landings (see \S \ref {esttime})
which we will put together here.
\begin{align}
&r_n(j_i)<c_{n-1}^{-14}, && \d=(j_1,...,j_m) \in LS(n),\\
&|\d|>c_n^{-1/2}, &&\d \in LS(n),\\
&l_n(\d)>c_{n-1}^{-1+\epsilon/2} |\d|, &&\d \in LS(n),\\
&C(\Omega \setminus LS(n),n)<c_n^{1/3}.
\end{align}

\subsection{The $l+1$ level}

Recall that $\epsilon$ was fixed in advance, so that $\delta \ll \epsilon$
if $l$ is big.

Assume that $\X_l(\d) \geq c_l^{-\epsilon}$, with $\d=(j_1,...,j_m)$. 
Then $\d \in A \cup B$, where $A$ is the set of $\d$ with
$X_l(j_i)>c_l^{-\epsilon/2}$ for some $j_i$ and $B$ is the
set of $\d$ with
$\#\{i|X_l(j_i)>0\} \geq c_l^{-\epsilon/2}$.

Since $I(\{r_l(j)>c_l^{-\epsilon/2}\},l)<e^{-c_l^{-\epsilon/3}}$, the
estimate of \S \ref {ld1} gives
\be
C(A,l) \leq e^{-c_l^{-\epsilon/4}}.
\ee

Let $k=c_l^{-\epsilon/4}$, $k_0=c_l^{-\epsilon/5}$.
Since $k k_0^{-1} c_l>l^2 c_l^{1-\delta} \geq
l^2 a^{1-\delta} \geq l^2 x_l$, the estimate of \S \ref
{ld2} gives
\be
C(B,l) \leq e^{-c_l^{-\epsilon/8}}.
\ee

We now conclude easily
\be
\x_l(c_l^{-\epsilon}) \leq C(A \cup B,l)<e^{-c_l^{-\epsilon/4}}+
e^{-c_l^{-\epsilon/8}} \leq e^{-c_l^{-\epsilon/11}} \leq c_l^{l^3}.
\ee
After pullback by
$R_l|I_{l+1}$ we get
\be \label {cll3}
x_{l+1}(c_l^{-2\epsilon}) \leq
x_{l+1}(c_l^{-\epsilon}+v_l)<e^{-c_l^{-\epsilon/12}} \leq c_l^{l^3}.
\ee

\comm{
Let $\tilde A=\{j|r_l(j) \geq c_l^{-\epsilon/2}\}$.  Then
\be
C^A_l \subset \bigcup_{\d \in \Omega} (R^\d_l)^{-1}(I^{\tilde A}_l).
\ee

On the other hand,
\be
I(\tilde A,l) \leq e^{-c_l^{-\epsilon/2}c_{l-1}^4}.
\ee
so
\be
C(A,l) \leq \sum_{r \geq 0} 2^{l+1}
\left (1-\frac {c_l} {4} \right )^r
e^{-c_l^{-\epsilon/2} c_{l-1}^4} \leq e^{-c_l^{-\epsilon/4}}.
\ee
}

\comm{
To estimate $C(B,l)$, notice that
\be
|\d|>c_l^{-\epsilon/2}, \d \in B,
\ee
\be
C(\{\d,|\d|>c_l^{-1-\epsilon/4}\},l)<
e^{-c_l^{-\epsilon/4} 2^{-l}}<
e^{-c_l^{-\delta_0/5}}.
\ee
So let us consider
\be
D=\{\d \in B,c_l^{-\epsilon/2}<|\d| \leq c_l^{-1-\epsilon/4}\}.
\ee

By (\ref {l-1} 1), $x_l<a^{1-\delta}<c_l^{1-\delta}$.  Notice that
\be
\#\{i,X_l(j_i)>0\}>c_l^{1-\epsilon/4} |\d|>c_l^{-\epsilon/5} x_l |\d|,
\quad \d=(j_1,...,j_m) \in D.
\ee
Thus a Large deviation argument gives
\be
C(B \cap D_2) \leq e^{-c_l^{-\epsilon/10}}.
\ee

Thus
\be
\x_l(c_l^{-\epsilon})<e^{-c_l^{-\epsilon/4}}+e^{-c_l^{-\epsilon/5}}+
e^{-c_l^{-\epsilon/10}} \leq e^{-c_l^{-\epsilon/11}} \leq c_l^{l^3}.
\ee
After pullback by
$R_l|I_{l+1}$ we get
\be
x_{l+1}(c_l^{-2\epsilon}) \leq
x_{l+1}(c_l^{-\epsilon}+v_l)<e^{-c_l^{-\epsilon/12}} \leq c_l^{l^3}.
\ee
}

\subsection{Levels $l+2$ and $l+3$}

We will need to consider two cases according to the size of $a$:
\begin{align}
&\text {Case A} \quad &&c_{l+1} \leq a<c_{l+1}^{\epsilon^{1/2}},\\
&\text {Case B} \quad
&&c_{l+1}^{\epsilon^{1/2}} \leq a<c_l.
\end{align}

\subsubsection{Case A}

Notice that $x_{l+1}<a^{1-\delta}$ in this case.

Let us say that $\d$ is a BAD landing (of level $l+1$) if
\be
\X_{l+1}(\d)>c_{l+1}^{-1+\epsilon}
a^{1-5 \epsilon^{1/2}}, \quad \text {BAD landing in Case A}.
\ee

\bigskip
\paragraph{}

Let us see that
\be \label {disb}
C(\{\d \text { is BAD}\},l+1) \leq e^{-c_{l+1}^{-\epsilon/8}}, \quad \text
{ in Case A}.
\ee

Indeed, if $\d$ is BAD and $\d=(j_1,...,j_m)$ then $\d \in A \cup B$ where
$A$ is the set of $\d$ with some $j_i$
with $r_{l+1}(j_i)>c_{l+1}^{-\epsilon/2}$ and $B$ is the set of $\d$
with at least $k$ entries in $\{X_{l+1}(j)>0\}$ where
$k=c_{l+1}^{-1+3\epsilon/2}a^{1-5 \epsilon^{1/2}}$.

As before, the estimate of \S \ref {ld1} gives
\be
C(A,l+1) \leq e^{-c_{l+1}^{-\epsilon/4}}.
\ee
Let $k_0=c_{l+1}^{-\epsilon}$.
Since
\be
k k_0^{-1}c_{l+1}>c_{l+1}^{3 \epsilon}
a^{1-5 \epsilon^{1/2}}>a^{1-2\epsilon^{1/2}}>(l+1)^2 x_{l+1}
\ee
the estimate of \S \ref {ld2} gives
\be
C(B,l+1) \leq e^{-c_{l+1}^{-\epsilon/4}},
\ee
and (\ref {disb}) follows.

\comm{
For the first case, the set of $\d$ such that there exists $j_i$ with
$r_n(j_i)>c_{l+1}^{-\epsilon/2}$ can be estimated by
$e^{-c_{l+1}^{-\epsilon/4}}$.

In the second case, since
$a^{1-5 \epsilon^{1/2}}>c_{l+1}^{1-5 \epsilon^{1/2}}$ we have
$|\d|>c_{l+1}^{4 \epsilon^{1/2}}$.
On the other hand,
\be
C(\{|\d| \geq c_{l+1}^{-1-\epsilon/4},l+1)<e^{-c_{l+1}^{-\epsilon/8}}.
\ee

Let
\be
D=\{\d \text {is BAD and }
c_{l+1}^{-\epsilon/4}<|\d|<c_{l+1}^{-1-\epsilon/4}\}.
\ee

We have for $\d=(j_1,...,j_m) \in D$,
\be
\#\{i,X_{l+1}(j_i)>0\}>
\frac {c_{l+1}^{-1+3\epsilon/2} a^{1-5 \epsilon^{1/2}}}
{c_{l+1}^{-1-\epsilon/4}}|\d|>c_{l+1}^{7 \epsilon/4} a^{1-5
\epsilon^{1/2}}>a^{1-3 \epsilon^{1/2}}>a^{-2 \epsilon^{1/2}}
x_l>c_{l+1}^{-2\epsilon} x_l.
\ee
Thus, the Large deviation argument gives
\be
c_{l+1}^{7 \epsilon/4} a^{1-2
\epsilon^{1/2}}>a^{1-\delta} c_{l+1}^{-\epsilon/10}
\ee

For any fixed $m<c_{l+1}^{1-\epsilon/4}$,
the probability of the union of $I^\d_{l+1}$ such that
$|\d|=m$ and the number of $1 \leq i \leq m$ such that $X_{l+1}(j_i)>0$
is at least $c_{l+1}^{-1+7\epsilon/4} a^{1-2\epsilon^{1/2}}$
can be estimated by
$e^{-c_{l+1}^{-\epsilon^{1/2}/4}}$
(since in this case they must have a frequency of $\{X_{l+1}>0\}$ bigger
than $c_{l+1}^{7\epsilon/4} a^{1-5\epsilon^{1/2}} \geq x_{l+1}
a^{-\epsilon^{1/2}}$.  Notice that $c_{l+1}^{7\epsilon/4}
a^{1-5\epsilon^{1/2}}<c_{l+1}^\epsilon$.
}

\bigskip
\paragraph{}

Define the set of BAD returns (of level $l+2$)
as the set of $j$ such that $R_{l+1}(I^j_{l+2}) \subset C^\d_{l+1}$ where
$\d$ is a BAD landing, so that
$I(\{j \text { is BAD}\},l+2)<e^{-c_{l+1}^{-\epsilon/11}}$.
Notice that a non BAD return $j$ satisfies
\be
X_{l+2}(j) \leq c_{l+1}^{-1+\epsilon} a^{1-5 \epsilon^{1/2}}+v_{l+1}, \quad
j \text { is not a BAD return}.
\ee

Let us define a set VB of landings (of level $l+2$) as the set of
$\d=(j_1,...,j_m)$ which are either
non-standard or $\#\{i,j_i \text { is a BAD return}\}>c_{l+1}^l|\d|$.
Notice that $c_{l+1}^l>(l+2)^2 I(\{j \text { is BAD}\}, l+2)$, so
using the estimate of \S \ref {ld3} (with $q_0=c_{l+1}^l>c_{l+2}^{1/5}$),
we conclude that
\be
C(VB,l+2) \leq c_{l+2}^{2/7}.
\ee

If $\d \notin VB$, all returns have time at most
$r_n(j_i)<c_{l+1}^{-14}$, so
\begin{align}
\X_{l+2}(\d) &\leq (c_{l+1}^{-14} c_{l+1}^l+c_{l+1}^{-1+\epsilon}
a^{1-5\epsilon^{1/2}}+v_{l+1}) |\d|\\
\nonumber
&\leq (c_{l+1}^{l-14}+c_{l+1}^{-1+\epsilon}
a^{1-5 \epsilon^{1/2}}+c_l^{-1-\delta}) c_{l+1}^{1-\epsilon/2} l_{l+2}(\d)
\leq a^{1-5 \epsilon^{1/2}} l_{l+2}(\d)
\end{align}
since $l_{l+2}(\d)>c_{l+1}^{-1+\epsilon/2} |\d|$ for a standard landing.

\bigskip
\paragraph{}

Define a VB return as the set of returns of level $l+3$ that fall in VB
landings.  Then $I(VB,l+3)<c_{l+2}^{1/20}$.
Each non VB return satisfies (notice that
$r_{l+3}(j) \geq c^{-1/2}_{l+2}$ when $j$ is not a VB return)
\be
X_{l+3}(j) \leq a^{1-5 \epsilon^{1/2}} r_{l+3}(j)+v_{l+2} \leq
a^{1-7 \epsilon^{1/2}} r_{l+3}(j).
\ee
By the estimate of \S \ref {l+4}, we conclude that
$\mu(I^\j_\n) \leq a^{1-10 \epsilon^{1/2}}$ in Case A.

\comm
{
By the Inductive Estimate, Lemma \ref {ier}, the time contribution
of non VB returns in a VG
return is less then $c_{l+3}^{1/200}$.  We conclude
\be
\X_n(j)<a^{1-2 \sqrt \epsilon} r_n(j) \quad \text {if } j \in VG(VB,l+3,n).
\ee

Since the majority of branches are VG(VB,l+3,n), we can use \S \ref {comp}
to conclude the main estimate in this case.
}

\subsubsection{Case B}

Let GOOD denote the set of $\d \in LS(l+1)$ such that
$\d=(j_1,...,j_m)$, and
\begin{align}
\label {GOOD1}
&\#\{i,\, X_{l+1}(j_i) \geq 1\}<a^{1-\epsilon} c_l^{-1} |\d|,
&&\text {GOOD landing, Case B},\\
\label {GOOD2}
&\#\{i,\, X_{l+1}(j_i) \geq c_l^{-3 \epsilon^{1/2}}\}<c_l^{l^3} |\d|,
&&\text {GOOD landing, Case B}.
\end{align}

By (\ref {xl+1}), (\ref {high}), and (\ref {cll3}) we have
that $a^{1-\epsilon} c_l^{-1} \geq (l+1)^2 x_{l+1}$ and
$c_l^{l^3} \geq (l+1)^2 x_{l+1}(c_l^{-3 \epsilon^{1/2}})$.

\bigskip
\paragraph{}

Let $D_1$ (respectively $D_2$)
be the set of $\d$ such that $|\d|>c_{l+1}^{-1/2}$ and which do not
satisfy (\ref {GOOD1}) (respectively (\ref {GOOD2})).
The argument of \S \ref {ld3} with $q_0=a^{1-\epsilon}c_l^{-1}$
(respectively, $q_0=c_l^{l^3}$) implies that (notice that in both cases
$q_0>c_{l+1}^{1/5}$)
\be
C(D_1,l+1),C(D_2,l+1)<e^{-c_{l+1}^{-1/4}}.
\ee
\comm{
\be
The measure of the union of $C^\d_{l+1}$ with $\d$ standard and
$\#\{i|X_{l+1}(j_i) \geq c_l^{-2\epsilon^{1/2}}\}>c_l^{l^3}$
can be bounded by
\be
\sum_{c_{l+1}^{-1/2}<m<c_{l+1}^{-2}} e^{-c_l^{l^3} m} <
e^{-c_{l+1}^{-1/3}}.
\ee

The measure of the union of $C^\d_{l+1}$ with $\d$ is standard and
$\#\{i|X_{l+1}(j_i) \geq 1\}>a^{1-50\epsilon} c_l^{-1}$
can be bounded by
\be
\sum_{c_{l+1}^{-1/2}<m<c_{l+1}^{-2}} e^{a^{1-50 \epsilon} m} <
e^{-c_{l+1}^{-1/3}}.
\ee
}
We conclude
\be
C(\{\d \text { is not GOOD}\},l+1) \leq c_{l+1}^{2/7}.
\ee

\bigskip
\paragraph{}

Let us say that a return $j$
of level $l+2$ is BAD if
$R_{l+1}(I^j_{l+2}) \subset C^\d_{l+1}$ where $\d$ is not GOOD.  Notice that
$I(\{j \text { is BAD}\},l+2)<c_{l+1}^{1/20}$.

If j is not BAD, with $R_{l+1}(I^j_{l+2}) \subset C^\d_{l+1}$, let us
consider two subcases.
If $c_{l+1}^{\epsilon^{1/2}} \leq a<c_l^{\sqrt l}$, using (\ref {GOOD1}) we
get
\be
X_{l+2}(j) \leq
(c_l^{-14} c_l^{-1} a^{1-\epsilon}) |\d|+v_{l+1} \leq a^{1-2 \epsilon} |\d|
\leq a^{1-2\epsilon} r_{l+2}(j).
\ee
If $c_l^{\sqrt l} \geq a>c_l$, a similar estimate can be
obtained using (\ref {GOOD1}) and (\ref {GOOD2})
\begin{align}
X_{l+2}(j) &\leq
(c_l^{l^3} c_l^{-14}+c_l^{-3 \epsilon^{1/2}} c_l^{-1} a^{1-\epsilon})
|\d|+v_{l+1}
\leq (c_l^{l^3-14}+c_l^{-1-3 \epsilon^{1/2}}
a^{1-\epsilon}+c_{l+1}^{1/2-\delta})|\d|\\
\nonumber
&\leq (c_l^{-1-5 \epsilon^{1/2}}
a^{1-\epsilon}) |\d| \leq a^{1-7\epsilon^{1/2}} l_{l+1}(\d)
\leq a^{1-7 \epsilon^{1/2}} r_{l+2}(j).
\end{align}
\comm{
If $c_l^{1+5 \sqrt \epsilon}>a>c_l$, (\ref {GOOD2}) gives
\be
X_{l+2}(j) \leq
(c_l^{l^3} c_l^{-14}+c_l^{-2\sqrt \epsilon})
|\d|+v_{l+1} \leq c_l^{1-3 \sqrt \epsilon} r_{l+2}(j)
\leq a^{1-10\sqrt \epsilon} r_{l+2}(j).
\ee
}

Thus we have in both subcases
\be
X_{l+2}(j) \leq a^{1-7 \epsilon^{1/2}} r_{l+2}(j), \quad
j \text { is not BAD, Case B}.
\ee

By the argument of \S \ref {l+4}, we conclude that
$\mu(I^\j_\n) \leq a^{1-10 \epsilon^{1/2}}$ in Case B.

This concludes the proof of the main estimate.

\comm{
For $n \geq l+2$ and for $j \in VG(BAD,l+2,n)$ we have
\be
X_n(j) \leq (c_{l+2}^{1/200}+a^{1-\sqrt \epsilon}) r_n(j) \leq
a^{1-2\sqrt \epsilon} r_n(j)
\ee
Since $\lim I(VG(BAD,l+2,n),n)=1$, the argument of \S \ref {comp} implies
that
\be
\mu(I^\j_\n) \leq a^{1-2 \sqrt \epsilon}
\ee
and we conclude the main estimate in Case B.
}

\comm{
Now $VG(BAD,l+2,n)$ are the majority of branches and any such branch has
$X_n(j)<a^{1-52 \epsilon}$ by our choice of $a$.  We conclude the main
estimate in Case B.

\subsubsection{Case B}

Let GOOD denote the set of $\d \in LS(l+1)$ such that
$\d=(j_1,...,j_m)$,
\be
\#\{i|X_{l+1}(j_i) \geq 1\}<a^{1-50 \epsilon} c_l^{-1} m \quad \text {GOOD
landing, Case B}.
\ee

The measure of the union of $C^\d_{l+1}$ with $\d$ is standard and
$\#\{i|X_{l+1}(j_i) \geq 1\}>a^{1-50\epsilon} c_l^{-1}$
can be bounded by
\be
\sum_{c_{l+1}^{-1/2}<m<c_{l+1}^{-2}} (1/2)^{a^{1-51\epsilon} m} <
e^{-c_{l+1}^{-1/3}}.
\ee
We conclude
\be
C(\{\d \text { is not GOOD}\},l+1) \leq c_{l+1}^{2/7}.
\ee

\bigskip
\paragraph{}

Let us say that a return of level $l+1$ is BAD if it falls in a not GOOD
landing.

If j is not BAD (falling in $\d$) then
\be
X_{l+2}(j) \leq
(c_l^{-14} c_l^{-1} a^{1-50\epsilon}) |\d|+v_{l+1} \leq a^{1-51\epsilon}
r_{l+2}(j).
\ee

Now $VG(BAD,l+2,n)$ are the majority of branches and any such branch has
$X_n(j)<a^{1-52 \epsilon}$ by our choice of $a$.  We conclude the main
estimate in Case B.

\subsubsection{Case 3}

Let GOOD denote the set of $\d \in LS(l+1)$ such that
$\d=(j_1,...,j_m)$,
\be
\#\{i|X_{l+1}(j_i) \geq 1\}<a^{1-50 \epsilon} c_l^{-1} m \quad \text {GOOD
landing Case 3}
\ee
(notice that
$a^{1-50 \epsilon} c_l^{-1} m \leq c_l^{\epsilon^{1/2}} m$
in Case 3) and
\be
\#\{i|X_{l+1}(j_i) \geq c_l^{-2 \sqrt\epsilon}\}<c_l^{l^3} m \quad \text {GOOD
landing Case 3}.
\ee

\bigskip
\paragraph{}

The measure of the union of $C^\d_{l+1}$ with $\d$ standard and
$\#\{i|X_{l+1}(j_i) \geq c_l^{-2\sqrt\epsilon}\}>c_l^{l^3}$
can be bounded by
\be
\sum_{c_{l+1}^{-1/2}<m<c_{l+1}^{-2}} e^{-c_l^{l^3} m} <
e^{-c_{l+1}^{-1/3}}.
\ee

The measure of the union of $C^\d_{l+1}$ with $\d$ is standard and
$\#\{i|X_{l+1}(j_i) \geq 1\}>a^{1-50\epsilon} c_l^{-1}$
can be bounded by
\be
\sum_{c_{l+1}^{-1/2}<m<c_{l+1}^{-2}} e^{a^{1-50 \epsilon} m} <
e^{-c_{l+1}^{-1/3}}.
\ee

We conclude
\be
C(\{\d \text { is not GOOD}\},l+1) \leq c_{l+1}^{2/7}.
\ee

\bigskip
\paragraph{}

Let us say that a return of level $l+1$ is BAD if it falls in a not GOOD
landing.

If $j$ is not BAD (falling in $\d$) then
\be
\X_{l+1}(\d) \leq
(c_l^{l^3} c_l^{-14}+c_l^{-2\sqrt \epsilon} c_l^{-1} a^{1-50\epsilon})
|\d|+v_{l+1}<a^{1-3 \sqrt \epsilon} r_{l+2}(j).
\ee

This implies the main estimate in Case 3.

\subsubsection{Case 4}

Let GOOD denote the set of $\d \in LS(l+1)$ such that
$\d=(j_1,...,j_m)$ and
\be
\#\{i|X_{l+1}(j_i) \geq c_l^{-2 \sqrt\epsilon}\}<c_l^{l^3} m \quad \text
{GOOD landing Case 4}.
\ee

The measure of the union of $C^\d_{l+1}$ with $\d$ standard and
$\#\{i|X_{l+1}(j_i) \geq c_l^{-2\sqrt\epsilon}\}>c_l^{l^3}$
can be bounded by
\be
\sum_{c_{l+1}^{-1/2}<m<c_{l+1}^{-2}} e^{-c_l^{l^3} m} <
e^{-c_{l+1}^{-1/3}}.
\ee

We conclude
\be
C(\{\d \text { is not GOOD}\},l+1) \leq c_{l+1}^{2/7}.
\ee

\bigskip
\paragraph{}

Let us say that a return of level $l+1$ is BAD if it falls in a not GOOD
landing.

If $j$ is not BAD (falling in $\d$) then
\be
X_{l+1}(j) \leq
(c_l^{l^3} c_l^{-14}+c_l^{-2\sqrt \epsilon})
|\d|+v_{l+1} \leq c_l^{1-3 \sqrt \epsilon} r_{l+2}(j)
\leq a^{1-10\sqrt \epsilon} r_{l+2}(j).
\ee

We conclude the main estimate in Case 4.
}

\section{Pathological laminations and other consequences}
\label {consequences}

\subsection{Laminations in spaces of analytic unimodal maps}

\setcounter{paragraph}{0}

\subsubsection{}

Let $\F$ be a Banach space.
A codimension-one {\it holomorphic lamination} $\LL$ on an open subset
$\WW \subset \F$ is a family of disjoint codimension-one  
Banach submanifolds of $\F$, called the {\it leaves} of the lamination
such that for any point $p \in \WW$, there exists a holomorphic local chart
$\Phi:\tilde \WW \to \VV \oplus \C$,
where $\tilde \WW \subset \WW$ is a neighborhood of $p$ and
$\VV$ is an open set in some complex Banach space $\E$,
such that for any leaf $L$ and any connected component $L_0$ of $L\cap \WW$,
the image $\Phi(L_0)$  is a graph of a holomorphic function $\VV \ra \C$.

The local theory of codimension-one holomorphic laminations coincide with
the theory of holomorphic motions (see \cite {ALM}, \S 2.5 and references
therein).  It follows from the $\lambda$-Lemma that holonomy maps of
holomorphic laminations have quasiconformal extensions.

\subsubsection{}

For $a>0$, let $\Omega_a \subset \C$ be the set of $z$ at distance at most
$a$ of $I$.  Let $\EE_a$ be the space of even holomorphic maps $f:\Omega_a
\to \C$, continuous up to $\overline \Omega_a$.  We endow $\EE_a$ with the
sup norm.  Let $\AAA_a=\{f \in \EE_a, f(-1)=f(1)=1\}$ and let $\AAA^\R_a=\{f
\in \AAA_a,\, f(z)=\overline {f(\overline z)}\}$.
Let $\UU_a$ be the space of
analytic quasiquadratic maps which belong to $\AAA^\R_a$.

One of the main results of \cite {ALM} is that the partition of $\UU_a$ on
topological conjugacy classes has the structure of a codimension-one
analytic lamination ``almost everywhere''.

\begin{thm}[Theorem A of \cite {ALM}] \label {8.1}

Let $f \in \U_a$ be a Kupka-Smale quasiquadratic map.  There exists a
neighborhood $\VV \subset \AAA_a$ of $f$ endowed with a codimension-one
holomorphic lamination $\LL$ (also called hybrid lamination)
with the following properties:

(1)\, the lamination is real-symmetric;

(2)\, if $g\in \VV\cap \AAA_a^\R$ is non-regular,
then the intersection of the leaf through 
$g$ with $\AAA_a^\R$ coincides with the intersection of the topological
conjugacy class of $g$ with $\VV$;

(3)\, Each $g \in \VV \cap \AAA^\R_a$ belongs to some leaf of $\LL$.

\end{thm}

(See also \cite {Av} for the non-quasiquadratic case.)

Notice that the set of non-Kupka-Smale maps is
contained on a countable union of codimension-one analytic submanifolds.

The lamination $\LL$ has automatically quasisymmetric holonomy. 
Quasisymmetric maps are not always absolutely continuous (even though
quasiconformal maps are).  It turns out that $\LL$ is very far from being
absolutely continuous, at least at the set of non-regular leaves (the
lamination restricted to regular maps is not uniquely defined, but can be
chosen in a quite natural way to be locally analytic).

\subsubsection{}

Let $\hat \LL$ be the lamination consisting of the non-isolated
non-regular leaves of $\LL$.  If $f_\lambda$ is a one-dimensional
family transversal to $\hat \LL$, it intersects $\hat \LL$ in a positive
measure set\footnote {Since preperiodic combinatorics are dense in $\hat
\LL$, and the generic unfolding of preperiodic combinatorics generates a
positive measure set of non-hyperbolic parameters (see for instance \cite
{T2} for a proof of a more general statement).}.

Let $X \subset \UU_a$ be the set of Collet-Eckmann
maps satisfying the conclusion of Theorem A.  Then $X$ intersects each
leaf of $\hat \LL$ in a set of maps which are analytically conjugate in the
attractor, and this
is a set of infinite codimension (possibly
empty): just notice that we can vary the exponent of any finite number
of periodic orbits independently.

Thus, $\hat \LL$ exhibits the same pathology described
by Milnor in \cite {Mi}: a full measure set intersecting the leaves
of a finite codimension lamination in tiny sets.
(The example described by Milnor is also
an analytic codimension-one lamination, on two dimensions, and the
intersections of the leaves with the full measure set are points.  In this
finite dimensional setting, this translates in the
complete failure of Fubini's Theorem).

\subsubsection{}

Although in our description we have to make use of transverse measures
(since our setting is infinite-dimensional), one can interpret this
pathology by taking finite dimensional sections as follows.

Let $\{f_\lambda\}_{\lambda \in \Lambda}$
be a small analytic $k$-dimensional transverse section to $\hat \LL$.
The lamination $\hat \LL$ induces a lamination $\hat \LL_\Lambda$ on
$\Lambda$.  Notice that $\hat \LL_\Lambda$ has positive $k$-dimensional
Lebesgue measure.

For $\lambda_1, \lambda_2 \in \Lambda$ distinct,
let $P(\lambda_1,\lambda_2)$ be the
number of periodic orbits of $f_{\lambda_1}$ which have the same exponent of
a periodic orbit of $f_{\lambda_2}$.
A transversality argument shows that, for most of those sections (actually
the complement has infinite-codimension), $P(\lambda_1,\lambda_2)<\infty$
whenever $\lambda_1 \neq \lambda_2$.  For any family $f_\lambda$ with this
property, we obtain the same phenomena, but in $k$ dimensions:
the set of parameters $X_\Lambda$ which are Collet-Eckmann and
satisfy the conclusion of Theorem A intersects each leaf of $\hat
\LL_\Lambda$ in at most one point.

\subsubsection{}
We point out that the set of {\it recurrent}
parameters which do not satisfy the conclusion
of Theorem A (or Theorem B) has Hausdorff dimension one in any
one-dimensional transversal $f_\lambda$ to $\hat \LL$.
Indeed, using the previous argument,
we can select another transversal $\tilde f_\lambda$
arbitrarily close to $f_\lambda$, and such that for any $\lambda_1$,
$\lambda_2$, the number of periodic orbits of $f_{\lambda_1}$ which have the
same exponent as some periodic orbit of $\tilde f_{\lambda_2}$ is finite.
Let $\tilde X$ be the set of parameters satisfying the conclusion of
Theorems A and B for $\tilde f_\lambda$.
Let $h$ be the holonomy
map from $f_\lambda$ to $\tilde f_\lambda$.
Then the quasisymmetric constant of $h$ (and thus the H\"older constant)
is close to $1$, provided $\tilde f_\lambda$ is close to
$f_\lambda$ (by the $\lambda$-Lemma).  Since $\tilde X$ has
positive Lebesgue measure, $h^{-1}(\tilde X)$ has Hausdorff
dimension close to one.  But parameters $\lambda$
in $h^{-1}(\tilde X)$ do not satisfy the conclusion of Theorem A (or Theorem
B, using \S \ref {geom}):
otherwise each of the infinitely many periodic orbits in the attractor
of $f_\lambda$ would have the same exponent of the corresponding
(by the topological conjugacy) periodic orbit for $\tilde f_{h(\lambda)}$.

\begin{rem}

It is easy to see that the conclusions of
Theorems A, B, and C fail for all Misiurewicz
(non-recurrent Collet-Eckmann) parameters.  This set of parameters has
Hausdorff dimension one in any one-dimensional transversal to $\hat \LL$.

\end{rem}

\comm{
\subsubsection{On universality and the holonomy method} \label {universa}

The key result from complex dynamics needed to obtain measure-theoretical
results for non-trivial analytic families of unimodal maps is the
the lamination structure for the partition in topological classes
of \cite {ALM}.  Since the holonomy of this lamination is
quasisymmetric, the parameter space of the quadratic family
has a universal quasisymmetric structure.

This was used by \cite {ALM} as a
shortcut to extend the regular or stochastic dichotomy, and in
\cite {AM2}, the same shortcut was used to
deduce the Collet-Eckmann and polynomial recurrence conditions from
statements for the quadratic family.

This {\it holonomy method}, consisting in the comparison between
parameter spaces of different families has to be applied to estimates which
on one hand are topological invariants.  More seriously, the set of
combinatorics concerned must have full measure {\it simultaneously} in all
non-trivial families of unimodal maps.

As we just saw, the quasisymmetric holonomy of the lamination is not
absolutely continuous.
It is clear that we cannot use the holonomy method to prove
lack of absolute continuity of the lamination itself, and in particular
Theorem A.

We should point out that the Phase-Parameter relation of \cite {Av} used in
this work is still based on the lamination structure.  However, one uses
local holonomy maps to obtain relate phase and parameter of the same family,
instead of using the global holonomy to relate phase and parameter between
some non-trivial family and the quadratic family, which introduces serious
distortion and lack of sharpness in the estimates.
}


\subsubsection{} \label {geom}

By \cite {MM}, stochastic
unimodal maps satisfying the conclusion of Theorem A are geometrically
rigid: two such maps are smoothly (and automatically
analytically) conjugate on the attractor.
The same conclusion can be obtained for
maps satisfying the conclusion of Theorem B.  Indeed, the asymptotic
distribution of the critical orbit, if it exists, is a topological
invariant, and Theorem B implies that
the conjugacy must be absolutely continuous on the support of the invariant
measure (the attractor).  The conjugacy is then easily promoted to being
smooth (and automatically analytic) by a well known argument,
see Exercise 3.1, Chapter V, page 375 in \cite {MS}.

We should point out that this conjugacy is not, in
general, analytic on the whole interval $I$,
as can be shown by a simple example.

\begin{example}

Let us consider the families
\be
f_a(x)=a x (1-x),
\ee
\be
g_a(x)=\frac {2} {\pi} \sin^{-1} \left ( \frac {\sqrt a} {2} \sin (\pi x)
\right ).
\ee

Then for $f_a$ and $g_a$ are analytic families of quasiquadratic
maps on the interval $[0,1]$ for $2<a<4$.  Notice that
\be
h_a \circ g_a=f_a \circ h_a, \quad h(x)=\frac {1-\cos (\pi x)} {2},
\ee
so $f_a$ and $g_a$ are analytically conjugate on $(0,1)$, and the holonomy
map between both families is trivial.  Whenever $f_a$
is Collet-Eckmann and satisfies the conclusion of Theorem A, $g_a$ also
does.  However, $f_a$ and $g_a$ are not analytically conjugate on $[0,1]$:
indeed, $Df_a(0)=a$ and $Dg_a(0)=\sqrt a$, so the exponent of the
fixed point $0$ is not preserved.

\end{example}

\subsection{Formula for the exponent of $\mu_f$} \label {lyapformula}

In order to compute the Lyapunov exponent of $\mu_f$ combinatorially, one
just has to find an expression for the Lyapunov exponent of the critical
value.  There are several ways to proceed, for instance, one can find
convenient approximations of the critical orbit by periodic orbits and apply
Theorem A.  However, there exists a very simple expression using the
combinatorics of the principal nest, whose proof only involves Corollary
\ref {lyap} and the asymptotic limits of \S \ref {values}.

\begin{thm}

Let $f_t$ be an analytic family of unimodal maps.  For almost every
non-regular parameter, the Lyapunov exponent of $\mu_{f_t}$ (which is
equal to the Lyapunov exponent of the critical value) is given
by an explicit combinatorial formula:
\be
\lambda(\mu_{f_t})=\lim_{n \to \infty} \frac {2 \ln v_{n+1}} {v_n},
\ee
where, as usual, $v_n$ is the return time of the critical point to the
$n$-th level of the principal nest.

\end{thm}

\begin{pf}

Notice that
\be
\lambda(f(0))=\lim \frac {\ln |Df^{v_n-1}(f(0))|} {v_n-1},
\ee
and by Lemma \ref {vn+1},
\be
\lim \frac {\ln v_{n+1}} {\ln c_n^{-1}}=1.
\ee
Thus, we only have to show that
\be
\lim \frac {\ln |Df^{v_n-1}(f(0))|} {\ln c_n^{-1}}=2.
\ee
Notice that
$f^{v_n-1}$ takes $f(I_{n+1})$ to $R_n(I_{n+1})$ with torrentially
small distortion.  By Lemma \ref {n1+delta},
we have $n^{-2}|I_n| \leq |R_n(I_{n+1})| \leq |I_n|$.
So we conclude
\be
\lim \frac {\ln |Df^{v_n-1}(f(0))|} {\ln c_n^{-1}}=
\lim \ln \left ( \frac {|I_n|} {|I_{n+1}|^2} \right )
\frac {1} {\ln c_n^{-1}}=
\lim \ln \left ( \frac {1} {|I_n| c_n^2} \right )
\frac {1} {\ln c_n^{-1}}=2,
\ee
since $|I_n|>c_{n-1}^{1+\delta}$ for $n$ big.
\end{pf}

\comm{
\appendix

\section{The complex phase-parameter relation}

\input{complex1.tex}

\input{lamination.tex}
}


\end{document}